\documentclass[12pt]{amsart}
\usepackage{amsmath}
\usepackage{amssymb}
\usepackage{latexsym}

\newcommand{\Zint}{\mathbb {Z}}    
     
\newcommand{\Rea}{\mathbb {R}}      
\newcommand{\Cplx}{\mathbb {C}}     

\newcommand{\halmos}{\rule{5pt}{5pt}}

\numberwithin{equation}{section}

\newtheorem{prop}{\bf Proposition}[section]

\newtheorem{thm}[prop]{\bf Theorem}
\newtheorem{lemma}[prop]{\bf Lemma}

\newtheorem{exa}{\bf Example}

\begin{document}

\title[Hermite-Krichever ansatz for Fuchsian equations]
{The Hermite-Krichever ansatz for Fuchsian equations with applications to the sixth Painlev\'e equation and to finite-gap potentials}
\author{Kouichi Takemura}
\address{Department of Mathematical Sciences, Yokohama City University, 22-2 Seto, Kanazawa-ku, Yokohama 236-0027, Japan.}
\email{takemura@yokohama-cu.ac.jp}


\begin{abstract}
Several results including integral representation of solutions and Hermite-Krichever Ansatz on Heun's equation are generalized to a certain class of Fuchsian differential equations, and they are applied to equations which are related with physics.
We investigate linear differential equations that produce Painlev\'e equation by monodromy preserving deformation and obtain solutions of the sixth Painlev\'e equation which include Hitchin's solution.
The relationship with finite-gap potential is also discussed.
We find new finite-gap potentials.
Namely, we show that the potential which is written as the sum of the Treibich-Verdier potential and additional apparent singularities of exponents $-1$ and $2$ is finite-gap, which extends the result obtained previously by Treibich.
We also investigate the eigenfunctions and their monodromy of the Schr\"odinger operator on our potential.
\end{abstract}

\maketitle

\section{Introduction}

It is well known that a Fuchsian differential equation with three singularities is transformed to a Gauss hypergeometric equation, and plays important roles in substantial fields in mathematics and physics. Several properties of solutions to the hypergeometric equation have been explained in various textbooks.

A canonical form of a Fuchsian equation with four singularities is written as 
\begin{equation}
\left( \! \left(\frac{d}{dw}\right) ^2 \! + \left( \frac{\gamma}{w}+\frac{\delta }{w-1}+\frac{\epsilon}{w-t}\right) \frac{d}{dw} +\frac{\alpha \beta w -q}{w(w-1)(w-t)} \right)\tilde{f}(w)=0
\label{Heun}
\end{equation}
with the condition 
\begin{equation}
\gamma +\delta +\epsilon =\alpha +\beta +1,
\label{Heuncond}
\end{equation}
and is called Heun's equation.
Heun's equation frequently appears in physics, i.e., general relativity \cite{STU}, fluid mechanics \cite{CH} and so on.
Despite that Heun's equation was resolved in the 19th century; several results of solutions have only been recently revealed.
Namely, integral representations of solutions, global monodromy in terms of hyperelliptic integrals, relationships with the theory of finite-gap potential and the Hermite-Krichever Ansatz for the case $\gamma, \delta, \epsilon , \alpha -\beta \in \Zint +1/2$ are contemporary (see \cite{BE,GW,Smi,Tak1,Tak2,Tak3,Tak4,TV} etc.), though they are not written in a textbook on Heun's equation \cite{Ron}.

In this paper, we consider differential equations which have additional apparent singularities to Heun's equation.
More precisely, we consider the equation
\begin{align}
& \left\{ \frac{d^2}{dw^2}+\left( \frac{\frac{1}{2}-l_1}{w}+  \frac{\frac{1}{2}-l_2}{w-1}+  \frac{\frac{1}{2}-l_3}{w-t}+ \sum _{{i'}=1}^M \frac{-r_{i'}}{w-\tilde{b}_{i'}} \right) \frac{d}{dw} \right. \label{Feqintro} \\
& \left.  +\frac{(\sum _{i=0}^3 l_i + \sum _{{i'}=1}^M r_{i'})(-1-l_0 +\sum _{i=1}^3 l_i + \sum _{{i'}=1}^M r_{i'})w+\tilde{p}+ \sum_{{i'}=1}^M \frac{\tilde{o}_{i'}}{w-\tilde{b}_{i'}}}{4w(w-1)(w-t)}\right\}\tilde{f}(w) =0 , \nonumber
\end{align}
for the case $l_i \in \Zint _{\geq 0}$ $(0 \leq i \leq 3)$, $r_{i'} \in \Zint _{>0}$ $(1 \leq i' \leq M)$ and the regular singular points $\tilde{b}_{i'}$ $(1 \leq i' \leq M)$ are apparent.

By a certain transformation, Eq.(\ref{Feqintro}) is rewritten in terms of elliptic functions such as
\begin{align}
 & \left\{ -\frac{d^2}{dx^2} + \sum_{i=0}^3 l_i(l_i+1)\wp (x+\omega_i) \right. \label{ellDE} \\
 & \left. + \sum_{i'=1}^M \left( \frac{r_{i'}}{2} \left( \frac{r_{i'}}{2}+1\right) (\wp (x-\delta _{i'}) + \wp (x+\delta _{i'})) +\frac{s_{i'}}{\wp (x) -\wp (\delta _{i'})} \right) -E \right\} f(x) =0, \nonumber
\end{align}
with the condition that logarithmic solutions around the singularities $x= \pm \delta _{i'}$ $(i'=1,\dots ,M)$ disappear.
We then establish that solutions to Eq.(\ref{ellDE}) have an integral representation and they are also written as a form of the Hermite-Krichever Ansatz. For details see Proposition \ref{prop:Linteg} and Theorem \ref{thm:alpha}.
Note that the results on the Hermite-Krichever Ansatz are related to Picard's theorem on differential equations with coefficients of elliptic functions \cite[\S 15.6]{Inc}.
By the Hermite-Krichever Ansatz, we can obtain information on the monodromy of solutions to differential equations.

Results on the integral representation and the Hermite-Krichever Ansatz are applied for particular cases.
One example is Painlev\'e equation.
For the case $M=1$ and $r_1=1$, it is known that Eq.(\ref{Feqintro}) produces the sixth Painlev\'e equation by monodromy preserving deformation (see \cite{IKSY}).
On the other hand, solutions to Eq.(\ref{ellDE}) are expressed as  a form of the Hermite-Krichever Ansatz for the case $l_i \in \Zint _{\geq 0}$ $(i=0,1,2,3)$, and we obtain an expression of monodromy.
Fixing monodromy corresponds to the monodromy preserving deformation; thus, we obtain solutions to the sixth Painlev\'e equation by fixing monodromy (see section \ref{sec:P6}). For the case $l_0=l_1=l_2=l_3=0$, we recover Hitchin's solution \cite{Hit}.
Note that the sixth Painlev\'e equation and the Hitchin's solution appear in topological field theory \cite{Man} and Einstein metrics \cite{Hit}.

Another example for application of the integral representation and the Hermite-Krichever Ansatz is finite-gap potential.
On solid-state physics, band structure of spectral is essential, and examples and properties of finite-gap (finite-band) potential could be applicable (e.g. see \cite{BEES}).

Recently several authors have been active in producing a variety of studies of finite-gap potential, and several results have been applied to the analysis of Schr\"odinger-type operators and so on.
Here we briefly review these results. 
Let $q(x)$ be a periodic, smooth, real function, $H$ be the operator $-d^2/dx^2+q(x)$, and $\sigma _b(H)$ be the set such that 
$$
E \in \sigma _b(H) \; \Leftrightarrow \mbox{ Every solution to }(H-E)f(x)=0 \mbox{ is bounded on }x \in \Rea .
$$
If the closure of the set $\sigma _b(H)$ can be written as
\begin{equation}
\overline{\sigma _b(H)}= [E_{0},E_{1}]\cup [E_{2},E_{3}]  \cup \dots \cup [E_{2g}, \infty ),
\end{equation}
where $E_0<E_{1}<\cdots <E_{2g}$, then $q(x)$ is called the finite-gap potential.

Let $\wp (x)$ be the Weierstrass $\wp $-function with periods $(2\omega _1, 2\omega _3)$. Ince \cite{I} established in 1940 that if $n \in \Zint _{\geq 1}$, $\omega _1 \in \Rea$ and $\omega _3 \in \sqrt{-1} \Rea$, then the potential of the Lam\'e's operator, 
\begin{equation}
-\frac{d^2}{dx^2}+n(n+1)\wp (x),
\end{equation}
is finite-gap.
From the 1960s, relationships among finite-gap potentials, odd-order commuting operators and soliton equations were investigated. 
If there exists an odd-order differential operator 
$A= \left( d/dx \right)^{2g+1} +  \! $ $ \sum_{j=0}^{2g-1}\! $ $ b_j(x) \left( d/dx \right)^{2g-1-j} $ such that $[A, -d^2/dx^2+q(x)]=0$, then $q(x)$ is called the algebro-geometric finite-gap potential.
Under the condition that $q(x) $ is real-valued, smooth and periodic, it is known that $q(x)$ is a finite-gap potential if and only if $q(x)$ is an algebro-geometric finite-gap potential. For a detailed historical review, see \cite{GW2} and the references therein.

In the late 1980s, Treibich and Verdier invented the theory of elliptic solitons, which is based on an algebro-geometric approach to soliton equations developed by Krichever \cite{Kri} among others, and found a new algebro-geometric finite-gap potential, which is now called the Treibich-Verdier potential (see \cite{TV}). This potential may be written in the form
\begin{equation}
v(x)=  \sum_{i=0}^3 l_i(l_i+1)\wp (x+\omega_i)
\label{TVpotent}
\end{equation}
for the Schr\"odinger operator $-d^2/dx^2+v(x)$, where $l_i$ $(i=0,1,2,3)$ are integers and $\omega _1$, $\omega _3$, $\omega _0 (=0)$, $\omega _2 (=\omega _1 -\omega _3)$ are half-periods.
Subsequently several studies \cite{GW,Wei,Smi,Tak1,Tak2,Tak3,Tak4} have further added to understanding of this subject.
Note that the function in Eq.(\ref{TVpotent}) corresponds to the potential of the Schr\"odinger operator as Eq.(\ref{ellDE}) for the case $M=0$, and it is closely related to Heun's equation.

Later, by following his joint work with Verdier, Treibich \cite{Tre} established that, if $l_0 , l_1, l_2 , l_3 \in \Zint _{\geq 0}$ and $\delta $ satisfy
\begin{equation}
\sum _{i=0}^3 (l_i +1/2)^2 \wp ' (\delta  +\omega _i ) =0,
\label{eq:Trepot}
\end{equation}
then the potential 
\begin{align}
v(x) = & 2 (\wp (x-\delta ) + \wp (x+\delta )) + \sum_{i=0}^3 l_i(l_i+1) \wp (x+\omega_i) 
\end{align}
for the Schr\"odinger operator $-d^2/dx^2+v(x)$ is algebro-geometric finite-gap. In \cite{Smi2} Smirnov presented  further results.

In this paper, we generalize the results of Treibich and Smirnov. In particular, we will find that, if $l_0 , l_1, l_2 , l_3 \in \Zint _{\geq 0}$, $\delta _j \not \equiv \omega _i$ mod $2\omega_1 \Zint \oplus 2\omega_3 \Zint$ $(0\leq i\leq 3, \; 1\leq j\leq M)$ and  $\delta _j \pm \delta _{j'} \not \equiv 0$  mod $2\omega_1 \Zint \oplus 2\omega_3 \Zint$ $(1\leq j< j' \leq M)$, and $\delta _1 ,\dots, \delta _M$ satisfy the equation
\begin{equation}
2\sum _{j' \neq j} (\wp ' (\delta _j -\delta _{j'} ) + \wp ' (\delta _j +\delta _{j'} ) )+ \sum _{i=0}^3 (l_i +1/2)^2 \wp ' (\delta _j +\omega _i ) =0 \quad (j=1,\dots ,M),
\label{eq:dsi}
\end{equation}
then the potential
\begin{align}
v(x) = & \sum_{i=0}^3 l_i(l_i+1) \wp (x+\omega_i) +2\sum_{i'=1}^M (\wp (x-\delta _{i'}) + \wp (x+\delta _{i'})) \label{fingapDEin}
\end{align}
for the Schr\"odinger operator $-d^2/dx^2+v(x)$ is algebro-geometric finite-gap.
Note that the potential in Eq.(\ref{fingapDEin}) corresponds to Eq.(\ref{ellDE}) with conditions $r_{i'}=2$ and $s_{i'}=0$ $(i'=1,\dots ,M)$.
For the special case $M=1$, we recover the Treibich's result \cite{Tre}.

Our approach differs from that of Treibich and Verdier and is elementary; we do not use knowledge of sophisticated algebraic geometry. The approach is based on writing the product of two specific eigenfunctions of the Schr\"odinger operator in the form of a doubly-periodic function for all eigenvalues $E$, which follows from the apparency of regular singularities of the Schr\"odinger operator (see section \ref{sec:HK}).
Using the doubly-periodic function, an odd-order commuting operator is constructed, and it follows that the potential is algebro-geometric finite-gap.
As a consequence, we obtain results concerning integral representations of solutions, monodromy formulae in terms of a hyperelliptic integral, the Bethe-Ansatz and the Hermite-Krichever Ansatz, as is shown in \cite{Tak3,Tak4} for Heun's equation.
We can also obtain two expression of monodromy.
By comparing the two expressions, we obtain hyperelliptic-to-elliptic integral reduction formulae.
Note that our approach can be related to the theory of Picard's potential, which is developed by Gesztesy and Weikard \cite{GW2}.

This paper is organized as follows. 
In section \ref{sec:HK}, we obtain integral representations of solutions to the differential equation of the class mentioned above and rewrite them to the form of the Hermite-Krichever Ansatz. To obtain an integral representation, we introduce doubly-periodic functions that satisfy a differential equation of order three. Some properties related with this doubly-periodic function are investigated, and we obtain another expression of solutions that looks like the form of the Bethe Ansatz.
In section \ref{sec:P6}, we consider the relationship with the sixth Painlev\'e equation. We show that solutions of the sixth Painlev\'e equation are obtained from solutions expressed in the form of the Hermite-Krichever Ansatz of linear differential equations considered in section \ref{sec:HK} by fixing monodromy. Some explicit solutions that include Hitchin's solution are displayed.
In section \ref{sec:FGP}, we discuss the relationship with the results on finite-gap potential.
In subsection \ref{sec:fg}, we show that the potential $v(x)$ in Eq.(\ref{fingapDEin}) is algebro-geometric finite-gap under the conditions of Eq.(\ref{eq:dsi}).
In subsection \ref{sec:hypell}, we express global monodromy of eigenfunctions of the Schr\"odinger operator in terms of a hyperelliptic integral.
In subsection \ref{sec:HKA}, we investigate the eigenfunctions and monodromy by the Bethe Ansatz and the Hermite-Krichever Ansatz. As a consequence, we are able to derive another monodromy formula. 
In subsection \ref{sec:red}, we obtain hyperelliptic-to-elliptic integral reduction formulae by comparing two expressions of monodromy.
In section \ref{sec:exa}, we consider several examples on finite-gap potential.
In section \ref{sec:rmk}, we give concluding remarks and present an open problem.
In the appendix, we note definitions and formulae for elliptic functions.

\section{Fuchsian differential equation and Hermite-Krichever Ansatz} \label{sec:HK}

\subsection{Fuchsian differential equation} \label{sec:FDE}
To begin with, we introduce the following differential equation;
\begin{align}
& \left\{ \frac{d^2}{dz^2}+\left( \sum _{i=1}^3 \frac{\frac{1}{2}-l_i}{z-e_i}+ \sum _{{i'}=1}^M \frac{-r_{i'}}{z-b_{i'}} \right) \frac{d}{dz} +\frac{N(N-2l_0-1)z+p+ \sum_{{i'}=1}^M \frac{o_{i'}}{z-b_{i'}}}{4(z-e_1)(z-e_2)(z-e_3)}\right\}\tilde{f}(z) =0 ,
\label{Feq} 
\end{align}
where $N=\sum _{i=0}^3 l_i + \sum _{{i'}=1}^M r_{i'}$. This equation is Fuchsian, i.e., all singularities $\{ e_i \} _{i=1,2,3 }$, $\{ b_{i'} \} _{ i' =1, \dots ,M}$ and $\infty$ are regular. The exponents at $z=e_i$ $(i=1,2,3)$ (resp. $z=b_{i'} $ $(i'=1,\dots ,M)$) are $0$ and $l_i +1/2$ (resp. $0$ and $r_{i'} +1$), and the exponents at $z= \infty$ are $N/2$ and $(N-2l_0-1)/2$. Conversely, any Fuchsian differential equation that has regular singularities at $\{ e_i \} _{i=1,2,3 }$, $\{ b_{i'} \} _{ i' =1, \dots ,M}$ and $\infty$ such that one of the exponents at $e_i$ and $b_{i'} $ for all $i\in \{1,2,3 \}$ and $i' \in \{1, \dots ,M\}$ are zero is written as Eq.(\ref{Feq}).
By the transformation $z \rightarrow z +\alpha $, we can change to the case $e_1 +e_2 +e_3 =0$. In this paper we restrict discussion to the case $e_1 +e_2 +e_3 =0$.
We remark that any Fuchsian equation with $M+4$ singularities is transformed to Eq.(\ref{Feq}) with the condition $e_1 +e_2 +e_3 =0$.

It is known that, if $e_1 +e_2 +e_3 =0$ and $e_1 \neq e_2 \neq e_3 \neq e_1$, then there exists some periods $(2\omega _1 ,2\omega _3)$ such that $\wp (\omega _1)= e_1$ and $\wp (\omega _3 )= e_3$, where $\wp (x)$ is the Weierstrass $\wp$-function with periods $(2\omega_1, 2\omega_3)$.
We set $\omega _0 =0$ and $\omega _2 =-\omega _1 -\omega _3$. Then we have $\wp (\omega _2 )= e_2$.

Now we rewrite Eq.(\ref{Feq}) in an elliptic form. We set
\begin{equation}
\Phi (z)= \prod_{i=1}^3 (z-e_i)^{-l_i/2} \prod_{{i'}=1}^M (z-b_{i'})^{-r_{i'}/2} , \quad z= \wp(x),
\end{equation}
and $\tilde{f}(z) \Phi (z)=f(x)$. Then we have
\begin{equation}
(H-E) f(x)= 0,
\label{eq:H}
\end{equation}
where $H$ is a differential operator defined by
\begin{align}
 H= & -\frac{d^2}{dx^2} + v(x), \label{Ino} \\
v(x) = & \sum_{i=0}^3 l_i(l_i+1)\wp (x+\omega_i) \label{Inopotent} \\
& + \sum_{i'=1}^M \frac{r_{i'}}{2} \left( \frac{r_{i'}}{2}+1\right) (\wp (x-\delta _{i'}) + \wp (x+\delta _{i'})) +\frac{s_{i'}}{\wp (x) -\wp (\delta _{i'})}, \nonumber
\end{align}
and
\begin{align}
& \wp (\delta _{i'})= b_{i'} , \quad ({i'}=1,\dots ,M) ,\\
& o_{i'} = -s_{i'} + r_{i'} \left\{ \frac{1}{8}r _{i'} (12b_{i'} ^2-g_2)   +\frac{1}{2}(4b_{i'} ^3- g_2 b_{i'} -g_3)\left(\sum_{i'' \neq {i'}}\frac{r_{i''}}{(b_{i'} -b_{i''})}\right) \right. \\
& \quad \quad \left. +2(l_1(b_{i'} -e_2)(b_{i'} -e_3)+l_2(b_{i'} -e_1)(b_{i'} -e_3)+l_3(b_{i'} -e_1)(b_{i'} -e_2))\right\} ,\nonumber \\
& p= E+(e_1l_1^2+e_2l_2^2+e_3l_3^2)  -2(l_1l_2e_3+l_2l_3e_1+l_3l_1e_2) -\frac{1}{2}\sum _{{i'}=1}^M b_{i'} r_{i'} ^2\\
&  \quad \quad  +2\sum _{{i'}=1}^M \sum _{i=1}^3 l_ir_{i'} (e_i+b_{i'} )+2\left(\sum _{{i'}=1}^M b_{i'} r_{i'} \right)\left(\sum _{{i'}=1}^M r_{i'} \right) ,\nonumber \\
& g_2=-4(e_1e_2+e_2e_3+e_3e_1), \quad g_3=4e_1e_2e_3. 
\end{align}
Conversely, Eq.(\ref{Feq}) is obtained from Eq.(\ref{eq:H}) by the transformation above.

We consider another expression. Set
\begin{align}
&  H_g= -\frac{d^2}{dx^2} + \sum_{i'=1}^M  \frac{r_{i'} \wp ' (x)}{\wp (x) -\wp (\delta _{i'})} \frac{d}{dx} + \left(l_0 + \sum_{i'=1}^M r_{i'}\right) \left(l_0 +1-  \sum_{i'=1}^M r_{i'}\right) \wp (x) \\
& \quad \quad +\sum_{i=1}^3 l_i(l_i+1) \wp (x+\omega_i) + \sum_{i'=1}^M \frac{\tilde{s}_{i'}}{\wp (x) -\wp (\delta _{i'})} ,\nonumber \\
& f_g(x) = f(x) \Psi _g (x), \quad \Psi _g (x)=\prod _{i'=1}^M (\wp (x) -\wp (\delta _{i'}))^{r_{i'} /2}. 
\end{align}
Then Eq.(\ref{eq:H}) is also equivalent to 
\begin{equation}
(H_g -E-C_g) f_g(x) =0,
\label{eq:Hg}
\end{equation}
where 
\begin{align}
& \tilde{s}_{i'} = s_{i'} -r_{i'} \left\{ \frac{1}{8} r _{i'} (12 b_{i'} ^2-g_2) +\frac{1}{2}(4b_{i'} ^3- g_2 b_{i'} -g_3)\left(\sum_{i'' \neq {i'}}\frac{r_{i''}}{(b_{i'} -b_{i''})}\right) \right\} ,\\
& C_g =- \frac{1}{2}\sum _{{i'}=1}^M b_{i'} r_{i'} ^2 +2\left(\sum _{{i'}=1}^M b_{i'} r_{i'} \right)\left(\sum _{{i'}=1}^M r_{i'} \right) .
\end{align}
Note that the exponents at $x=\pm \delta _{i'} $ $(i'=1,\dots ,M)$ are $0$ and $r_{i'} +1$.

In this paper, we consider solutions to Eq.(\ref{Feq}), which is equivalent to Eq.(\ref{eq:H}) or Eq.(\ref{eq:Hg}) for the case $l_i \in \Zint $, and the regular singular point $z= b_{i'}$ is apparent for all ${i'}$. Here, a regular singular point $x=a$ of a linear differential equation of order two is said to be apparent, if and only if the differential equation does not have a logarithmic solution at $x=a$ and the exponents at $x=a$ are integers. It is known that the regular singular point $x=a$ is apparent, if and only if the monodromy matrix around $x=a$ is a unit matrix.
Note that Smirnov investigated solutions to Eq.(\ref{eq:H}) in \cite{Smi2} with the assumptions $s_{i'}=0$, $r_{i'} \in 2\Zint $ for all ${i'}$.

Now we study the condition that the regular singular point $x=a$ is apparent.
More precisely, we describe the condition that a differential equation of order two does not have logarithmic solutions at a regular singular point $x=a$ $(a \neq \infty)$ for the case $\alpha _2 -\alpha _1 \in \Zint$, where $\alpha _1$ and  $\alpha _2$ are exponents at $x=a$.
If $\alpha _1 = \alpha _2$, then the differential equation has logarithmic solutions at $x=a$. We assume that the exponents satisfy $\alpha _2 -\alpha _1 =n \in \Zint _{\geq 1}$.
Since the point $x=a$ is a regular singular, the differential equation is written as
\begin{equation}
\left\{ \frac{d^2}{dx^2} + \sum _{j=0}^{\infty} p_j (x-a)^{j-1} \frac{d}{dx} + \sum _{j=0}^{\infty} q_j (x-a)^{j-2} \right\} f(x)=0, \label{eq:Feqxa}
\end{equation}
for some $p_j , q_j \in \Cplx$ $(j=0,1, \dots )$.
Let $F(t)$ be the characteristic polynomial at the regular singular point $x=a$.
Since exponents at $x=a$ are $\alpha _1$ and $\alpha _2$, $F(t)$ is written as $F(t)= t^2+(p_0-1)t+q_0=(t-\alpha _1)(t-\alpha _2)$.
We now calculate solutions to Eq.(\ref{eq:Feqxa}) in the form 
\begin{equation}
f(x) = \sum _{j=0}^{\infty} c_j (x-a)^{\alpha _1 + j},
\end{equation}
where $f(x)$ is normalized to satisfy $c_0=1$. By substituting it into Eq.(\ref{eq:Feqxa}) and comparing the coefficients of $(x-a)^{\alpha _1 + j-2}$, we obtain the relations 
\begin{equation}
F(\alpha _1 +j) c_j + \sum_{j'=0}^{j-1} \{ (\alpha _1 +j') p_{j-j'} +q_{j-j'}\} c_{j'} =0. \label{eq:recj}
\end{equation}
If the positive integer $j$ satisfies $F(\alpha _1 +j) \neq 0$ (i.e. $j\neq 0,n$), then the coefficient $c_j$ is determined recursively.
For the case $j=n$, we have $F(\alpha _1 +n)=0$ and 
\begin{equation}
\sum_{j'=0}^{n-1} \{ (\alpha _1 +j') p_{n-j'} +q_{n-j'}\} c_{j'} =0. \label{eq:recn}
\end{equation}
Eq.(\ref{eq:recn}) with recursive relations (\ref{eq:recj}) for $j=1, \dots ,n-1$ is a necessary and sufficient condition that Eq.(\ref{eq:Feqxa}) does not have a logarithmic solution for the case $\alpha _2 -\alpha _1 =n \in \Zint _{\geq 1}$. In fact, if $p_0 ,q_0 , \dots , p_n ,q_n$ satisfy Eq.(\ref{eq:recn}), then there exist solutions to Eq.(\ref{eq:Feqxa}) that include two parameters $c_0$ and $c_n$. Thus any solutions are not logarithmic at $x=a$. Conversely, if Eq.(\ref{eq:recn}) is not satisfied, there exists a logarithmic solution written as $f(x)= \sum _{j=0}^{\infty} c_j (x-a)^{\alpha _1 + j} + \log (x-a) \sum _{j=n}^{\infty} \tilde{c}_j (x-a)^{\alpha _1 + j}$.

It follows from $\wp (\delta _{i'})= b_{i'}$, $\wp '(\delta _{i'}) \neq 0$ and holomorphy of $\prod _{i=1}^3 (\wp (x)-e_i) ^{-l_i/2}$ at $x=\pm \delta _{i'}$ that, the monodromy matrix to Eq.(\ref{Feq}) around a regular singular point $z= b_{i'}$ is a unit matrix, if and only if the monodromy matrix to Eq.(\ref{eq:Hg}) around a regular singular point $x= \pm \delta _{i'}$ is a unit matrix.
It is obvious that, if the monodromy matrix to Eq.(\ref{Feq}) around a regular singular point $z= b_{i'}$ is a unit matrix, then we have $r_{i'} \in \Zint _{\neq 0}$. In this paper we assume that $r_{i'} \in \Zint _{>0}$ for all $i'$.

\subsection{Integral representation and the Hermite-Krichever Ansatz}

We introduce doubly-periodic functions to obtain an integral expression of solutions to Eq.(\ref{eq:H}) (or Eq.(\ref{eq:Hg})) for the case $l_i \in \Zint _{\geq 0}$ $(i=0,1,2,3)$, $r_{i'} \in \Zint _{> 0}$ $(i'=1,\dots ,M)$ and the regular singular points $z=b_{i'}$ $(i'=1,\dots ,M)$ of Eq.(\ref{Feq}) are apparent.

\begin{prop} \label{prop:prod}
Let $v(x)$ be the function defined in Eq.(\ref{Inopotent}).
If $l_i \in \Zint _{\geq 0}$ $(i=0,1,2,3)$, $r_{i'} \in \Zint _{> 0}$ $(i'=1,\dots ,M)$ and regular singular points $z=b_{i'}$ $(i'=1,\dots ,M)$ of Eq.(\ref{Feq}) are apparent, then the equation
\begin{align}
& \left\{ \frac{d^3}{dx^3}-4\left( v(x) -E\right)\frac{d}{dx} -2\frac{dv(x)}{dx} \right\} \Xi (x)=0,
\label{prodDE}
\end{align}
has an even nonzero doubly-periodic solution that has the expansion
\begin{equation}
\Xi (x)=c_0+\sum_{i=0}^3 \sum_{j=0}^{l_i-1} b^{(i)}_j \wp (x+\omega_i)^{l_i-j} + \sum _{i'=1}^M  \sum_{j=0}^{r_{i'}-1} \frac{d^{(i')}_j}{(\wp (x)-\wp (\delta _{i'}))^{r_{i'} -j}}.
\label{Fx}
\end{equation}
\end{prop}
\begin{proof}
First, we show a lemma that is related to the monodromy of solutions to Eq. (\ref{eq:Hg}).
\begin{lemma} \label{prop:locmonod}
If $l_0, l_1, l_2, l_3  \in \Zint _{\geq 0}$, then the monodromy matrix of Eq.(\ref{eq:Hg}) around a point $x=n_1 \omega_1 + n_3 \omega_3$ $(n_1, n_3 \in \Zint)$ is a unit matrix.
\end{lemma} 
\begin{proof}
Due to periodicity, it is sufficient to consider the case $x=\omega _i$ $(i=0,1,2,3)$. We first deal with the case $i=1,2,3$.
The exponents at the singular point $x=\omega_i$ ($i=1,2,3$) are $-l_i$ and $l_{i}+1$. Because Eq.(\ref{eq:Hg}) is invariant under the transformation $x -\omega _i \rightarrow -(x-\omega _i) $ and the gap of the exponents at $x=\omega_i$ (i.e. $l_i+1-(-l_i)$) is odd, there exist solutions in the form $f_{i,1}(x)=(x-\omega _i)^{-l_i}(1+\sum_{j=1}^{\infty}a_j (x-\omega _i)^{2j})$ and $f_{i,2}(x)=(x-\omega _i)^{l_i+1}(1+\sum_{j=1}^{\infty}a'_j (x-\omega _i)^{2j})$.
Since the functions $f_{i,1}(x)$ and $f_{i,2}(x)$ form a basis for solutions to Eq.(\ref{eq:Hg}) and they are non-branching around the point $x=\omega _i$, the monodromy matrix around $x=\omega _i$ is a unit matrix.
For the case $i=0$, the exponents at $x=0$ are $-l_0 - \sum_{i'=1}^M r_{i'}$ and $l_{0}+1 - \sum_{i'=1}^M r_{i'}$, the gap of the exponents is odd, and similarly it is shown that the monodromy matrix around the point $x=0$ is a unit matrix.
Hence we obtain the lemma.
\end{proof}
We continue the proof of Proposition \ref{prop:prod}.
Let $M_j$ $(j=1,3)$ be the transformations
obtained by the analytic continuation $x \rightarrow x+2\omega _j$.
It follows from double-periodicity of Eq.(\ref{eq:Hg}) that, if $f_g(x)$ is a solutions to Eq.(\ref{eq:Hg}), then $M_j f_g(x)$ $(j=1,3)$ is also a solution to Eq.(\ref{eq:Hg}).
From the assumption that regular singular points $z=b_{i'}$ are apparent for all ${i'}$, the monodromy matrix to Eq.(\ref{eq:Hg}) around a regular singular point $x= \pm \delta _{i'}$ is a unit matrix for all $i'$.
By combining with Lemma \ref{prop:locmonod}, it follows that all local monodromy matrices around any singular points are units. Hence the transformations $M_j$ do not depend on the choice of paths. From the fact that the fundamental group of the torus is commutative, we have $M_1 M_3=M_3 M_1$.
Recall that the operators $M_j$ act on the space of solutions to Eq.(\ref{eq:Hg}) for each $E$, which is two dimensional.
By the commutativity $M_1 M_3=M_3 M_1$, there exists a joint eigenvector $\tilde{\Lambda }_g (x)$ for the operators $M_1$ and $M_3$. 
It follows from Proposition \ref{prop:locmonod} and the apparency of singular points that the function $\tilde{\Lambda } _g (x)$ is single-valued and satisfies equations $(H_g-E-C_g) \tilde{\Lambda } _g(x)=0$, $M_1\tilde{\Lambda } _g(x)=\tilde{m}_1\tilde{\Lambda } _g(x)$ and $M_3\tilde{\Lambda } _g(x)=\tilde{m}_3\tilde{\Lambda } _g(x)$ for some $\tilde{m}_1,\tilde{m}_3 \in \Cplx \setminus \{0\}$.
By changing parity $x \leftrightarrow -x$, it follows immediately that $(H_g-E-C_g) \tilde{\Lambda } _g(-x)=0$, $M_1\tilde{\Lambda } _g(-x)=\tilde{m}_1^{-1}\tilde{\Lambda } _g(-x)$ and $M_3\tilde{\Lambda } _g(-x)=\tilde{m}_3^{-1}\tilde{\Lambda } _g(-x)$. Then the function $\tilde{\Lambda } _g(x)\tilde{\Lambda } _g(-x)$ is single-valued, even and doubly-periodic.
We set $\tilde{\Lambda } (x)= \tilde{\Lambda } _g(x) /\Psi _g (x) $. Then $\tilde{\Lambda } (x)$ and $\tilde{\Lambda } (-x)$ are solutions to Eq.(\ref{eq:H}).

Now consider the function $\Xi (x)=\tilde{\Lambda } _g(x)\tilde{\Lambda } _g(-x) /\Psi _g (x)^2$. Since the function $\Psi _g (x)^2$ is single-valued, even and doubly-periodic, the function $\Xi (x)$ is single-valued, even (i.e. $\Xi (x)=\Xi (-x)$), doubly-periodic (i.e. $\Xi (x+2\omega _1)=\Xi (x+2\omega _3)=\Xi (x)$), and satisfies the equation
\begin{align}
& \left\{ \frac{d^3}{dx^3}-4\left( v(x) -E\right)\frac{d}{dx} -2\frac{dv(x)}{dx} \right\} \Xi (x)=0
\nonumber
\end{align}
that the products of any pair of solutions to Eq.(\ref{eq:H}) satisfy.

Since the function $\Xi (x)$ is an even doubly-periodic function that satisfies the differential equation (\ref{prodDE}) and the exponents of Eq.(\ref{prodDE}) at $x=\omega _i$ $(i=0,\dots ,3)$ (resp. $x=\pm \delta _{i'}$ $(i'=1,\dots ,M)$) are $-2l_i,1,2l_i+2$ (resp. $-r _{i'} ,1,r _{i'} +2 $), it is written as a rational function of variable $\wp (x)$, and it admits the expansion as Eq.(\ref{Fx}) by considering exponents.
\end{proof}

The function $\Xi (x)$ is calculated by substituting Eq.(\ref{Fx}) into the differential equation (\ref{prodDE}) and solving simultaneous equations for the coefficients. We introduce an integral formula for a solution to the differential equation Eq.(\ref{eq:H}) in use of the function $\Xi (x)$.
Set
\begin{align}
& Q= \Xi (x)^2\left( E- v(x)\right) +\frac{1}{2}\Xi (x)\frac{d^2\Xi (x)}{dx^2}-\frac{1}{4}\left(\frac{d\Xi (x)}{dx} \right)^2. \label{const}
\end{align}
It follows from Eq.(\ref{prodDE}) that
\begin{align}
\frac{dQ}{dx} = \frac{1}{2} \Xi (x) \left( 4\frac{d\Xi (x)}{dx} (E-v(x) )- 2 \Xi (x) \frac{dv(x)}{dx} + \frac{d^3\Xi (x)}{dx^3} \right) =0 .
\label{pfconst}
\end{align}
Hence the value $Q$ is independent of $x$.
\begin{prop} \label{prop:Linteg}
Let $\Xi (x)$ be the doubly-periodic function defined in Proposition \ref{prop:prod} and $Q$ be the value defined in Eq.(\ref{const}).
Then the function 
\begin{equation}
\Lambda ( x)=\sqrt{\Xi (x)}\exp \int \frac{ \sqrt{-Q}dx}{\Xi (x)},
\label{integ1}
\end{equation}
is a solution to the differential equation (\ref{eq:H}), and the function
\begin{equation}
\Lambda _g ( x)=\Psi _g (x) \sqrt{\Xi (x)}\exp \int \frac{ \sqrt{-Q}dx}{\Xi (x)},
\label{integ1g}
\end{equation}
is a solution to the differential equation (\ref{eq:Hg}).
\end{prop}
\begin{proof}
From Eqs.(\ref{integ1}, \ref{const}) we have
\begin{align}
& \frac{\Lambda ' (x)}{\Lambda (x)} =\frac{1}{2} \frac{\Xi ' (x)}{\Xi (x)} + \frac{\sqrt{-Q}}{\Xi (x)} , \label{eq:llp} \\
& \frac{\Lambda '' (x)}{\Lambda (x)} = \frac{1}{2} \frac{\Xi '' (x)}{\Xi (x)} - \frac{1}{4} \left( \frac{\Xi ' (x)}{\Xi (x)} \right)^2 - \frac{Q}{\Xi (x)^2} = v(x) -E.
\end{align}
Hence we have $-\frac{d^2}{dx^2} \Lambda (x) +v(x)\Lambda (x) =E\Lambda (x)$.
It follows from the equivalence of Eq.(\ref{eq:H}) and Eq.(\ref{eq:Hg}) that the function $\Lambda _g ( x)$ is a solution to Eq.(\ref{eq:Hg}). 
\end{proof}

\begin{prop} \label{prop:indep}
If $Q\neq 0$, then the functions $\Lambda (x) $ and $\Lambda (-x) $ are linearly independent and any solution to Eq.(\ref{eq:H}) is written as a linear combination of $\Lambda (x) $ and $\Lambda (-x) $.
\end{prop}
\begin{proof}
It follows from Eq.(\ref{eq:llp}) and the evenness of the function $\Xi (x)$ that 
\begin{equation}
\frac{\frac{d}{dx} \Lambda (-x)}{\Lambda (-x) }  = \frac{1}{2} \frac{\Xi ' (x)}{\Xi (x)} - \frac{\sqrt{-Q}}{\Xi (x)}.
\end{equation}
Hence we have
\begin{equation}
\Lambda (-x) \frac{d}{dx} \Lambda (x) -\Lambda (x) \frac{d}{dx} \Lambda (-x) = \Lambda (x) \Lambda (-x) \frac{2\sqrt{-Q}}{\Xi (x)}.
\label{eq:L-L}
\end{equation}
If $\Lambda (x) $ and $\Lambda (-x) $ are linearly dependent, then the l.h.s. of Eq.(\ref{eq:L-L}) must be zero; however, this is impossible because $Q \neq 0$. Hence the functions $\Lambda (x)$ and $\Lambda (-x)$ are linearly independent.
It follows from the invariance of Eq.(\ref{eq:H}) with respect to the transformation $x \leftrightarrow -x$ that $\Lambda (-x) $ is also a solution to Eq.(\ref{eq:H}).

Since solutions to Eq.(\ref{eq:H}) form a two-dimensional vector space and the functions $\Lambda (x)$ and $\Lambda (-x)$ are linearly independent, the functions $\Lambda (x)$ and $\Lambda (-x)$ form a basis of the space of solutions to Eq.(\ref{eq:H}), and any solution to Eq.(\ref{eq:H}) is written as a linear combination of $\Lambda (x) $ and $\Lambda (-x) $.
\end{proof}

It follows from Proposition \ref{prop:indep} that, if $Q\neq 0$, then the functions $\Lambda _g(x) $ and $\Lambda _g(-x) $ are linearly independent, and any solution to Eq.(\ref{eq:Hg}) is written as a linear combination of $\Lambda _g(x) $ and $\Lambda _g(-x) $.

From the formulae (\ref{integ1}, \ref{integ1g}) and the doubly-periodicity of the functions $\Xi (x)$ and $\Psi_g (x)^2$, we have 
\begin{align}
& \Lambda (x+2\omega _j)=\pm \Lambda (x) \exp \int _{0+ \varepsilon }^{2\omega _j+\varepsilon }\frac{\sqrt{-Q}dx}{\Xi (x)}, \quad (j=1,3), \label{eqn:Lam01} \\
& \Lambda _g (x+2\omega _j)=\pm \Lambda _g (x) \exp \int _{0+ \varepsilon }^{2\omega _j+\varepsilon }\frac{\sqrt{-Q}dx}{\Xi (x)}, \quad (j=1,3), \label{eqn:Lam01g}
\end{align}
with $\varepsilon $ a constant determined so as to avoid passing through the poles while integrating. The sign $\pm$ is determined by the analytic continuation of the function $\sqrt{\Xi (x)}$, and the integrations in Eqs.(\ref{eqn:Lam01}, \ref{eqn:Lam01g}) may depend on the choice of the path.
The function $\Lambda (x)$ may have branching points, althought the function $\Lambda _g (x)$ does not have branching points and is meromorphic on the complex plane, because $\Lambda _g (x)$ is a solution to Eq.(\ref{eq:Hg}) and any singularity of Eq.(\ref{eq:Hg}) is apparent.
It follows from Eq.(\ref{eqn:Lam01g}) that there exists $m_1, m_3 \in \Cplx $ such that 
\begin{equation}
\Lambda _g (x+2\omega _j) =\exp (\pi \sqrt{-1} m_j ) \Lambda _g (x), \quad (j=1,3).
\label{eqn:Lam010}
\end{equation}

We now show that a solution to Eq.(\ref{eq:Hg}) can be expressed in the form of the Hermite-Krichever Ansatz.
We set
\begin{equation}
\Phi _i(x,\alpha )= \frac{\sigma (x+\omega _i -\alpha ) }{ \sigma (x+\omega _i )} \exp (\zeta( \alpha )x), \quad \quad (i=0,1,2,3),
\label{Phii}
\end{equation}
where $\sigma (x)$ (resp. $\zeta (x)$) is the Weierstrass sigma (resp. zeta) function.
Then we have 
\begin{equation}
\left( \frac{d}{dx} \right) ^{k} \Phi _i(x+2\omega _{j} , \alpha ) = \exp (-2\eta _{j} \alpha +2\omega _{j} \zeta (\alpha )) \left( \frac{d}{dx} \right) ^{k} \Phi _i(x, \alpha )
\label{ddxPhiperiod}
\end{equation}
for $i=0,1,2,3$, $j=1,2,3$ and $k \in \Zint _{\geq 0}$, where $\eta _j =\zeta (\omega _j)$ $(j=1,2,3)$.
\begin{thm} \label{thm:alpha}
Set $\tilde{l} _0 = l_0 +\sum _{i'=1}^M r_{i'}$ and $\tilde{l}_i =l_i$ $(i=1,2,3)$.
The function $\Lambda _g (x)$ in Eq.(\ref{integ1g}) is expressed as
\begin{align}
& \Lambda _g (x) = \exp \left( \kappa x \right) \left( \sum _{i=0}^3 \sum_{j=0}^{\tilde{l}_i-1} \tilde{b} ^{(i)}_j \left( \frac{d}{dx} \right) ^{j} \Phi _i(x, \alpha ) \right)
\label{Lalpha}
\end{align}
for some values $\alpha $, $\kappa$ and $\tilde{b} ^{(i)}_j$ $(i=0,\dots ,3, \: j= 0,\dots ,\tilde{l}_i-1)$, or
\begin{align}
& \Lambda _g (x) = \exp \left( \bar{\kappa } x \right) \left( \bar{c} +\sum _{i=0}^3 \sum_{j=0}^{\tilde{l}_i-2} \bar{b} ^{(i)}_j \left( \frac{d}{dx} \right) ^{j} \wp (x+\omega _i) +\sum_{i=1}^3 \bar{c}_i \frac{\wp '(x)}{\wp (x)-e_i} \right)
\label{Lalpha0}
\end{align}
for some values $\bar{\kappa }$, $\bar{c}$, $\bar{c}_i$ $(i=1,2,3)$ and $\bar{b} ^{(i)}_j$ $(i=0,\dots ,3, \: j= 0,\dots ,\tilde{l}_i-2)$.

If the function $\Lambda _g (x)$ is expressed as Eq.(\ref{Lalpha}), then
\begin{align}
& \Lambda _g (x+2\omega _j) = \exp (-2\eta _j \alpha +2\omega _j \zeta (\alpha ) +2 \kappa \omega _j ) \Lambda _g (x) , \quad  (j=1,3), \label{ellint} 
\end{align}
else
\begin{align} 
& \Lambda _g (x+2\omega _j) = \exp (2 \bar{\kappa } \omega _j ) \Lambda _g (x) , \quad  (j=1,3). \label{ellint0} 
\end{align}		
\end{thm}

\begin{proof}
Set
\begin{align}
& \alpha  = -m_1 \omega _3 +m_3 \omega _1 \label{al},
\end{align}
where $m_1$ and $m_3 $ are determined in Eq.(\ref{eqn:Lam010}).

If $\alpha \not \equiv 0$ $($mod $2\omega_1 \Zint \oplus 2\omega_3 \Zint)$, then we set
\begin{align}
 & \kappa = \zeta (m_1 \omega _3 -m_3 \omega _1 ) -m_1 \eta _3 +m_3 \eta _1 . \label{kapp}
\end{align}
It follows from Legendre's relation $\eta _1 \omega _3 - \eta_3 \omega _1 =\pi \sqrt{-1} /2$ and the relation $\zeta (-\alpha )=-\zeta (\alpha )$ that 
\begin{align}
& \exp (\kappa  (x+2\omega _{j}) ) \left( \frac{d}{dx} \right) ^{k} \Phi _i(x+2\omega _{j} , \alpha ) \label{eq:periodj'} \\
& = \exp (-2\eta _{j} \alpha +2\omega _{j} (\zeta (\alpha ) + \kappa )) \exp (\kappa x )\left( \frac{d}{dx} \right) ^{k} \Phi _i(x, \alpha ) \nonumber \\
& = \exp (2m_1 (\eta _{j} \omega _3- \eta _3\omega _j) + 2m_3 (\eta _{1} \omega _j- \eta _j \omega _1)) \exp (\kappa x )\left( \frac{d}{dx} \right) ^{k} \Phi _i(x, \alpha ) \nonumber \\
& = \exp (\pi \sqrt{-1} m_{j}) \exp (\kappa x )\left( \frac{d}{dx} \right) ^{k} \Phi _i(x, \alpha ) \nonumber
\end{align}
for $i=0,1,2,3$, $j=1,3$ and $k \in \Zint _{\geq 0}$.
Hence the function $\Lambda _g (x)$ and the functions $\exp (\kappa x )\left( \frac{d}{dx} \right) ^{k} \Phi _i(x, \alpha )$ have the same periodicity with respect to periods $(2\omega_1, 2\omega _3)$.
Since the meromorphic function $\Lambda _g (x)$ satisfies Eq.(\ref{eq:Hg}), the regular singular point $x=\pm \delta _{i'}$ $(i'=1,\dots ,M)$ is apparent, and the exponents at $x=\pm \delta _{i'} $ are $0$ and $r_{i'} +1$, it is holomorphic except for $\Zint \omega _1 \oplus \Zint \omega _3 $ and has a pole of degree $\tilde{l} _i$ or zero of degree $\tilde{l} _i +1$ at $x= \omega _i$ $(i=0,1,2,3)$.
The function $\exp(\kappa x) \left( \frac{d}{dx} \right) ^{k} \Phi _i(x , \alpha ) $ has a pole of degree $k+1$ at $x=\omega _i$.
By subtracting the functions $\exp(\kappa x) \left( \frac{d}{dx} \right) ^{k} \Phi _i(x , \alpha ) $ from the function $\Lambda _g (x)$ to erase the poles, we obtain a holomorphic function that has the same periods as $\Phi _0(x , \alpha ) $, and must be zero, because if we denote the holomorphic function by $f(x)$, then $f(x)/(\exp (\kappa x) \Phi _0 (x))$ is doubly-periodic and have only one pole of degree one in a fundamental domain, and $f(x)$ must be zero.
Hence we obtain the expression (\ref{Lalpha}). 
The periodicity (see Eq.(\ref{ellint})) follows from Eq.(\ref{eq:periodj'}).

If $\alpha \equiv 0$ $($mod $2\omega_1 \Zint \oplus 2\omega_3 \Zint)$ (i.e. $m_1 \omega _3 \equiv m_3 \omega _1$ $($mod $2\omega_1 \Zint \oplus 2\omega_3 \Zint)$), then we set
\begin{align}
 & \bar{\kappa }= -m_1 \eta _3 +m_3 \eta _1 .
\end{align}
The function $\Lambda _g (x)$ and the function $\exp (\bar{\kappa }x )$ have the same periodicity with respect to periods $(2\omega_1, 2\omega _3)$.
Hence the function $\Lambda _g (x) \exp (-\bar{\kappa }x )$ is doubly periodic, and we obtain the expression (\ref{Lalpha0}) by considering the poles. Periodicity (see Eq.(\ref{ellint0})) follows immediately.
\end{proof}

We investigate the situation that Eq.(\ref{eq:Hg}) has a non-zero solution of an elliptic function.
Let ${\mathcal F}_{\epsilon _1 , \epsilon _3 }$ and $\tilde{\mathcal F} _{\epsilon _1 , \epsilon _3 }$ $(\epsilon _1 , \epsilon _3 \in \{ \pm 1 \})$ be the spaces defined by
\begin{align}
& {\mathcal F} _{\epsilon _1 , \epsilon _3 }=\{ f(x) \mbox{; meromorphic }| f(x+2\omega_1)= \epsilon _1 f(x), \; f(x+2\omega_3)= \epsilon _3 f(x) \} , \\
& \tilde{\mathcal F} _{\epsilon _1 , \epsilon _3 }=\left\{ f(x) \: \left| 
\begin{array}{l}
f(x) \Psi _g (x) \in {\mathcal F}_{\epsilon _1 , \epsilon _3 }, \; \; f(x) \Psi _g (x)  \mbox{ is holomorphic}\\
\mbox{except for } \Zint \omega _1 \oplus \Zint \omega _3 \mbox{, and the degree of the pole at }\\
 x=\omega _i \mbox{ is no more than } \left\{ 
\begin{array}{ll}
l_i, &   i=1,2,3, \\
l_0+\sum _{i'=1}^M r_{i'}, & i=0.
\end{array}
\right.
\end{array}
\right. \right\} ,
\end{align}
where $(2\omega_1, 2\omega_3)$ are basic periods of elliptic functions. Then $\tilde{\mathcal F} _{\epsilon _1 , \epsilon _3 }$ is a finite-dimensional vector space. 
Note that, if a solution $f(x)$ to Eq.(\ref{eq:H}) satisfies the condition $f(x+2\omega _1) \Psi _g (x +2\omega _1) =\epsilon _1 f(x) \Psi _g (x) $ and $f(x+2\omega _3) \Psi _g (x +2\omega _3) =\epsilon _3 f(x) \Psi _g (x) $ for some $\epsilon _1 , \epsilon _3 \in \{ \pm 1 \}$, then we have $f(x) \in {\mathcal F} _{\epsilon _1 , \epsilon _3 }$, because the position of the poles and their degree are restricted by the differential equation.

\begin{prop} \label{prop:disttwoch} 
Assume that Eq.(\ref{eq:H}) has a non-zero solution in the space $\tilde{\mathcal F} _{\epsilon _1 , \epsilon _3 }$ for some $\epsilon _1 , \epsilon _3 \in \{ \pm 1 \}$. Then the signs $(\epsilon _1 , \epsilon _3 )$ are determined uniquely for each $E$, $\tilde{s}_{i'}$ $(i'=1, \dots ,M)$ etc.
\end{prop}
\begin{proof}
Assume that Eq.(\ref{eq:H}) has a non-zero solution in both the spaces $\tilde{\mathcal F} _{\epsilon _1 , \epsilon _3 }$ and $\tilde{\mathcal F} _{\epsilon '_1 , \epsilon '_3 }$. Let $f_1 (x)$ (resp. $f_2(x)$) be the solution to the differential equation (\ref{eq:H}) in the space $\tilde{\mathcal F} _{\epsilon _1 , \epsilon _3 }$ (resp. the space $\tilde{\mathcal F} _{\epsilon '_1 , \epsilon '_3 }$). 
Then periodicity of the function $f_1 (x)\Psi _g (x)$ and $f_2(x)\Psi _g (x)$ is different, more precisely there exists $j \in \{ 1,3\}$ such that 
\begin{equation}
\left\{ 
\begin{array}{ll}
f_1 (x +2\omega _j) \Psi _g (x +2\omega _j)= \pm f_1 (x) \Psi _g (x),\\
f_2 (x +2\omega _j) \Psi _g (x +2\omega _j)= \mp f_2 (x) \Psi _g (x).
\end{array}
\right.
\label{f1f2antp}
\end{equation}
Then the functions $f_1(x)$ and $f_2(x)$ are linearly independent.
Since the functions $f_1(x)$ and $f_2(x)$ satisfy Eq.(\ref{eq:H}), we have $\frac{d}{dx} \left( f_2(x) f'_1(x) -f_1(x) f'_2(x) \right) = f_2(x) f''_1(x) -  f_1(x) f''_2(x)= 0$.
Therefore $f_2(x) f'_1 (x) - f_1(x) f'_2 (x) =C$ for constants $C$, and $C$ is non-zero, which follows from linear independence. By Eq.(\ref{f1f2antp}), the function $(f_2(x) f'_1 (x) - f_1(x) f'_2 (x) )\Psi _g (x)^2$ is anti-periodic with respect to the period $2\omega _j$, but it contradicts to $C\neq 0$.
Hence, we proved that Eq.(\ref{eq:H}) does not have a non-zero solution in both the spaces $\tilde{\mathcal F} _{\epsilon _1 , \epsilon _3 }$ and $\tilde{\mathcal F} _{\epsilon '_1 , \epsilon '_3 }$.
\end{proof}

\begin{prop} \label{prop:Q0F}
If $Q=0$, then we have $\Lambda (x) \in \tilde{\mathcal F} _{\epsilon _1 , \epsilon _3 }$ for some $\epsilon _1 , \epsilon _3 \in \{ \pm 1 \}$.
\end{prop}
\begin{proof}
It follows from Eq.(\ref{integ1}) and the double-periodicity of the function $\Xi (x)\Psi _g (x)^2$ that
\begin{equation}
(\Lambda (x+2\omega _j) \Psi _g (x+2\omega _j))^2 = \Xi (x+2\omega _j)\Psi _g (x+2\omega _j)^2=\Xi (x)\Psi _g (x)^2= (\Lambda (x)\Psi _g (x) )^2,
\end{equation}
for $j=1,3$. Hence $\Lambda (x+2\omega _j) \Psi _g (x+2\omega _j)= \pm \Lambda (x)\Psi _g (x) $ $(j=1,3)$
and we have $\Lambda (x) \in  \tilde{\mathcal F} _{\epsilon _1 , \epsilon _3 }$ for some $\epsilon _1 , \epsilon _3 \in \{ \pm 1 \}$.
\end{proof}

It follows from Proposition \ref{prop:prod} that the dimension of the space of solutions to Eq.(\ref{prodDE}), which are even doubly-periodic, is no less than one.
Since the exponents of Eq.(\ref{prodDE}) at $x=0$ are $-2l_0$, $1$ and $2l_0 +2$, the dimension of the space of even solutions to Eq.(\ref{prodDE}) is at most two.
Hence, the dimension of the space of solutions to Eq.(\ref{prodDE}), which are even doubly-periodic, is one or two.

\begin{prop} \label{prop:twoF}
Assume that the dimension of the space of solutions to Eq.(\ref{prodDE}), which are even doubly-periodic, is two for some eigenvalue $E$. Then all solutions to Eq.(\ref{eq:H}) for the eigenvalue $E$ are contained in the space $\tilde{\mathcal F} _{\epsilon _1 , \epsilon _3 }$ for some $\epsilon _1 , \epsilon _3 \in \{ \pm 1 \}$.
\end{prop}
\begin{proof}
Since the differential equation (\ref{eq:H}) is invariant under the change of parity $x \leftrightarrow -x$ and exponents at $x=0$ are even one and odd one, a basis of the solutions to Eq.(\ref{eq:H}) is taken as $f_e (x)$ and $f_o (x)$ such that $f_ e(x)$ (resp. $f_o (x)$) satisfies $f_ e(-x)=f_e (x)$ (resp. $f_ o(-x)=-f_o (x)$).
Then the functions $f_ e(x)^2$ and $f_ o(x)^2$ are even and they are solutions to Eq.(\ref{prodDE}). Since the dimension of the space of even solutions to Eq.(\ref{prodDE}) is at most two, and the dimension of the space of solutions to Eq.(\ref{prodDE}), which are even doubly-periodic, is two, the even functions $f_e(x)^2$ and $f_o(x)^2$ must be doubly-periodic. Hence $(f_e(x+2\omega _j) \Psi _g (x+2\omega _j))^2 =(f_e(x) \Psi _g (x))^2$ $(j=1,3)$ and it follows that $f_e(x+2\omega _j) \Psi _g (x+2\omega _j)= \pm f_e(x)\Psi _g (x)$ $(j=1,3)$. Therefore we have $f_e (x) \in \tilde{\mathcal F} _{\epsilon _1 , \epsilon _3 }$ for some $\epsilon _1 , \epsilon _3 \in \{ \pm 1 \}$. Similarly we have $f_ o(x) \in \tilde{\mathcal F} _{\epsilon '_1 , \epsilon '_3 }$ for some $\epsilon '_1 , \epsilon '_3 \in \{ \pm 1 \}$, and it follows from Proposition \ref{prop:disttwoch} that $\epsilon '_j =\epsilon _j$ $(j=1,3)$. Since $f_e (x)$ and $f_o (x)$ are a basis of solutions to Eq.(\ref{eq:H}), all solutions to Eq.(\ref{eq:H}) are contained in the space $\tilde{\mathcal F} _{\epsilon _1 , \epsilon _3 }$.
\end{proof}

\begin{prop} \label{prop:Xionedim}
If $M=0$ or ($M=1$ and $r_1 =1$), then the dimension of the space of solutions to Eq.(\ref{prodDE}), which are even doubly-periodic, is one for all $E$.
\end{prop}
\begin{proof}
Assume that the dimension of the space of solutions to Eq.(\ref{prodDE}), which are even doubly-periodic, is two.
From Proposition \ref{prop:twoF}, all solutions to Eq.(\ref{eq:H}) are contained in the space $\tilde{\mathcal F} _{\epsilon _1 , \epsilon _3 }$ for some $\epsilon _1 , \epsilon _3 \in \{ \pm 1 \}$. Since the differential equation (\ref{eq:Hg}) is invariant under the change of parity $x \leftrightarrow -x$ and exponents at $x=0$ are even one and odd one, a basis of the solutions to Eq.(\ref{eq:Hg}) can be taken as $f_e(x)$ and $f_o (x)$ such that $f_e (x)$ (resp. $f_o (x)$) is even (resp. odd) function.
From the assumption that $l_i \in \Zint $ $(i=0,1,2,3)$ and that regular singular points $b_{i'}$ are apparent ($i'=1,\dots ,M$), the functions $f_e (x)$ and $f_o (x)$ are meromorphic.
Since the function $f_e (x)$ (resp. $f_o (x)$) satisfies Eq.(\ref{eq:Hg}), it does not have poles except for $\Zint \omega _1 \oplus \Zint \omega _3$.
Hence the function $f_e (x)$ admits the expression $f_e (x)= \wp _1 (x) ^{\tilde{\beta }_1}\wp _2 (x) ^{\tilde{\beta }_2}  \wp _3 (x) ^{\tilde{\beta }_3} (P^{(1)}(\wp (x))+\wp '(x) P^{(2)}(\wp (x)))$, where $\wp _i (x)$ $(i=1,2,3)$ are co-$\wp $ functions and $P^{(1)}(z)$, $P^{(2)} (z)$ are polynomials in $z$. Since the function $f_e (x) $ is even, we have $P^{(1)} (z)=0$ or $P^{(2)} (z)=0$. By combining with the relation $\wp '(z) = -2\wp _1 (z) \wp _2 (z) \wp _3 (z) $, the function $f_e (x)$ is expressed as 
\begin{equation}
f_e (x)= \wp _1 (x) ^{\beta _1}\wp _2 (x) ^{\beta _2}  \wp _3 (x) ^{\beta _3} P_e (\wp (x)),
\end{equation}
where $P_e(z)$ is a polynomial in $z$.
Because the exponents of Eq.(\ref{eq:Hg}) at $x=\omega _i$ $(i=1,2,3)$ are $-l_i$ and $l_{i} +1$, we have $\beta _i \in \{ -l_i ,l _{i} +1 \}$ $(i=1,2,3)$.
Similarly the function $f_o (x)$ is expressed as 
\begin{equation}
f_o(x)= \wp _1 (x) ^{\beta '_1}\wp _2 (x) ^{\beta '_2}  \wp _3 (x) ^{\beta '_3} P_o (\wp (x)),
\end{equation}
where $P_o(z)$ is a polynomial in $z$ and $\beta ' _i \in \{ -l_i ,l _{i} +1 \}$.

Since $ \wp (-x) =\wp (x)$, $ \wp _i (-x) =-\wp _i (x) $ $(i=1,2,3)$ and the parity of functions $f_e (x)$ and $f_o (x)$ is different, we have $\beta _1 +\beta _2 +\beta _3 \not \equiv \beta '_1 +\beta '_2 +\beta '_3 $ (mod $2$). Since $f_e ( x +2\omega _1) = (-1) ^{\beta _2 +\beta _3 } f_e (x)$, $f_o ( x +2\omega _1) = (-1) ^{\beta '_2 +\beta '_3 } f_o (x)$, $ f_e ( x +2\omega _3) = (-1) ^{\beta _1 +\beta _2 } f_e (x)$, $ f_o ( x +2\omega _3) = (-1) ^{\beta '_1 +\beta '_2 } f_o (x)$, we have $\beta _2 +\beta _3 \equiv \beta '_2 +\beta '_3 $ (mod $2$) and $\beta _1 +\beta _2 \equiv \beta '_1 +\beta '_2 $ (mod $2$). Hence we have $\beta _i \not \equiv \beta '_i $ (mod $2$) for $i=1,2,3$. Therefore $(\beta _i , \beta '_i ) = (-l_i , l _{i} +1)$ or $(\beta _i , \beta '_i ) = (l _{i} +1, -l_i)$ for each $i \in \{ 1,2,3 \}$.
Let $\beta _0$ (resp. $\beta '_0$) be the exponent of the function $f_e (x)$ (resp. $f_o (x)$) at $x=0$. Since the parity of functions $f_e (x)$ and $f_o (x)$ is different and the exponents of Eq.(\ref{eq:Hg}) at $x=0$ are $-l_0- \sum_{i'=1}^M r_{i'}$ and $l _{0} +1- \sum_{i'=1}^M r_{i'}$, we have $(\beta _0 , \beta '_0 ) = (-l_0 - \sum_{i'=1}^M r_{i'}, l _0 +1- \sum_{i'=1}^M r_{i'})$ or $(\beta _0, \beta '_0 ) = ( l _0 +1- \sum_{i'=1}^M r_{i'} ,-l_0 - \sum_{i'=1}^M r_{i'})$.

Since the function $f_e (x)$ is doubly-periodic with periods $(4\omega _1 , 4\omega _3)$, the sum of degrees of zeros of $f_e (x)$ on the basic domain is equal to the sum of degrees of poles of $f_e (x)$. Since the function $f_e (x)$ does not have poles except for $\Zint \omega _1 \oplus \Zint \omega _3$, we have $\sum _{i=0}^3 \beta _i \leq 0$. Similarly we have $\sum _{i=0}^3 \beta '_i \leq 0$.
Hence $0 \geq \sum _{i=0}^3 (\beta _i +\beta '_i) = 4-2 \sum_{i'=1}^M r_{i'}$.
Therefore we have $\sum_{i'=1}^M r_{i'} \geq 2$. 

Thus we obtain that, if $M=0$ or ($M=1$ and $r_1 =1$), then the dimension of the space of solutions to Eq.(\ref{prodDE}), which are even doubly-periodic, is one.
\end{proof}

Note that the case $M=0$ corresponds to Heun's equation, and the case $M=1$ and $r_1 =1$ is related with the sixth Painlev\'e equation.

\begin{exa}
Let us consider the following differential equation:
\begin{equation}
\left\{ -\left( \frac{d}{dx} \right) ^2 +\left( \frac{\wp ' (x)}{\wp (x) + \sqrt{\frac{g_2}{12}}}+\frac{\wp ' (x)}{\wp (x) - \sqrt{\frac{g_2}{12}}} \right) \frac{d}{dx} \right\} f(x)=0.
\label{deq:twodim}
\end{equation}
This equation corresponds to the case $l_0=1$, $l_1=l_2=l_3=0$, $M=2$ and $r_1=r_2=1$, if $g_2 \neq 0$.
From the relation
\begin{equation}
 \frac{\wp ' (x)}{\wp (x) + \sqrt{\frac{g_2}{12}}}+\frac{\wp ' (x)}{\wp (x) - \sqrt{\frac{g_2}{12}}} = \frac{\wp ''' (x)}{\wp '' (x)},
\end{equation}
a basis of the solutions to Eq.(\ref{deq:twodim}) is $1$, $\wp '(x)$.
The dimension of the solutions to Eq.(\ref{prodDE}), which are even doubly-periodic, is two, and a basis of the solutions to Eq.(\ref{prodDE}) is written as $1/\wp ''(x)$, $\wp '(x)^2/\wp ''(x)$, $\wp' (x)/\wp ''(x)$.
\end{exa}

\begin{prop} \label{prop:char}
Assume that the dimension of the space of the solutions to Eq.(\ref{prodDE}), which are even doubly-periodic, is one. 
Let $c_0$ and $b^{(i)}_j$ be constants defined in Eq.(\ref{Fx}).\\
(i) If there exists a non-zero solution to Eq.(\ref{eq:H}) in the space $\tilde{\mathcal F} _{\epsilon _1 , \epsilon _3 }$ for some $\epsilon _1 , \epsilon _3 \in \{ \pm 1 \}$, then we have $Q=0$.\\
(ii) If $Q \neq 0$ and $l_i \neq 0$, then the function $\Lambda (x)$ has a pole of degree $l_i$ and $b^{(i)}_0 \neq 0$.\\
(iii) If $Q \neq 0$ and $l_i =0$, then $\Lambda (\omega _i )\neq 0$ and $\Xi (\omega _i )\neq 0$. If $Q \neq 0$ and $l_0=l_1=l_2=l_3 =0$, then $c_0 \neq 0$.
\end{prop}
\begin{proof}
First we prove (i). Suppose that there exists a non-zero solution to Eq.(\ref{eq:H}) in the space $\tilde{\mathcal F} _{\epsilon _1 , \epsilon _3 }$ and $Q \neq 0$. From the condition $Q \neq 0$, the functions $\Lambda (x) $ and $\Lambda (-x) $ form the basis of the space of the solutions to the differential equation (\ref{eq:H}).
Since there is a non-zero solution to Eq.(\ref{eq:H}) in the space $\tilde{\mathcal F} _{\epsilon _1 , \epsilon _3 }$, there exist constants $(C_1, C_2) \neq (0,0)$ such that $C_1 \Lambda (x) +C_2 \Lambda (-x) \in \tilde{\mathcal F} _{\epsilon _1 , \epsilon _3 }$.
By shifting $x \rightarrow x +2\omega _j$ ($j=1,3$), it follows from Eq.(\ref{eqn:Lam010}) that
\begin{align}
& \quad (C_1 \Lambda (x +2\omega _j ) +C_2 \Lambda (-(x+2\omega _j) )) \Psi _g (x+2\omega _j) \\
& = C_1 \Lambda (x +2\omega _j ) \Psi _g (x+2\omega _j) \pm C_2 \Lambda (-x-2\omega _j) \Psi _g (-x-2\omega _j) \nonumber \\
& = C_1 \exp (\pi \sqrt{-1} m_j ) \Lambda (x) \Psi _g (x) \pm C_2 \exp (-\pi \sqrt{-1} m_j ) \Lambda (-x) \Psi _g (-x) \nonumber \\
& = (C_1 \exp (\pi \sqrt{-1} m_j ) \Lambda (x) +C_2 \exp (-\pi \sqrt{-1} m_j ) \Lambda (-x) ) \Psi _g (x) , \nonumber
\end{align}
where the sign $\pm $ is determined by the branching of the function $\Psi _g (x)$, and the function $C_1 \exp (\pi \sqrt{-1} m_j ) \Lambda (x) +C_2 \exp (-\pi \sqrt{-1} m_j ) \Lambda (-x)$ also satisfies Eq.(\ref{eq:H}).
On the other hand, it follows from the definition of the space $\tilde{\mathcal F} _{\epsilon _1 , \epsilon _3 }$ that $(C_1 \Lambda (x +2\omega _j ) +C_2 \Lambda (-(x+2\omega _j) )) \Psi _g (x+2\omega _j)= (C_1 \Lambda (x ) +C_2 \Lambda (-x ))\Psi _g (x)$ or $(C_1 \Lambda (x +2\omega _j ) +C_2 \Lambda (-(x+2\omega _j) )) \Psi _g (x+2\omega _j)= -(C_1 \Lambda (x ) +C_2 \Lambda (-x ))\Psi _g (x)$. By comparing two expressions, we have $\exp (\pi \sqrt{-1} m_j ) \in \{ \pm 1 \}$ ($j=1,3$) and the periodicities of the functions $\Lambda (x ) \Psi _g (x)$  and $(C_1 \Lambda (x ) +C_2 \Lambda (-x ))\Psi _g (x)$ coincide. Thus $ \Lambda (x ), \: \Lambda (-x ) \in \tilde{\mathcal F} _{\epsilon _1 , \epsilon _3 }$. The functions $\Lambda (x )^2$ and $\Lambda (-x )^2$ are even doubly-periodic function and satisfy Eq.(\ref{prodDE}), because they are the products of a pair of solutions to Eq.(\ref{eq:H}). Since the functions $\Lambda (x)$ and $\Lambda (-x)$ are linearly independent, the functions $\Lambda (x)^2$ and $\Lambda (-x)^2$ are linearly independent. Hence the dimension of the space of solutions to Eq.(\ref{eq:H}), which are even doubly-periodic, is no less than two, and contradict the assumption of the proposition.
Therefore the supposition $Q\neq 0$ is false, and we obtain (i).

Next we show (ii). Assume that $l_i \neq 0$. Since the exponents of Eq.(\ref{eq:H}) at $x=\omega_{i}$ are $-l_i$ or $l_i +1$, the function $\Lambda (x)$ has a pole of degree $l_i$ or a zero of degree $l_i +1$ at $x=\omega_{i}$. It follows from the periodicity (see Eq.(\ref{eqn:Lam010})) that, if the function $\Lambda (x)$ has a zero at $x=\omega_{i}$, then $\Lambda (x)$ has also a zero at $x=-\omega_{i}$.
Hence the function $\Lambda (-x)$ has a zero at $x=\omega_{i}$. From the assumption $Q \neq 0$, any solution to Eq.(\ref{eq:H}) is written as a linear combination of functions $\Lambda (x)$ and $\Lambda (-x)$. But it contradicts that one of the exponents at $x=\omega_{i}$ is $-l_i$. Hence the function $\Lambda (x)$ has a pole of degree $l_i$.
Since the dimension of the space of the solutions to Eq.(\ref{eq:H}), which are even doubly-periodic, is one, we have $\Xi (x)=C \Lambda (x) \Lambda (-x) $ for some non-zero constant $C$ and $b^{(i)}_0 \neq 0$.

(iii) is proved similarly by showing that the function $\Lambda (x)$ does not have zero at $x=\omega _i$.
\end{proof}

By combining Propositions \ref{prop:Q0F} and \ref{prop:char} (i) we obtain the following proposition:
\begin{prop} \label{prop:zeros}
Assume that the dimension of the space of solutions to Eq.(\ref{prodDE}), which are even doubly-periodic, is one.
Then the condition $Q=0$ is equivalent to that there exists a non-zero solution to Eq.(\ref{eq:H}) in the space $\tilde{\mathcal F} _{\epsilon _1 , \epsilon _3 }$ for some $\epsilon _1 , \epsilon _3 \in \{ \pm 1 \}$.
\end{prop}

We show that the function $\Lambda (x) $ admits an expression of the Bethe Ansatz type.
\begin{prop} \label{prop:BA}
Set $l=\sum _{i=0}^3 l_i +\sum _{i'=1}^M r_{i'} $, $\tilde{l} _0 = l_0 +\sum _{i'=1}^M r_{i'}$ and $\tilde{l}_i =l_i$ $(i=1,2,3)$. Assume that  $Q \neq 0$ and the dimension of the space of the solutions to Eq.(\ref{prodDE}), which are even doubly-periodic, is one.\\
(i) The function $\Lambda _g(x) $ in Eq.(\ref{integ1g}) is expressed as 
\begin{align}
& \Lambda _g(x) = \frac{C_0 \prod_{j=1}^l \sigma(x-t_j)}{\sigma(x)^{\tilde{l}_0}\sigma_1(x)^{\tilde{l}_1}\sigma_2(x)^{\tilde{l}_2}\sigma_3(x)^{\tilde{l}_3}}\exp \left(cx \right), 
\label{eq:tilL}
\end{align}
for some $t_1, \dots , t_l$, $c$ and $C_0 (\neq 0)$ such that $t_j \equiv 0$ $($mod $\omega_1 \Zint \oplus \omega_3 \Zint)$ for all $j$, where $\sigma _i (x)$ $(i=1,2,3)$ are co-sigma functions.\\
(ii) $t_j + t_{j'} \not \equiv 0$ $($mod $2\omega_1 \Zint \oplus 2\omega_3 \Zint)$ for $1\leq j<j'\leq l$.\\
(iii) If $t_j \not \equiv \pm \delta _{i'}$ $($mod $2\omega_1 \Zint \oplus 2\omega_3 \Zint)$ for all $i' \in \{1,\dots ,M\}$, then we have $t_j \not \equiv t_{j'}$ $($mod $2\omega_1 \Zint \oplus 2\omega_3 \Zint)$ for all $j' (\neq j)$.\\
(iv) If $t_j \equiv \pm \delta _{i'}$ $($mod $2\omega_1 \Zint \oplus 2\omega_3 \Zint)$, then $\# \{ j' \: | \: t_j \equiv t_{j'}$ $($mod $2\omega_1 \Zint \oplus 2\omega_3 \Zint) \} = r_{i'}+1$.
(v) If $l_0 \neq 0$ (resp. $l_0=0$), then we have $c=\sum_{i=1}^l \zeta(t_j)$ (resp. $c=\sum_{i=1}^l \zeta(t_j)+ \sqrt{-Q}/\Xi (0)$). (Note that it follows from Proposition \ref{prop:char} (iii) that $\sqrt{-Q}/\Xi (0)$ is finite.)\\
(vi) Set $z= \wp (x)$ and $z_j=\wp (t_j)$. Then 
\begin{equation}
\left. \frac{d\Xi (x)}{dz}\right| _{z=z_j}=\frac{2\sqrt{-Q}}{\wp'(t_j)}.
\label{signpptj}
\end{equation}
\end{prop}
\begin{proof}
Let $\alpha $ be the value defined in Eq.(\ref{al}).
First, we consider the case $\alpha \not \equiv 0$ $($mod $2\omega_1 \Zint \oplus 2\omega_3 \Zint)$. Let $\kappa $ be the value defined in Eq.(\ref{kapp}).
Then the function $\Lambda _g (x) / \left( \exp (\kappa x ) \Phi _0(x, \alpha ) \right) $ is meromorphic and doubly-periodic. Hence there exists $a_1 ,\dots ,a_{l'}$, $b_1 ,\dots ,b_{l'}$ such that $a_1 + \dots + a_{l'} = b_1 +\dots + b_{l'}$ and 
$$
\Lambda _g (x) / \left( \exp (\kappa x ) \Phi _0(x, \alpha ) \right)  = \frac{\prod_{j=1}^{l'} \sigma (x- a_j )}{\prod_{j=1}^{l'} \sigma (x- b_j )}.
$$
For the case $\alpha \equiv 0$ $($mod $2\omega_1 \Zint \oplus 2\omega_3 \Zint)$ the function $\Lambda _g (x) / \exp (\bar{\kappa }x)$ is similarly expressed as
$$
\Lambda _g (x) / \exp (\bar{\kappa }x) = \frac{\prod_{j=1}^{l'} \sigma (x- a_j )}{\prod_{j=1}^{l'} \sigma (x- b_j )}.
$$

Since the function $\Lambda _g (x)$ satisfies Eq.(\ref{eq:Hg}), it does not have poles except for $\omega_1 \Zint \oplus \omega_3 \Zint$ and we have
\begin{align}
& \Lambda _g (x) = \frac{C_0 \prod_{j=1}^l \sigma(x-t_j)}{\sigma(x)^{\tilde{l}_0}\sigma_1(x)^{\tilde{l}_1}\sigma_2(x)^{\tilde{l}_2}\sigma_3(x)^{\tilde{l}_3}}\exp \left(cx \right), 
\label{eq:tilL0}
\end{align}
for some $t_1, \dots , t_l$, $c$ and $C_0 (\neq 0)$.
It follows from Proposition \ref{prop:char} (ii) and (iii) that $t_j \not \equiv 0$ $($mod $\omega_1 \Zint \oplus \omega_3 \Zint)$.
Therefore we obtain (i).

Suppose that $t_j + t_{j'} \equiv 0$ $($mod $2\omega_1 \Zint \oplus 2\omega_3 \Zint)$ for some $j $ and $j'(\neq j)$,
From Eq.(\ref{eq:tilL}) and $-t_j \equiv t_{j'} $ $($mod $2\omega_1 \Zint \oplus 2\omega_3 \Zint)$, we have $\Lambda _g (t_j) = \Lambda  _g (-t_{j}) =0$. Since $Q\neq 0$, all solutions to Eq.(\ref{eq:Hg}) are written as linear combinations of $\Lambda  _g (x)$ and $\Lambda  _g (-x) $. Hence $t_j$ is a zero for all solutions to Eq.(\ref{eq:Hg}), but they contradict that one of the exponents at $x=t_j$ is zero. Therefore we obtain (ii).

If $t_j \not \equiv \pm \delta _{i'}, \omega _i$ $($mod $2\omega_1 \Zint \oplus 2\omega_3 \Zint)$ for all $i$ and $i'$, then the exponents of Eq.(\ref{eq:Hg}) at $x=t_j$ are $0$ and $1$, and $x=t_j$ is a zero of $\Lambda _g (x)$ of degree one.
Incidentally, the exponents of Eq.(\ref{eq:Hg}) at $x=\pm \delta _{i'}$ are $0$ and $r_{i'}+1$. Hence, if $t_j \equiv \pm \delta _{i'}$ $($mod $2\omega_1 \Zint \oplus 2\omega_3 \Zint)$, then $x=t_j$ is a zero of $\Lambda _g (x)$ of degree $r_{i'}+1$. Thus we obtain (iii) and (iv).

It follows from Eq.(\ref{eq:tilL}) and $\Lambda _g(x)=\Psi _g(x) \Lambda (x)$ that
\begin{equation}
\frac{\Lambda ' (x)}{\Lambda (x)} = c- \tilde{l}_0 \frac{\sigma '(x)}{\sigma (x)} -\sum _{i=1}^3 \tilde{l}_i \frac{\sigma ' _i (x)}{\sigma _i (x)}+\sum_{j=1}^l \frac{\sigma '(x-t_j)}{\sigma (x-t_j)} -\sum _{i' =1}^M \frac{r_{i'}}{2} \frac{\wp '(x)}{\wp (x) -\wp (\delta _{i'})}.
\end{equation}
By expanding Eq.(\ref{eq:llp}) at $x=0$ and observing coefficient of $x^0$, we obtain
\begin{equation}
c -\sum_{j=1}^l \zeta(t_j )= \left. \frac{\sqrt{-Q}}{\Xi (x)}\right| _{x=0},
\end{equation}
because the functions $\sigma '(x)/\sigma (x)$, $\sigma ' _i (x)/\sigma _i (x)$, $\wp '(x)/(\wp (x) -\wp (\delta _{i'}))$ and $\Xi '(x)/\Xi (x)$ are odd and $\sigma '(-t)/\sigma (-t)= -\zeta (t)$.  
It follows from $Q \neq 0$ and Proposition \ref{prop:char} that, if $l_0 \neq 0$, then $\left. \sqrt{-Q}/\Xi (x)\right| _{x=0} =0$, and if $l_0 = 0$, then $\left. \sqrt{-Q}/\Xi (x)\right| _{x=0} $ is finite. Thus we obtain (v).

We show (vi). 
The function $\Lambda (x) \Lambda (-x) $ is even doubly-periodic and satisfies Eq.(\ref{prodDE}), because it is a product of the solutions to Eq.(\ref{eq:H}). Since the dimension of the space of the solutions to Eq.(\ref{eq:H}), which are even doubly-periodic, is one, we have $\Xi (x)=C \Lambda (x) \Lambda (-x) $ for some non-zero constant $C$. Hence we have $\Xi (t_j ) = \Xi (-t_j )=0$. 
On the other hand, we have $\Lambda (-t_j) \neq 0$ from (ii).
At $x=-t_j$, the l.h.s. of Eq.(\ref{eq:llp}) is finite, and the denominator of the r.h.s. is zero. Therefore we have 
\begin{equation}
\Xi ' (x)| _{x=-t_j} +2\sqrt{-Q} =0.
\end{equation}
By changing the variable $z=\wp (x)$ and the oddness of the function $\wp '(x)$, we obtain (vi).
\end{proof}
Note that Gesztesy and Weikard \cite{GW3} obtained a similar expression to Eq.(\ref{eq:tilL}) in the framework of Picard's potential.

\section{The case $M=1$, $r_1 =1$ and Painlev\'e equation} \label{sec:P6}
\subsection{}
We consider Eq.(\ref{eq:Hg}) for the case $M=1$, $r_1 =1$.
For this case, Eq.(\ref{eq:Hg}) is written as 
\begin{equation}
(H_g-\tilde{E})f_g (x)=0,
\label{Hgkr1}
\end{equation}
where
\begin{align}
&  H_g= -\frac{d^2}{dx^2} + \frac{\wp ' (x)}{\wp (x) -\wp (\delta _{1})} \frac{d}{dx} + \frac{\tilde{s}_{1}}{\wp (x) -\wp (\delta _{1})} +\sum_{i=0}^3 l_i(l_i+1) \wp (x+\omega_i) . \label{HgP6}
\end{align}
We set
\begin{align}
& \Psi _g (x)= \sqrt{\wp (x) -\wp (\delta _1)} , \quad b_1 =\wp (\delta _1) , \label{pgkr0} \\
& \mu _1= \frac{-\tilde{s}_{1}}{4 b_1^3 -g_2 b_1 -g_3} +\sum _{i=1}^3 \frac{l_i}{2(b_1-e_i)}, \\
& p= \tilde{E} -2(l_1l_2e_3+l_2l_3e_1+l_3l_1e_2) +\sum _{i=1}^3 l_i (l_i e_i +2(e_i+b_1 )) .  \label{pgkr}
\end{align}
The condition that, the regular singular points $x=\pm \delta _1$ is apparent, is written as
\begin{align}
& p=(4 b_1^3 -g_2 b_1 -g_3) \left\{ -\mu _{1} ^2 +\sum_{i=1}^3\frac{l_i+\frac{1}{2}}{b_1-e_i} \mu _{1} \right\} \label{pgkrap} \\
& \quad \quad -b_1 (l_1+l_2+l_3-l_0)(l_1+l_2+l_3+l_0+1) .\nonumber
\end{align}
From now on we assume that $l_0, l_1, l_2, l_3 \in \Zint _{\geq 0}$ and the eigenvalue $\tilde{E}$ satisfies Eqs.(\ref{pgkr}, \ref{pgkrap}). Then the assumption in Proposition \ref{prop:prod} is true, and propositions and theorem in the previous section are valid. The function $\Xi (x)$ in Proposition \ref{prop:prod} is written as 
\begin{equation}
\Xi (x)=c_0+\frac{d_0}{(\wp (x)-\wp (\delta _1))}+\sum_{i=0}^3 \sum_{j=0}^{l_i-1} b^{(i)}_j \wp (x+\omega_i)^{l_i-j} .
\label{Fxkr1}
\end{equation}
It follows from Proposition \ref{prop:Xionedim} that the function $\Xi (x)$ is determined uniquely up to multiplicative constant.
Ratios of the coefficients $c_0/d_0$ and $b^{(i)}_j/d_0$ $(i=0,1,2,3, \: j=0,\dots ,l_i-1)$ are written as rational functions in variables $b_1$ and $\mu _{1} $, because the coefficients $b^{(i)}_j$, $c_0$ and $d_0$ satisfy linear equations whose coefficients are rational functions in $b_1$ and $\mu _1$, which are obtained by substituting Eq.(\ref{Fxkr1}) into Eq.(\ref{prodDE}). The value $Q$ is calculated by Eq.(\ref{const}) and it is expressed as a rational function in $b_1$ and $\mu _{1} $ multiplied by $d_0 ^2$. It is shown by observing asymptotic $\mu _1 \rightarrow \infty$ that $Q$ is not identically zero. By an appropriate choice of $d_0$, $Q$ is expressed as a polynomial in $b_1$ and $\mu _1$. We set 
\begin{equation}
\Lambda _g ( x)=\Psi _g(x) \sqrt{\Xi (x)}\exp \int \frac{ \sqrt{-Q}dx}{\Xi (x)}.
\label{integ1P6}
\end{equation}
Due to Proposition \ref{prop:Linteg}, the function $\Lambda _g (x)$ is a solution to the differential equation (\ref{Hgkr1}). 
By Theorem \ref{thm:alpha}, the eigenfunction $\Lambda _g (x)$ is also expressed in the form of the Hermite-Krichever Ansatz. Namely, it is expressed as 
\begin{align}
& \Lambda _g (x) = \exp \left( \kappa x \right) \left( \sum _{i=0}^3 \sum_{j=0}^{\tilde{l}_i-1} \tilde{b} ^{(i)}_j \left( \frac{d}{dx} \right) ^{j} \Phi _i(x, \alpha ) \right)
\label{LalphaP6}
\end{align}
or
\begin{align}
& \Lambda _g (x) = \exp \left( \bar{\kappa } x \right) \left( \bar{c} +\sum _{i=0}^3 \sum_{j=0}^{\tilde{l}_i-2} \bar{b} ^{(i)}_j \left( \frac{d}{dx} \right) ^{j} \wp (x+\omega _i) +\sum_{i=1}^3 \bar{c}_i \frac{\wp '(x)}{\wp (x)-e_i} \right)
\label{Lalpha0P6}
\end{align}
where $\tilde{l} _0 = l_0 +1$ and $\tilde{l}_i =l_i$ $(i=1,2,3)$. Now we investigate the values $\alpha $ and $\kappa $ in Eq.(\ref{LalphaP6}). Note that, if $\alpha \not \equiv 0$ (mod $2\omega _1 \Zint \oplus 2\omega _3 \Zint $), then the function $\Lambda _g (x)$ is expressed as Eq.(\ref{LalphaP6}) and we have
\begin{align}
& \Lambda _g (x+2\omega _j) = \exp (-2\eta _j \alpha +2\omega _j \zeta (\alpha ) +2 \kappa \omega _j ) \Lambda _g (x) , \quad  (j=1,3). \label{ellint00} 
\end{align}
\begin{prop} \label{prop:P6HK}
Assume that $M=1$, $r_1 =1$, $l_0, l_1, l_2, l_3 \in \Zint _{\geq 0}$ and the value $p$ satisfies Eq.(\ref{pgkrap}). Let $\alpha $ and $\kappa $ be the values determined by the Hermite-Krichever Ansatz  (see Eq.(\ref{LalphaP6})).
Then $\wp (\alpha )$ is expressed as a rational function in variables $b_1$ and $\mu _{1} $, $\wp ' (\alpha )$ is expressed as a product of $\sqrt {-Q}$ and a rational function in variables $b_1$ and $\mu _{1} $, and $\kappa $ is expressed as a product of $\sqrt {-Q}$ and a rational function in variables $b_1$ and $\mu _{1} $.
\end{prop}
\begin{proof}
It follows from Eqs.(\ref{eq:tilL}, \ref{periods}, \ref{rel:sigmai}) that
\begin{align} 
& \Lambda _g (x+2\omega _j ) = \exp \left( 2\eta _j \left( - \sum_{j'=1}^l t_{j'} + \sum_{i=1}^3 l_{i}\omega_{i} \right)  +2 \omega _j \left( c - \sum _{i=1}^3 l_i \eta_i \right) \right) \Lambda _g (x), \label{ellinttj} 
\end{align}
for $j=1,3$, where $l= l_0+l_1+l_2+l_3 +1$. By comparing with Eq.(\ref{ellint00}), we have
\begin{align}
& -2\eta _1 \alpha +2\omega _1 (\zeta (\alpha ) + \kappa )= -2\eta _1 \left( \sum_{j'=1}^l t_{j'} - \sum_{i=1}^3 l_{i}\omega_{i}\right)  +2\omega _1 \left( c - \sum _{i=1}^3 l_i \eta_i \right) +2\pi \sqrt{-1}n_1 \label{eq:n1} ,\\
& -2\eta _3 \alpha +2\omega _3 (\zeta (\alpha ) + \kappa )= -2\eta _3 \left(\sum_{j'=1}^l t_{j'} - \sum_{i=1}^3 l_{i}\omega_{i}\right)  +2\omega _3 \left( c - \sum _{i=1}^3 l_i \eta_i \right)  +2\pi \sqrt{-1}n_3, \label{eq:n3}
\end{align}
for integers $n_1$, $n_3$. It follows that
\begin{align}
& \left( \alpha - \left(\sum_{j'=1}^l t_{j'} - \sum_{i=1}^3 l_{i}\omega_{i}\right) \right)(-2\eta_1 \omega _3 + 2\eta_3 \omega _1 ) = 2\pi \sqrt {-1} (n_1\omega _3 -n_3 \omega _1 ) , \label{eq:aln} \\
& \left( \zeta (\alpha ) + \kappa -c+ \sum _{i=1}^3 l_i \eta_i \right)(2\eta_3 \omega _1 -2\eta_1 \omega _3 ) = 2\pi \sqrt {-1} (n_1\eta _3 -n_3 \eta _1 ) . \label{eq:kan} 
\end{align}
From Legendre's relation $\eta_1 \omega _3 - \eta_3 \omega _1 = \pi \sqrt {-1}/2$, we have 
\begin{equation}
\alpha \equiv \sum_{j'=1}^l t_{j'} - \sum_{i=1}^3 l_{i}\omega_{i} \quad \quad  (\mbox{mod }2\omega_1 \Zint \oplus 2\omega_3 \Zint).
\end{equation}
Combining Eqs.(\ref{eq:aln}, \ref{eq:kan}) with Proposition \ref{prop:BA} (v) and relations $\zeta (\alpha +2\omega _j )= \zeta (\alpha ) +2\eta _j$ $(j=1,3)$, we have 
\begin{equation}
\kappa = -\zeta \left(\sum_{j'=1}^{l} t_{j'} - \sum_{i=1}^3 l_i \omega_i \right) + \sum_{j'=1}^{l} \zeta (t_{j'})- \sum_{i=1}^3 l_i \eta_i + \delta _{l_0,0} \frac{\sqrt{-Q}}{\Xi (0)}.
\end{equation}

Next, we investigate values $\wp (\alpha )$, $\wp '(\alpha )$ and $\kappa $.
The functions $\wp (\sum_{j=1}^l t_j - \sum_{i=1}^3 l_{i}\omega_{i} )$, $\wp '(\sum_{j=1}^l t_j - \sum_{i=1}^3 l_{i}\omega_{i} )$ and 
$\zeta (\sum_{j=1}^{l} t_j - \sum_{i=1}^3 l_{i}\omega_{i}) - \sum_{j=1}^{l} \zeta (t_j) + \sum_{i=1}^3 l_{i}\eta_{i}$ are doubly-periodic in variables $t_1 ,\dots ,t_l$.
Hence by applying addition formulae of elliptic functions and considering the parity of functions $\wp (x)$, $\wp '(x)$ and $\zeta (x)$, we obtain the expression
\begin{align}
& \wp \left(\sum_{j=1}^l t_j - \sum_{i=1}^3 l_{i}\omega_{i}\right) = \sum_{j_1<j_2<\dots <j_m\atop{m:\mbox{\scriptsize{ even}}}}f^{(1)}_{j_1, \dots ,j_m}(\wp(t_1), \dots ,\wp (t_l)) \wp'(t_{j_1}) \dots \wp'(t_{j_l}), \\
& \wp '\left(\sum_{j=1}^l t_j - \sum_{i=1}^3 l_{i}\omega_{i}\right) = \sum_{j_1<j_2<\dots <j_m\atop{m:\mbox{\scriptsize{ odd}}}}f^{(2)}_{j_1, \dots ,j_m}(\wp(t_1), \dots ,\wp (t_l)) \wp'(t_{j_1}) \dots \wp'(t_{j_l}), \nonumber \\
& \zeta \left(\sum_{j=1}^{l} t_j - \sum_{i=1}^3 l_{i}\omega_{i}\right) - \sum_{j=1}^{l} \zeta (t_j) + \sum_{i=1}^3 l_{i}\eta_{i} \nonumber \\
& = \sum_{j_1<j_2<\dots <j_m\atop{m:\mbox{\scriptsize{ odd}}}}f^{(3)}_{j_1, \dots ,j_m}(\wp(t_1), \dots ,\wp (t_l)) \wp'(t_{j_1}) \dots \wp'(t_{j_l}), \nonumber
\end{align}
where $f^{(i)}_{j_1, \dots ,j_m}(x_1, \dots , x_l)$ $(i=1,2,3)$ are rational functions in $x_1 ,\dots ,x_l$.
From Eq.(\ref{signpptj}), the function  $\wp '(t_j)/\sqrt{-Q}$ is expressed as a rational function in $b_1$, $\mu _{1} $ and $\wp (t_j)$.
Hence, $\wp (\sum_{j=1}^l t_j - \sum_{i=1}^3 l_{i}\omega_{i})$, $\wp '(\sum_{j=1}^l t_j - \sum_{i=1}^3 l_{i}\omega_{i})/\sqrt{-Q}$ and $(\zeta (\sum_{j=1}^{l} t_j- \sum_{i=1}^3 l_{i}\omega_{i}) - \sum_{j=1}^{l} \zeta (t_j) + \sum_{i=1}^3 l_{i}\eta_{i})/\sqrt{-Q}$ are expressed as rational functions in the variable $\wp(t_1), \dots ,\wp (t_l)$, $b_1$ and $\mu _{1} $, and they are symmetric in $\wp(t_1), \dots ,\wp (t_l)$.

Since the dimension of the space of the solutions to Eq.(\ref{eq:H}), which are even doubly-periodic, is one, we have $\Xi (x)=C \Lambda (x) \Lambda (-x) $ for some non-zero scalar $C$.
Hence, we have the following expression;
\begin{equation}
\Xi (x) \Psi _g(x) ^2 =\frac{D \prod_{j=1}^{l}(\wp(x)-\wp(t_j))}{(\wp(x)-e_1)^{l_1}(\wp(x)-e_2)^{l_2}(\wp(x)-e_3)^{l_3}} \label{Fxtj}
\end{equation}
for some value $D (\neq 0)$. Thus
\begin{equation}
\prod_{j=1}^{l}(\wp(x)-\wp(t_j))= \Xi (x) \Psi _g (x) ^2 (\wp(x)-e_1)^{l_1}(\wp(x)-e_2)^{l_2}(\wp(x)-e_3)^{l_3}/D.
\end{equation}
Hence, the elementary symmetric functions $ \sum _{j_1<\dots <j_{l'}} \wp(t_{j_1})\dots \wp (t_{j_{l'}})$ ($l'=1,\dots ,l$) are expressed as rational functions in $b_1$ and $\mu _{1} $.
By substituting elementary symmetric functions into the symmetric expressions of $\wp (\sum_{j=1}^l t_j - \sum_{i=1}^3 l_{i}\omega_{i})$, $\wp (\sum_{j=1}^l t_j - \sum_{i=1}^3 l_{i}\omega_{i})$ and $(\zeta (\sum_{j=1}^{l} t_j- \sum_{i=1}^3 l_{i}\omega_{i}) - \sum_{j=1}^{l} \zeta (t_j) + \sum_{i=1}^3 l_{i}\eta_{i})/\sqrt{-Q}$, it follows that $\wp (\sum_{j=1}^l t_j - \sum_{i=1}^3 l_{i}\omega_{i})$, $\wp '(\sum_{j=1}^l t_j - \sum_{i=1}^3 l_{i}\omega_{i})/\sqrt{-Q}$ and $(\zeta (\sum_{j=1}^{l} t_j- \sum_{i=1}^3 l_{i}\omega_{i}) - \sum_{j=1}^{l} \zeta (t_j) + \sum_{i=1}^3 l_{i}\eta_{i})/\sqrt{-Q}$ are expressed as rational functions in $b_1$ and $\mu _{1} $.
Hence, $\wp (\alpha )$, $\wp '(\alpha )/ \sqrt{-Q}$ and $\kappa / \sqrt{-Q}$ are expressed as rational functions in variables $b_1$ and $\mu _{1} $.
\end{proof}

We now discuss the relationship between the monodromy preserving deformation of Fuchsian equations and the sixth Painlev\'e equation. For this purpose we recall some definitions and results of Painlev\'e equation.

The sixth Painlev\'e equation is a non-linear ordinary differential equation written as
\begin{align}
\frac{d^2\lambda }{dt^2} = & \frac{1}{2} \left( \frac{1}{\lambda }+\frac{1}{\lambda -1}+\frac{1}{\lambda -t} \right) \left( \frac{d\lambda }{dt} \right) ^2 -\left( \frac {1}{t} +\frac {1}{t-1} +\frac {1}{\lambda -t} \right)\frac{d\lambda }{dt} \label{eq:P6eqn} \\
& +\frac{\lambda (\lambda -1)(\lambda -t)}{t^2(t-1)^2}\left\{ \frac{\kappa _{\infty}^2}{2} -\frac{\kappa _{0}^2}{2}\frac{t}{\lambda ^2} +\frac{\kappa _{1}^2}{2}\frac{(t-1)}{(\lambda -1)^2} +\frac{(1-\kappa _{t}^2)}{2}\frac{t(t-1)}{(\lambda -t)^2} \right\}. \nonumber
\end{align}
A remarkable property of this differential equation is that its solutions do not have movable singularities other than poles.
This equation is also written in terms of a Hamiltonian system by adding the variable $\mu$, which is called the sixth Painlev\'e system:
\begin{equation}
\frac{d\lambda }{dt} =\frac{\partial H_{VI}}{\partial \mu}, \quad \quad
\frac{d\mu }{dt} =-\frac{\partial H_{VI}}{\partial \lambda},
\label{eq:Psys}
\end{equation}
with the Hamiltonian 
\begin{align}
H_{VI} = & \frac{1}{t(t-1)} \left\{ \lambda (\lambda -1) (\lambda -t) \mu^2 \right. \label{eq:P6} \\
& \left. -\left\{ \kappa _0 (\lambda -1) (\lambda -t)+\kappa _1 \lambda (\lambda -t) +(\kappa _t -1) \lambda (\lambda -1) \right\} \mu +\kappa (\lambda -t)\right\} ,\nonumber
\end{align}
where $\kappa = ((\kappa _0 +\kappa _1 +\kappa _t -1)^2- \kappa _{\infty} ^2)/4$.
The sixth Painlev\'e equation for $\lambda $ is obtained by eliminating $\mu $ in Eq.(\ref{eq:Psys}).
Set $\omega _1=1/2$, $\omega _3=\tau /2$ and write
\begin{equation}
t= \frac{e_3- e_1}{e_2-e_1}, \quad \lambda = \frac{\wp (\delta )- e_1}{e_2-e_1}.
\end{equation}
Then the sixth Painlev\'e equation is equivalent to the following equation (see \cite{Man,Tks}):
\begin{equation}
\frac{d^2 \delta }{d \tau ^2} = -\frac{1}{4\pi ^2} \left\{ \frac{\kappa _{\infty}^2}{2} \wp ' \left(\delta  \right) + \frac{\kappa _{0}^2}{2} \wp ' \left(\delta +\frac{1}{2} \right) + \frac{\kappa _{1}^2}{2} \wp ' \left(\delta +\frac{\tau +1}{2} \right) +  \frac{\kappa _{t}^2}{2} \wp ' \left(\delta +\frac{\tau }{2}\right) \right\}, \label{eq:P6ellip}
\end{equation}
where $\wp ' (z ) = (\partial /\partial z ) \wp (z)$.

It is widely known that the sixth Painlev\'e equation is obtained by the monodnomy preserving deformation of a certain linear differential equation.
Let us introduce the following Fuchsian differential equation:
\begin{equation}
\frac{d^2y}{dw^2} + p_1 (w) \frac{dy}{dw} +p_2(w) y=0, \label{eq:mpdP6}
\end{equation}
where 
\begin{align}
& p_1 (w) = \frac{1-\kappa _0}{w} + \frac{1-\kappa _1}{w-1} + \frac{1-\kappa _t}{w-t} -\frac{1}{w-\lambda}, \\
& p_2 (w) = \frac{\kappa }{w(w-1)} -\frac{t(t-1) H_{VI}}{w(w-1)(w-t)} + \frac{\lambda (\lambda -1) \mu}{w(w-1)(w-\lambda)}.
\end{align}
This equation has five regular singular points $\{ 0,1,t,\infty ,\lambda \}$ and the exponents at $w=\lambda $ are $0$ and $2$.
It follows from Eq.(\ref{eq:P6}) that the regular singular point $w=\lambda $ is apparent.
Then the sixth Painlev\'e equation is obtained by the monodromy preserving deformation of Eq.(\ref{eq:Psys}), i.e., the condition that the monodromy of Eq.(\ref{eq:mpdP6}) is preserved as deforming the variable $t$ is equivalent to that $\mu $ and $\lambda $ satisfy the Painlev\'e system (see Eq.(\ref{eq:Psys})), provided $\kappa _0, \kappa _1, \kappa _t , \kappa _{\infty} \not \in \Zint$. For details, see \cite{IKSY}.

Now we transform Eq.(\ref{eq:mpdP6}) into the form of Eq.(\ref{HgP6}). We set 
\begin{align}
& w=\frac{\wp (x) -e_1}{e_2-e_1}, \quad y= f_g (x) \prod _{i=1}^3 (\wp (x)-e_i)^{l_i/2}, \label{eq:wwpx} \\
&  \quad t=\frac{e_3-e_1}{e_2-e_1}, \quad \lambda =\frac{b_1 -e_1}{e_2 -e_1}, \quad \wp (\delta _1) =b_1.
\end{align}
Then we obtain Eq.(\ref{HgP6}) by setting
\begin{align}
& \kappa _0 =l_1 +1/2, \quad \kappa _1 =l_2 +1/2, \quad \kappa _t =l_3 +1/2, \quad \kappa _{\infty} =l_0 +1/2, \label{eq:kili} \\
& \mu = (e_2-e_1)\mu _1, \quad \kappa = (l_1+l_2+l_3+l_0 +1)(l_1+l_2+l_3-l_0) , \\
& H_{VI}=\frac{1}{t(1-t)} \left\{ \frac{p+\kappa e_3}{e_2-e_1} +\lambda (1- \lambda)\mu \right\},
\end{align}
(see Eqs.(\ref{pgkr0}--\ref{pgkr})), and Eq.(\ref{eq:P6}) is equivalent to Eq.(\ref{pgkrap}), that means that the apparency of regular singularity is inheritted.
Mapping from the variable $x$ to the variable $w$ (see Eq.(\ref{eq:wwpx})) is a double covering from the punctured torus $(\Cplx / (2\omega _1 \Zint + 2\omega _3 \Zint)) \setminus \{ 0, \omega_1 , \omega _2, \omega _3 \} $ to the punctured Riemann sphere ${\mathbb P}^1 \setminus \{ 0,1,t,\infty \}$. A solution $y(w)$ to Eq.(\ref{eq:mpdP6}) corresponds to a solution $f _g (x) $ to Eq.(\ref{HgP6}) by $y(w)= f_g (x) \prod _{i=1}^3 (\wp (x)-e_i)^{l_i/2}$. Hence the monodromy preserving deformation of Eq.(\ref{eq:mpdP6}) in $t$ corresponds to the monodromy preserving deformation of Eq.(\ref{HgP6}) in $\tau $.

Now we consider monodromy preserving deformation in the variable $\tau$ ($\omega _1 =1/2, \omega _3=\tau /2$) by applying solutions obtained by the Hermite-Krichever Ansatz for the case $l_i \in \Zint _{\geq 0}$ $(i=0,1,2,3)$.
Let $\alpha $ and $\kappa $ be values determined by the Hermite-Krichever Ansats (see Eq.(\ref{LalphaP6})). We consider the case $Q\neq 0$. Then a basis for solutions to Eq.(\ref{eq:Hg}) is given by $\Lambda _g(x)$ and $\Lambda _g(-x)$, and the monodromy matrix with respect to the cycle $x \rightarrow x+2\omega _j $ ($j=1,3$) is diagonal. The elements of the matrix are obtained from Eq.(\ref{ellint00}).
Hence, the eigenvalues $\exp (\pm ( -2\eta _j \alpha +2\omega _j \zeta (\alpha ) +2 \kappa \omega _j))$ $(j=1,3)$ of the monodromy matrices are preserved by the monodromy preserving deformation.
We set 
\begin{align}
& -2\eta _1 \alpha +2\omega _1 \zeta (\alpha ) +2 \kappa \omega _1 = \pi \sqrt{-1} C_1, \\
& -2\eta _3 \alpha +2\omega _3 \zeta (\alpha ) +2 \kappa \omega _3 = \pi \sqrt{-1} C_3, 
\end{align}
for contants $C_1$ and $C_3$.
By Legendre's relation, we have 
\begin{align}
& \alpha  = C_3 \omega _1 -C_1 \omega _3  \label{al00},\\
& \kappa = \zeta (C_1 \omega _3 -C_3 \omega _1 ) +C_3 \eta _1 -C_1 \eta _3  , \label{kapp00}
\end{align}
(see Eqs.(\ref{al}, \ref{kapp})). From Proposition \ref{prop:P6HK}, the value $\wp (\alpha)(=\wp (C_3 \omega _1-C_1 \omega _3 ))$ is expressed as a rational function in variables $b_1$ and $\mu _1$, the value $\wp '(\alpha)(=\wp '(C_3 \omega _1-C_1 \omega _3 ))$ is expressed as a product of $\sqrt {-Q}$ and a rational function in variables $b_1$ and $\mu _1$, and the value $\kappa (= \zeta (C_1 \omega _3 -C_3 \omega _1 ) +C_3 \eta _1 -C_1 \eta _3)$ is expressed as a product of $\sqrt {-Q}$ and rational function in variables $b_1$ and $\mu _1$.
By solving these equations for $b_1$ and $\mu _1$ and evaluating them into Eq.(\ref{HgP6}), the monodromy of the solutions on the cycles $x \rightarrow x +2\omega _j$ $(j=1,3)$ is preserved for the fixed values $C_1$ and $C_3$. 
Let $\gamma _0$ be the path in the $x$-plane which is obtained by the pullback of the cycle turning the origin around anti-clockwise in the $w$-plane, where $x$ and $w$ are related with $w=(\wp (x) -e_1)/(e_2 -e_1)$.
Then the monodromy matrix on $\gamma _0$ with respect to the basis $(\Lambda _g(x) ,\Lambda _g(-x) )$ is written as 
\begin{equation}
(\Lambda _g(x) ,\Lambda _g(-x) ) \rightarrow (\Lambda _g(-x) ,\Lambda _g(x) ) = (\Lambda _g(x) ,\Lambda _g(-x) ) 
\left(
\begin{array}{cc}
0 & 1 \\
1 & 0
\end{array}
\right) ,
\end{equation}
and does not depend on $\tau $.
Since the fundamental group on the punctured Riemann sphere ${\mathbb P}^1 \setminus \{ 0,1,t,\infty \}$ is generated by the images of $\gamma _0$ and the cycles $x \rightarrow x +2\omega _j$ $(j=1,3)$, Eqs.(\ref{al00}, \ref{kapp00}) describe the condition for the monodromy preserving deformation on the punctured Riemann sphere by rewriting the variable $\tau$ to $t$.
Summarizing, we have the following proposition.
\begin{prop} \label{prop:P6}
We set $\omega _1= 1/2$, $\omega _3 =\tau /2$ and assume that $l_i \in \Zint _{\geq 0}$ $(i=0,1,2,3)$ and $Q\neq 0$.
By solving the equations in Proposition \ref{prop:P6HK} in variable $b_1 =\wp (\delta _1)$ and $\mu _1$, we express $\wp (\delta _1)$ and $\mu _1$ in terms of $\wp (\alpha)$, $\wp '(\alpha)$ and $\kappa $, and we replace $\wp (\alpha)$, $\wp '(\alpha)$ and $\kappa $ with $\wp (C_3 \omega _1-C_1 \omega _3 )$, $\wp '(C_3 \omega _1-C_1 \omega _3 )$ and $\zeta (C_1 \omega _3 -C_3 \omega _1) +C_3 \eta _1 -C_1 \eta _3$. Then $\delta _1$ satisfies the sixth Painlev\'e equation in the elliptic form
\begin{equation}
\frac{d^2 \delta _1}{d \tau ^2} = -\frac{1}{8\pi ^2} \left\{ \sum _{i=0}^3 (l_i +1/2)^2 \wp '( \delta _1 + \omega _i) \right\}. \label{eq:P6ellipl}
\end{equation}
\end{prop}
We observe the expressions of $b_1$ and $\mu _1$ in detail for the cases $l_0=l_1=l_2=l_3=0$ and $l_0=1$, $l_1=l_2=l_3=0$.

\subsection{The case $M=1$, $r_1 =1$, $l_0=l_1=l_2=l_3=0$}

We investigate the case $M=1$, $r_1 =1$, $l_0=l_1=l_2=l_3=0$ in detail.
The differential equation (\ref{Hgkr1}) is written as 
\begin{equation}
\left\{ -\frac{d^2}{dx^2} + \frac{\wp ' (x)}{\wp (x) -b_1} \frac{d}{dx} - \frac{ \mu _{1} (4 b_1^3 -g_2 b_1 -g_3)}{\wp (x) -b_1} -p \right\} f_g (x)=0.
\label{Hgkr1l00}
\end{equation}
We assume that $b_1 \neq e_1, e_2, e_3$. The condition that the regular singular points $x= \pm \delta _1$ $(\wp (\delta _{1})=b_1 )$ are apparent is written as
\begin{align}
& p=- (4 b_1^3 -g_2 b_1 -g_3) \mu _{1} ^2 +(6b_1^2 -g_2/2) \mu _{1}, \label{pgkrapl00} 
\end{align}
(see Eq.(\ref{pgkrap})). The doubly-periodic function $\Xi (x)$ (see Eq.(\ref{Fxkr1})) which satisfies Eq.(\ref{prodDE}) is calculated as 
\begin{equation}
\Xi (x)= 2\mu _1 +\frac{1}{\wp(x)-b_1} .
\end{equation}
The value $Q$ (see Eq.(\ref{const})) is calculated as
\begin{align}
& Q= 2\mu _1(2\mu _1 (e_1-b_1)+1)(2(e_2-b_1) \mu _1+1)(2\mu _1(e_3-b_1)+1). 
\end{align}
We set
\begin{equation}
\Lambda _g(x) = \sqrt{\Xi (x) (\wp (x) - b_1)} \exp \int \frac{ \sqrt{-Q}dx}{\Xi (x)},
\label{integ1P6l00}
\end{equation}
(see Eq.(\ref{integ1P6})).
Then a solution to Eq.(\ref{Hgkr1l00}) is written as $\Lambda _g (x)$, and is expressed in the form of the Hermite-Krichever Ansatz as
\begin{align}
& \Lambda _g (x) = \bar{b} ^{(0)}_0 \exp (\kappa x) \Phi _0 (x, \alpha )
\end{align}
for generic $(\mu_1 , b_1)$.
The values $\alpha $ and $\kappa $ are determined as
\begin{align}
& \wp (\alpha )= b_1 - \frac{1}{2 \mu_1}, \quad \wp '(\alpha )= -\frac{\sqrt{-Q}}{2\mu _1^2} , \quad \kappa = \frac{\sqrt{-Q}}{2\mu _1}.
\end{align}
Hence we have
\begin{align}
& \mu _1 = -\frac{\kappa  }{\wp ' (\alpha )} ,\quad  b_1 = \wp (\alpha ) -\frac{\wp ' (\alpha )}{2\kappa }.
\end{align}
From Proposition \ref{prop:P6}, the function $\delta _1$ determined by
\begin{align}
\wp (\delta _1) =  b_1 & = \wp (C_3 \omega _1-C_1 \omega _3 ) -\frac{\wp ' (C_3 \omega _1-C_1 \omega _3 )}{2(\zeta (C_1 \omega _3 -C_3 \omega_1 ) -C_1 \eta _3 +C_3 \eta _1) } \label{P6sol0000} \\
& = \wp (C_1 \omega _3 -C_3 \omega_1 ) +\frac{\wp ' (C_1 \omega _3 -C_3 \omega_1 )}{2(\zeta (C_1 \omega _3 -C_3 \omega_1 ) -(C_1 \eta _3 -C_3 \eta _1)) } \nonumber
\end{align}
is a solution to the sixth Painlev\'e equation in the elliptic form (see Eq.(\ref{eq:P6ellipl})). This solution coincides with the one found by Hitchin \cite{Hit} when he studied Einstein metrics and isomonodromy deformations.

Note that in \cite{DIKZ,KK}, solutions in terms of theta functions are obtained.

Now we consider the case $Q=0$. If $Q=0$, then $\mu _1=0$ or $\mu _1 =1/(2(b_1 -e_i))$ for some $i \in \{1,2,3 \}$. 

If $\mu _1 =0$, then a solution to Eq.(\ref{Hgkr1l00}) is $1 (= \Lambda _g (x))$ and another solution is written as 
\begin{equation}
\textstyle \zeta (x) +b_1 x (=\int -(\wp (x) -b_1) dx ) .
\end{equation} 
We investigate the monodromy preserving deformation on the basis $s _1 (x)= B(\tau )$ and $s_2(x )=  \zeta (x) +b_1 x$, where $B(\tau )$ is a constant that is independent of $x$.
The monodromy matrix with respect to the path $\gamma _0$ is written as diag$(1,-1)$.
Since $s_2 (x + 2\omega _j) = s_2 (x) +2(\eta _j + \omega _j b_1)$ $(j=1,3)$, the monodromy matrix with respect to the basis $(s_1 (x), s _2 (x))$ on the cycle $x \rightarrow x+2\omega _j$ $(j=1,3)$ is written as 
\begin{equation}
\left( 
\begin{array}{cc}
1 & 2(\eta _j + \omega _j b_1)/B(\tau ) \\
0 & 1
\end{array}
\right).
\end{equation}
To preserve monodromy, the matrix elements should be constants of the variable $\tau (=\omega _3/\omega _1)$ up to simultaneous change of basis. Hence we obtain 
\begin{align}
& 2(\eta _1 + \omega _1 b_1) =D_1 B(\tau ) ,\\
& 2(\eta _3 + \omega _3 b_1) =D_3 B(\tau ) , \nonumber
\end{align}
for some constants $D_1$ and $D_3$. By using Legendre's relation, we obtain that $ B(\tau ) = \pi \sqrt{-1} /(D_1 \omega _3 -D_3 \omega _1)$ and 
\begin{equation}
\wp (\delta _1) =b_1 =-\frac{D_1 \eta _3 -D_3 \eta _1}{D_1 \omega _3 -D_3 \omega _1}.
\label{b1mu10}
\end{equation}
Since Eq.(\ref{b1mu10}) is obtained by monodromy preserving deformation, the function $\delta _1$ satisfies the sixth Painlev\'e equation.

If $\mu _1 =1/(2(b_1 -e_i))$  for some $i \in \{1,2,3 \}$, then $\wp _i (x) (= \Lambda _g (x))$ is a solution to Eq.(\ref{Hgkr1l00}), and another solution is written as
\begin{equation}
\textstyle \wp _i (x) \left\{ \frac{e_i- b_1}{(e_i-e_{i'})(e_i-e_{i''})} \zeta (x +\omega _i)  +\right. \left. (1- \frac{e_i- b_1}{(e_i-e_{i'})(e_i-e_{i''})}) x \right\} (= \wp _i (x) \int \frac{\wp (x) -b}{\wp (x)-e_i} dx ),
\end{equation}
where $i'$ and $i''$ are elements in $\{ 1,2,3 \}$ such that $i' \neq i$, $i'' \neq i$ and $i' <i''$.
By calculating similarly to the case $\mu _1=0$, we obtain that the function $\delta _1$, which is determined by
\begin{equation}
\wp (\delta _1) =b_1 =\frac{(g_2/4 -2e_i^2)(D_1 \omega _3 -D_3 \omega _1) +e_i (D_1 \eta _3 -D_3 \eta _1) }{e_i (D_1 \omega _3 -D_3 \omega _1) +(D_1 \eta _3 -D_3 \eta _1) },
\label{b1mui}
\end{equation}
is a solution to the sixth Painlev\'e equation for constants $D_1$ and $D_3$.

We now show that Eqs.(\ref{b1mu10} ,\ref{b1mui}) are obtained by suitable limits from Eq.(\ref{P6sol0000}). Set $(C_1 ,C_3 )= (CD_1, CD_3)$ in Eq.(\ref{P6sol0000}) and consider the limit $C \rightarrow 0$, then we recover Eq.(\ref{b1mu10}). Similarly, set $(C_1 ,C_3) =(C D_1 ,-1+CD_3 )$ (resp. $(C_1 ,C_3) =(-1 +C D_1 ,1+ CD_3 )$, $(C_1 ,C_3) =(1+C D_1 , CD_3 )$) and consider the limit $C \rightarrow 0$, then we recover Eq.(\ref{b1mui}) for the case $i=1$ (resp. $i=2$, $i=3$).
Hence the space of the parameters of the solutions to the sixth Painlev\'e equation (i.e. the space of initial conditions) for the case $l_0=l_1=l_2=l_3=0$ is obtained by blowing up four points on the surface $\Cplx /2 \Zint \times \Cplx /2 \Zint $, and this reflects the $A_1 \times A_1 \times A_1 \times A_1$ structure of Riccati solutions by Saito and Terajima \cite{ST}.

\subsection{The case $M=1$, $r_1 =1$, $l_0=1$, $l_1=l_2=l_3=0$}
The differential equation (\ref{Hgkr1}) for this case is written as 
\begin{equation}
\left\{ -\frac{d^2}{dx^2} + \frac{\wp ' (x)}{\wp (x) -b_1} \frac{d}{dx} - \frac{ \mu _{1} (4 b_1^3 -g_2 b_1 -g_3)}{\wp (x) -b_1} + 2\wp (x) -p \right\} f_g (x)=0,
\label{Hgkr1l01}
\end{equation}
We assume that $b_1 \neq e_1, e_2, e_3$. The condition that the regular singular points $x= \pm \delta _1$ $(\wp (\delta _{1})=b_1 )$ are apparent is written as
\begin{align}
& p=- (4 b_1^3 -g_2 b_1 -g_3) \mu _{1} ^2 +(6b_1^2 -g_2/2) \mu _{1} +2b_1, \label{pgkrapl01} 
\end{align}
(see Eq.(\ref{pgkrap})). The doubly-periodic function $\Xi (x)$ (see Eq.(\ref{Fxkr1})), which satisfies Eq.(\ref{prodDE}), is calculated as 
\begin{align}
\Xi (x)=& \wp (x) +( (-4b_1^3+b_1g_2+g_3)\mu _1^2+(6b_1^2-g_2/2)\mu _1 -b_1) \\
& +  ((-4b_1^3+b_1g_2+g_3)\mu _1 /2+3b_1^2-g_2/4)/(\wp (x) -b_1) .\nonumber
\end{align}
The value $Q$ (see Eq.(\ref{const})) is calculated as
\begin{align}
 Q= -& ((2(4 b_1^3-b_1g_2-g_3)\mu _1^3-(12b_1^2-g_2)\mu _1^2+4) (2(b_1^2+e_1b_1+e_2e_3)\mu _1-2b_1-e_1) \\
& (2(b_1^2+e_2b_1+e_1e_2)\mu _1-2b_1-e_2) (2(b_1^2+e_3b_1+e_1e_3)\mu _1-2b_1-e_3) \nonumber .
\end{align}
We set
\begin{equation}
\Lambda _g(x) = \sqrt{\Xi (x) (\wp (x) - b_1)} \exp \int \frac{ \sqrt{-Q}dx}{\Xi (x)}, 
\label{integ1P6l01}
\end{equation}
(see Eq.(\ref{integ1P6})).
Then a solution to Eq.(\ref{Hgkr1l01}) is written as $\Lambda _g (x)$, and it is expressed in the form of the Hermite-Krichever Ansatz as
\begin{align}
& \Lambda _g (x) = \exp (\kappa x) \left\{ \bar{b} ^{(0)}_0 \Phi _0 (x, \alpha ) +\bar{b} ^{(0)}_1 \frac{d}{dx} \Phi _0 (x, \alpha ) \right\} 
\end{align}
for generic $(\mu_1 , b_1)$.
The values $\alpha $ and $\kappa $ are determined as
\begin{align}
& \wp (\alpha )  = \frac{2(4b_1^3-b_1g_2-g_3)b_1 \mu _1^3+(-24b_1^3+4g_2b_1+3g_3)\mu _1^2+(24b_1^2-2g_2)\mu _1-8b_1}{2(4b_1^3-b_1 g_2-g_3)\mu _1^3 -(12b_1^2-g_2)\mu _1^2+4},\\
& \wp '(\alpha )= \frac{-4((4b_1^3-b_1g_2-g_3)\mu _1^3-(12b_1^2-g_2)\mu _1^2+12b_1\mu _1-4)}{(2(4b_1^3-b_1 g_2-g_3)\mu _1^3 -(12b_1^2-g_2)\mu _1^2+4)^2}\sqrt{-Q}  ,\\
& \kappa = \frac{2\mu _1}{2(4b_1^3-b_1 g_2-g_3)\mu _1^3 -(12b_1^2-g_2)\mu _1^2+4}\sqrt{-Q}.
\end{align}
Hence we have
\begin{align}
& b_1 = \frac{2\wp (\alpha ) \kappa ^3-3\wp '(\alpha )\kappa ^2+(6\wp (\alpha ) ^2 -g_2)\kappa -\wp (\alpha ) \wp '(\alpha )}{2(\kappa ^3-3\wp (\alpha ) \kappa +\wp '(\alpha ))} ,\\
& \mu _1 = \frac{2(\kappa ^3-3\wp (\alpha ) \kappa +\wp '(\alpha ))\kappa }{-2\wp '(\alpha )\kappa ^3+(12\wp (\alpha ) ^2-g_2)\kappa ^2-6\wp (\alpha ) \wp '(\alpha )\kappa +\wp '(\alpha )^2}.
\end{align}
From Proposition \ref{prop:P6}, the function $\delta _1$ determined by
\begin{align}
& \wp (\delta _1) = b_1 = \label{P6sol1000} \\ 
& \quad  \frac{2\wp (\omega ) (\zeta (\omega )- \eta )^3+3\wp '(\omega )(\zeta (\omega )- \eta )^2+(6\wp (\omega ) ^2 -g_2)(\zeta (\omega )- \eta )+\wp (\omega ) \wp '(\omega )}{2((\zeta (\omega )- \eta )^3-3\wp (\omega ) (\zeta (\omega )- \eta ) -\wp '(\omega ))} , \nonumber \\
& (\omega = C_1 \omega _3 -C_3 \omega_1 , \quad  \eta = C_1 \eta _3 -C_3 \eta _1) ,\nonumber
\end{align}
is a solution to the sixth Painlev\'e equation in the elliptic form (see Eq.(\ref{eq:P6ellipl})).
In the sixth Painlev\'e equation, it is known that the case $(\kappa _{0}, \kappa _{1}, \kappa _{t}, \kappa _{\infty}) =(1/2, 1/2, 1/2, 3/2 ) $ is linked to the case $(\kappa _{0}, \kappa _{1}, \kappa _{t}, \kappa _{\infty}) =(1/2, 1/2, 1/2, 1/2 ) $ by B\"acklund transformation. 
For a table of B\"acklund transformation of the sixth Painlev\'e equation, see \cite{TOS}.
By transformating the solution in Eq.(\ref{P6sol0000}) of the case $(\kappa _{0}, \kappa _{1}, \kappa _{t}, \kappa _{\infty}) =(1/2, 1/2, 1/2, 1/2 )$ to the one of the case  $(\kappa _{0}, \kappa _{1}, \kappa _{t}, \kappa _{\infty}) =(1/2, 1/2, 1/2, 3/2 )$, we recover the solution in Eq.(\ref{P6sol1000}).

Now we consider the case $Q=0$. If $Q=0$, then $\mu_1  $ is a solution to the equation $2(4 b_1^3-b_1g_2-g_3)\mu _1^3-(12b_1^2-g_2)\mu _1^2+4 =0$ or $\mu _1=(2b_1+e_i)/(2(b_1^2+e_ib_1+e_i^2-g_2/4))$ for some $i \in \{1,2,3 \}$.
We set $\omega = D_1 \omega _3 -D_3 \omega_1$ and $\eta = D_1 \eta _3 -D_3 \eta _1$, where $D_1$ and $D_3$ are constants. For the case that $\mu_1  $ is a solution to the equation $2(4 b_1^3-b_1g_2-g_3)\mu _1^3-(12b_1^2-g_2)\mu _1^2+4 =0$, the corresponding solutions to the sixth Painlev\'e equation are written as the function $\delta _1$, where 
\begin{equation}
\wp( \delta _1)= b_1 = \frac{4 \eta ^3 +g_2 \omega ^2 \eta -2 g_3 \omega ^3}{\omega (g_2 \omega ^2-12\eta ^2 )}. \label{b1mu1l01}
\end{equation}
For the case $\mu _1=(2b_1+e_i)/(2(b_1^2+e_ib_1+e_i^2-g_2/4))$ ($i \in \{1,2,3 \}$), we have 
\begin{equation}
\wp( \delta _1)= b_1 = \frac{-g_2 e_i \omega /2 +(6e_i^2 -g_2)\eta}{(6e_i^2 -g_2)\omega -6e_i\eta }. \label{b1muil01}
\end{equation}
Note that these solutions are also obtained by suitable limits from Eq.(\ref{P6sol1000}), and Eq.(\ref{b1mu1l01}) (resp. Eq.(\ref{b1muil01})) is transformed by B\"acklund transformation from Eq.(\ref{b1mu10}) (resp. Eq.(\ref{b1mui})).

\section{Relationship with finite-gap potential} \label{sec:FGP}

\subsection{Finite-gap property} \label{sec:fg}

We investigate the condition that the potential in Eq.(\ref{Inopotent}) is finite-gap.

If $M=0$ ($M$ is the number of additional apparent singularities) and $l_0, l_1, l_2, l_3 \in \Zint _{\geq 0}$, then the potential is called the Treibich-Verdier potential, and it is algebro-geometric finite-gap.

Next we consider the case $M=1$ and $r_1=2$.
Set $b_1 =\wp (\delta _1)$. The condition that the regular singularity $x =\pm \delta _{1} $ of Eq.(\ref{eq:H}) is apparent (which is equivalent to that the regular singularity $z= b_1$ of Eq.(\ref{Ino}) is apparent) is written as 
\begin{align}
& s_1^3+ (12b_1^2- g_2)s_1^2  +(4(4b_1^3-g_2 b_1 -g_3)E+ f_1(b_1))s_1  +f_0(b_1) =0.
\end{align}
where $f_1(b_1) $ and $f_0(b_1)$ are given by
\begin{align}
f_1(b_1) =& -2(2l_0^2+2l_0+5)b_1 (4b_1^3-g_2 b_1 -g_3) +(6b_1 ^2 -g_2/2)^2 \\
& -8(2l_1^2+2l_1+1)(b_1-e_2)(b_1-e_3)(e_1 b_1+e_1^2+e_2e_3 ) \nonumber \\
& -8(2l_2^2+2l_2+1)(b_1-e_1)(b_1-e_3)(e_2 b_1+e_2^2+e_1e_3 ) \nonumber \\
& -8(2l_3^2+2l_3+1)(b_1-e_1)(b_1-e_2)(e_3 b_1+e_3^2+e_1e_2 ) ,\nonumber \\
f_0(b_1) =& (2l_0+1)^2(4b_1^3-g_2 b_1 -g_3)^2 \\
& - 16(2l_1+1)^2 (e_1-e_2)(e_1-e_3)(b_1-e_2)^2(b_1-e_3)^2 \nonumber \\
& - 16(2l_2+1)^2 (e_2-e_1)(e_2-e_3)(b_1-e_1)^2(b_1-e_3)^2  \nonumber \\
& - 16(2l_3+1)^2 (e_3-e_1)(e_3-e_2)(b_1-e_1)^2(b_1-e_2)^2 . \nonumber 
\end{align}
If $s_1 =0$, then we obtain an equation
\begin{align}
& f_0(b_1) =0.
\label{eq:f0b10}
\end{align}
Remarkably, the value $b_1$ determined by this equation does not depend on the value $E$.
It is shown by Treibich \cite{Tre} that, if $M=1$, $r_1=2$, $s_1=0$, $l_0, l_1, l_2, l_3 \in \Zint _{\geq 0}$ and $b_1 (=\wp (\delta _1))$ satisfies Eq.(\ref{eq:f0b10}) (which is equivalent to Eq.(\ref{eq:Trepot}) by setting $\delta _1 =\delta$), then the potential is algebro-geometric finite-gap.

In this section, we investigate the differential equation
\begin{align}
&  \left( -\frac{d^2}{dx^2} + v(x)\right) f(x) =Ef(x), \label{fingapDE} \\
v(x) = & \sum_{i=0}^3 l_i(l_i+1) \wp (x+\omega_i) +2\sum_{i'=1}^M (\wp (x-\delta _{i'}) + \wp (x+\delta _{i'})) , \nonumber
\end{align}
for the case when the regular singular points $x =\pm \delta _{i'}$ ($i'=1,\dots ,M$) of Eq.(\ref{fingapDE}) are apparent, $\delta _j \not \equiv \omega _i$ mod $2\omega_1 \Zint \oplus 2\omega_3 \Zint$ $(0\leq i\leq 3, \; 1\leq j\leq M)$ and  $\delta _j \pm \delta _{j'} \not \equiv 0$  mod $2\omega_1 \Zint \oplus 2\omega_3 \Zint$ $(1\leq j< j' \leq M)$. Note that the potential in Eq.(\ref{fingapDE}) corresponds to the one in Eq.(\ref{Inopotent}) with conditions $r_{i'}=2$ and $s_{i'}=0$ $(i'=1,\dots ,M)$.
The exponents of Eq.(\ref{fingapDE}) at the regular singularity $x =\pm \delta _{i'}$ $(i'=1,\dots ,M)$ are $-1$ and $2$.
\begin{prop}
Assume that $l_0 , l_1, l_2 , l_3 \in \Zint _{\geq 0}$, $\delta _j \not \equiv \omega _i$ mod $2\omega_1 \Zint \oplus 2\omega_3 \Zint$ $(0\leq i\leq 3, \; 1\leq j\leq M)$ and  $\delta _j \pm \delta _{j'} \not \equiv 0$  mod $2\omega_1 \Zint \oplus 2\omega_3 \Zint$ $(1\leq j< j' \leq M)$. If the values $\delta _1 ,\dots, \delta _M$ satisfy the equation
\begin{equation}
2\sum _{j' \neq j} (\wp ' (\delta _j -\delta _{j'} ) + \wp ' (\delta _j +\delta _{j'} ) )+ \sum _{i=0}^3 (l_i +1/2)^2 \wp ' (\delta _j +\omega _i ) =0 \quad (j=1,\dots ,M),
\label{eq:ds}
\end{equation}
then the regular singular points $x =\pm \delta _{i'}$ ($i'=1,\dots ,M$) of Eq.(\ref{fingapDE}) are apparent.
\end{prop}
\begin{proof}
We show that, if $\delta _1 ,\dots, \delta _M$ satisfy Eq.(\ref{eq:ds}), then the regular singular points $x = \delta _j$ $(j=1,\dots ,M)$ are apparent.
The coefficients of expansion around $x=\delta _j$ of Eq.(\ref{fingapDE}) written in the form of Eq.(\ref{eq:Feqxa}) ($a= \delta _j$) are given by
\begin{align}
& p_k=0 \quad (k\in \Zint _{\geq 0}), \quad q_0=-2, \quad q_1=0, \\
& q_2=E - \left( 2\wp ( 2\delta _j) + 2\sum _{j' \neq j} (\wp (\delta _j -\delta _{j'} ) + \wp (\delta _j +\delta _{j'} ) + \sum _{i=0}^3 l_i(l_i +1) \wp (\delta _j +\omega _i )\right) , \nonumber \\
& q_3= -\left(2 \wp '( 2\delta _j) + 2\sum _{j' \neq j} (\wp ' (\delta _j -\delta _{j'} ) + \wp ' (\delta _j +\delta _{j'} ) )+ \sum _{i=0}^3 l_i(l_i +1) \wp ' (\delta _j +\omega _i )\right) ,\nonumber
\end{align}
and the characteristic polynomial $F(t)$ at $x =\delta _{j}$ is written as $F(t)=(t-2)(t+1)$. 
In section \ref{sec:FDE}, we obtained a condition for apparency of a regular singular point. On the case $x =\delta _{j}$, it is written as
\begin{equation}
(-p_3+q_3)c_0+q_2c_1 +(p_1+q_1)c_2 =0,
\end{equation} 
with $c_0=1$, $-2c_1+(-p_1+q_1)c_0=0$, $-2c_2 +(-p_2+ q_2)c_0+ q_1 c_1=0$ (see Eqs.(\ref{eq:recn}, \ref{eq:recj})). Hence $c_1=0$ and the condition that $x=\delta _j$ is apparent is written as $q_3=0$, i.e.,
\begin{equation}
2 \wp '( 2\delta _j) + 2\sum _{j' \neq j} (\wp ' (\delta _j -\delta _{j'} ) + \wp ' (\delta _j +\delta _{j'} ) )+ \sum _{i=0}^3 l_i(l_i +1) \wp ' (\delta _j +\omega _i ) =0. \label{eq:ds2}
\end{equation}
From the identity
\begin{equation}
 8 \wp '( 2x) = \wp ' (x) + \wp ' (x +\omega _1) + \wp ' (x +\omega _2) + \wp ' (x +\omega _3),
\end{equation}
Eq.(\ref{eq:ds2}) is equivalent to Eq.(\ref{eq:ds}).
The condition that $x = -\delta _j$ is apparent is given by
\begin{equation}
2 \wp '( -2\delta _j) + 2\sum _{j' \neq j} (\wp ' (-\delta _j -\delta _{j'} ) + \wp ' (-\delta _j +\delta _{j'} ) )+ \sum _{i=0}^3 l_i(l_i +1) \wp ' (-\delta _j +\omega _i ) =0, \label{eq:ds2'}
\end{equation}
and it is equivalent to Eq.(\ref{eq:ds}) by the oddness and the double-periodicity of the function $\wp '(x)$. Therefore, if Eq.(\ref{eq:ds}) is satisfied, then the points $x=\pm \delta _j$ $(j=1,\dots ,M)$ are apparent.
\end{proof}
It is remarkable that Eq.(\ref{eq:ds}) does not contain the variable $E$. We examine this equation with the introduction of
\begin{equation}
\Phi (\delta _1, \dots ,\delta _M) = 2\sum _{1\leq j_1 < j_2 \leq M} (\wp (\delta _{j_1} -\delta _{j_2} ) + \wp (\delta _{j_1} +\delta _{j_2} ) )+ \sum _{j=1}^M \sum _{i=0}^3 (l_i +1/2)^2 \wp (\delta _j +\omega _i ),
\end{equation}
in which case Eq.(\ref{eq:ds}) is equivalent to the equations
\begin{equation}
\frac{\partial }{\partial \delta _j} \Phi (\delta _1, \dots ,\delta _M) =0 \quad (j=1,\dots ,M).
\label{eq:dsPhi}
\end{equation}

We will now show that Eq.(\ref{eq:ds}) has a good solution.
\begin{prop}
Assume that $l_0 , l_1 \in \Rea \setminus \{-1/2\}$, $l_2 , l_3 \in \Rea $, $\omega _1 \in \Rea _{>0 }$ and $\omega _3 \in \sqrt{-1} \Rea _{>0 }$. Then Eq.(\ref{eq:ds}) has a solution such that $\delta _j \in \Rea$, $\delta _j \not \equiv \omega _i$ mod $2\omega_1 \Zint \oplus 2\omega_3 \Zint$ $(0\leq i\leq 3, \; 1\leq j\leq M)$ and  $\delta _j \pm \delta _{j'} \not \equiv 0$  mod $2\omega_1 \Zint \oplus 2\omega_3 \Zint$ $(1\leq j< j' \leq M)$.
\end{prop}
\begin{proof}
From the assumption $\omega _1 \in \Rea _{>0 }$ and $\omega _3 \in \sqrt{-1} \Rea _{>0 }$, the functions $\wp (x+\omega _i)$ $(i=0,1,2,3)$ are real-valued for $x \in \Rea \setminus \omega _1 \Zint $ and $\lim _{x \rightarrow 0 , x \in \Rea} \wp (x) = \lim _{x \rightarrow \omega _1 , x \in \Rea} \wp (x+\omega _1) = +\infty$.
Now we consider the function $\Phi (\delta _1, \dots ,\delta _M)$ on the real domain $D= \{ (\delta _1 ,\dots , \delta _M) \in \Rea ^M | \: 0< \delta _1 <\dots <\delta _M <\omega _1 , \; \delta _j +\delta _{j'} <\omega _1 \: (\forall j,j' $ s.t. $j<j'  )\}$. Then $\Phi (\delta _1, \dots ,\delta _M)$ is real-valued and continuous on the domain $D$.
As $(\delta _1, \dots ,\delta _M)$ tends to the boundary of the domain $D$, the value $\Phi (\delta _1, \dots ,\delta _M)$ tends to $+\infty$ by the assumption $l_0 \neq -1/2$, $l_1 \neq -1/2$.
Since $\wp (x+\omega _i) \geq $min$ (e_1,e_2,e_3)=e_3$ for $x\in \Rea $ and $i=0,1,2,3$, we have $\Phi (\delta _1, \dots ,\delta _M) \geq (2M-2+ \sum _{i=0}^3 (l_i+1/2)^2)Me_3 $.
Therefore the function $\Phi (\delta _1, \dots ,\delta _M)$ has a minimum value at $\exists (\delta _1^0, \dots ,\delta _M^0) \in D$.
Since $(\delta _1^0, \dots ,\delta _M^0)$ is an extremal point of the function $\Phi (\delta _1, \dots ,\delta _M)$, it satisfies Eq.(\ref{eq:dsPhi}).
Hence $(\delta _1^0, \dots ,\delta _M^0)$ is a solution to Eq.(\ref{eq:ds}). 

Because $(\delta _1^0, \dots ,\delta _M^0) \in D$, it satisfies $\delta _j^0 \not \equiv \omega _i$ mod $2\omega_1 \Zint \oplus 2\omega_3 \Zint$ $(0\leq i\leq 3, \; 1\leq j\leq M)$ and  $\delta _j^0 \pm \delta _{j'}^0 \not \equiv 0$  mod $2\omega_1 \Zint \oplus 2\omega_3 \Zint$ $(1\leq j< j' \leq M)$.
\end{proof}

Upon introducing $b_j=\wp (\delta _j)$ $(j=1 ,\dots ,M)$, it follows from the relations
\begin{align}
& \wp '(x+y) +\wp '(x-y) =-\frac{\wp '(x) \wp ''(y) }{(\wp (x)-\wp (y))^2}-\frac{2\wp '(x) \wp '(y)^2 }{(\wp (x)-\wp (y))^3}, \\
& \wp '(x +\omega _i) =-\frac{3e_i^2-g_2/4}{(\wp (x)-e_i)^2} \wp '(x), \quad (i=1,2,3),
\end{align}
that Eq.(\ref{eq:ds}) may be expressed in the algebraic form 
\begin{align}
& \sum _{j' \neq j}  \left\{ \frac{12b_{j'}^2-g_2}{(b_j-b_{j'})^2} + \frac{4(4b_{j'}^3-g_2 b_{j'}-g_3)}{(b_j-b_{j'})^3} \right\} \\
& \quad \quad = (l_0 +1/2)^2 - \sum _{i=1}^3 (l_i +1/2)^2 \frac{3e_i^2-g_2/4}{(b_j-e_i)^2}, \quad (j=1,\dots ,M), \nonumber
\end{align}
under the condition $b_j \neq e_1, e_2, e_3$ $(j=1,\dots ,M)$. For the case $M=1$, it is written as Eq.(\ref{eq:f0b10}).

If $\delta _1, \dots ,\delta _M$ satisfy Eq.(\ref{eq:ds}), then the regular singular points $x= \pm \delta _{i'}$ $(i'=1,\dots ,M)$ of Eq.(\ref{fingapDE}) are apparent for all $E$, and it follows from Proposition \ref{prop:prod} that there exists a non-zero solution to the third-order differential equation satisfied by products of two solutions to Eq.(\ref{fingapDE}). Namely, if $l_0, l_1, l_2, l_3 \in \Zint _{\geq 0}$, $ \delta _1 ,\dots ,\delta _M$ satisfy Eq.(\ref{eq:ds}), $\delta _j \not \equiv \omega _i$ mod $2\omega_1 \Zint \oplus 2\omega_3 \Zint$ $(0\leq i\leq 3, \; 1\leq j\leq M)$ and $\delta _j \pm \delta _{j'} \not \equiv 0$  mod $2\omega_1 \Zint \oplus 2\omega_3 \Zint$ $(1\leq j< j' \leq M)$, then Eq.(\ref{prodDE})
has an even non-zero doubly-periodic solution that has the expansion
\begin{equation}
\Xi (x)=c_0+\sum_{i=0}^3 \sum_{j=0}^{l_i-1} b^{(i)}_j \wp (x+\omega_i)^{l_i-j} + \sum _{i'=1}^M \left(  \frac{d^{(i')}_0}{(\wp (x)-\wp (\delta _{i'}))^2} +\frac{d^{(i')}_1}{(\wp (x)-\wp (\delta _{i'}))}\right) ,
\label{Fx0}
\end{equation}
for all $E$.
On the present situation, we can improve Proposition \ref{prop:prod}.
\begin{prop} \label{prop:prodFG} $ $\\
(i) For each $l_0, l_1, l_2, l_3 \in \Zint _{\geq 0}$, periods $(2\omega _1, 2\omega _3)$ and values $ \delta _1 ,\dots ,\delta _M$, the number of eigenvalues $E$, such that the dimension of the space of even doubly-periodic solutions to Eq.(\ref{prodDE}) is no less than two, is finite.\\
(ii) If $l_0, l_1, l_2, l_3 \in \Zint _{\geq 0}$, $\delta _1 ,\dots ,\delta _M $ satisfy Eq.(\ref{eq:ds}), $\delta _j \not \equiv \omega _i$ mod $2\omega_1 \Zint \oplus 2\omega_3 \Zint$ $(0\leq i\leq 3, \; 1\leq j\leq M)$ and $\delta _j \pm \delta _{j'} \not \equiv 0$ mod $2\omega_1 \Zint \oplus 2\omega_3 \Zint$ $(1\leq j< j' \leq M)$, then Eq.(\ref{prodDE}) has a unique non-zero doubly-periodic solution $\Xi (x,E)$, which has the expansion
\begin{equation}
\Xi (x,E)=c_0(E)+\sum_{i=0}^3 \sum_{j=0}^{l_i-1} b^{(i)}_j (E)\wp (x+\omega_i)^{l_i-j}+ \sum _{i'=1}^M d^{(i')}(E)(\wp (x+\delta _{i'}) +\wp (x- \delta _{i'})),
\label{Fx1}
\end{equation}
where the coefficients $c_0(E)$, $b^{(i)}_j(E)$ and $d^{(i')}(E)$ are polynomials in $E$ such that these polynomials do not share any common divisors and the polynomial $c_0(E)$ is monic. 
We set $g=\deg_E c_0(E)$. Then the coefficients satisfy $\deg _E b^{(i)}_j(E)<g$ for all $i$ and $j$, and $\deg _E d^{(i')}(E)<g$ for all $i'$.
\end{prop}
\begin{proof}
By substituting Eq.(\ref{Fx0}) into Eq.(\ref{prodDE}), we derive linear equations in coefficients $c_0$, $b^{(i)}_j$, $d^{(i')}_0$ and $d^{(i')}_1$ to satisfy Eq.(\ref{prodDE}). We replace $c_0$, $b^{(i)}_j$, $d^{(i')}_0$, $d^{(i')}_1$ with $\tilde{c} _1, \dots ,\tilde{c} _{M'}$ $(M'=1+l_0+l_1+l_2+l_3+2M)$. Then the linear equations are written as 
\begin{equation}
\sum _{i=1}^{M'} (m_{k,i} E+n_{k,i}) \tilde{c}_i =0, \quad (k=1, \dots , M''),
\label{eq:Xirel}
\end{equation}
where $M''$ is the number of equations.
It follows from Proposition \ref{prop:prod} that there exists a non-zero solution to Eq.(\ref{prodDE}). Hence all minors of the matrix $(m_{k,i} E+n_{k,i}) _{k,i}$ of rank $M'$ are identically zero.

Now we assume that there exists infinitely-many values $E$ such that the dimension of the space of even doubly-periodic solutions to Eq.(\ref{prodDE}) is no less than two. Since any minors of the matrix $(m_{k,i} E+n_{k,i}) _{k,i}$ of rank $M'-1$ are written as polynomials in $E$, and they must be zero at infinitely-many values of $E$ by the assumption, they are identically zero. Hence the dimension of the space of even doubly-periodic solutions to Eq.(\ref{prodDE})  is no less than two for all $E$.
Because the coefficients of Eq.(\ref{eq:Xirel}) are written as polynomials in $E$, there exist linearly independent functions $\Xi ^{(1)}(x,E)$ and $\Xi ^{(2)} (x,E)$  which solve Eq.(\ref{prodDE}) and may be expressed in the form 
\begin{equation}
\Xi ^{(k)} (x,E) = \sum _{j=0}^{g_k} a_j ^{(k)}(x) E^{{g_k} -j} , \;  a_0^{(k)} (x) \neq 0, \quad (k=1,2).
\end{equation}
\begin{lemma} \label{lemma:Xi}
If $\tilde{\Xi }(x)$ is a solution of Eq.(\ref{prodDE}) written in the form $\tilde{\Xi }(x) = \sum _{i=0}^{\tilde{g}} \tilde{a}_i (x) E^{\tilde{g} -i}$ $(\tilde{a}_0 (x) \neq 0)$, then $\tilde{a}_0 (x)$ is independent of $x$.
\end{lemma}
\begin{proof}
By substituting $\tilde{\Xi }(x)$ into Eq.(\ref{prodDE}) and considering the coefficients of $E^{\tilde{g} +1}$, we obtain that $\tilde{a}'_0 (x) =0$. Hence $\tilde{a}_0 (x)$ is independent of $x$.
\end{proof}  
By the lemma, $a_0 ^{(0)}(x)$ and $a_0 ^{(1)}(x)$ are independent of $x$, and so we may denote them by $a_0 ^{(0)}$ and $a_0 ^{(1)}$. 
If $g_1 \geq g_2$ (resp. $g_1<g_2$), then we set $\tilde{\Xi }^ {(1)} (x,E) =  \Xi ^{(1)} (x,E) - (a_0^{(1)} /a_0^{(2)} )\Xi ^{(2)} (x,E) E^{g_1-g_2}$, $\tilde{\Xi }^ {(2)} (x,E) =  \Xi ^{(2)} (x,E)$ (resp. $\tilde{\Xi }^ {(1)} (x,E) =  \Xi ^{(1)} (x,E)$, $\tilde{\Xi }^ {(2)} (x,E) =  \Xi ^{(2)} (x,E) - (a_0^{(2)} /a_0^{(1)} )\Xi ^{(1)} (x,E) E^{g_2-g_1}$). Then the degree of either $\tilde{\Xi } ^{(1)} (x,E)$ or $\tilde{\Xi } ^{(2)} (x,E)$ in $E$ decreases from the one of $\Xi ^{(1)} (x,E)$ or $\Xi ^{(2)} (x,E)$, and so $\tilde{\Xi } ^{(1)} (x,E)$ and $\tilde{\Xi } ^{(2)} (x,E)$ are linearly-independent solutions to Eq.(\ref{prodDE}).
From Lemma \ref{lemma:Xi}, the top terms of $\tilde{\Xi }^ {(1)} (x,E)$ and $\tilde{\Xi }^ {(2)} (x,E)$ in $E$ are non-zero constants.
By the same procedure, we can construct functions $\tilde{\tilde{\Xi }}^ {(1)} (x,E)$ and $\tilde{\tilde{\Xi }}^ {(2)} (x,E)$ which are linearly independent solutions to Eq.(\ref{prodDE}) and the degree of either $\tilde{\tilde{\Xi }} ^{(1)} (x,E)$ or $\tilde{\tilde{\Xi }} ^{(2)} (x,E)$ in $E$ decreases from the one of $\tilde{\Xi }^{(1)} (x,E)$ or $\tilde{\Xi }^{(2)} (x,E)$.
By repeating this decreasing procedure, we find that there exist linearly-independent solutions to Eq.(\ref{prodDE}) such that their degrees in $E$ are zero. This is a contradiction, because if a solution, $f(x)$, to Eq.(\ref{prodDE}) is independent of $E$, then $f'(x)=v'(x)=0$.
Therefore we have that the number of eigenvalues $E$, such that the dimension of the space of even doubly-periodic solutions to Eq.(\ref{prodDE}) is no less than two, is finite at most. 

From Proposition \ref{prop:prod}, there exist non-zero solutions to Eq.(\ref{prodDE}). Thus Eq.(\ref{eq:Xirel}) has a non-zero solution for all $E$. Because the coefficients in Eq.(\ref{eq:Xirel}) are polynomial in $E$, a solution to Eq.(\ref{eq:Xirel}) is written in terms of rational functions in $E$. By multiplying by an appropriate term, a solution to Eq.(\ref{eq:Xirel}) (i.e., $c_0$, $b^{(i)}_j$, $d^{(i')}_0$, $d^{(i')}_1$) may be expressed by polynomials in $E$ which do not share a common divisor, and they are determined uniquely up to scalar multiplication, because the dimension of solutions is one. We denote the doubly-periodic function uniquely determined in this way by $\Xi (x,E)$. By combining with the relation
\begin{align}
& \wp (x+\delta )+ \wp (x-\delta )=  2\wp (\delta ) +\frac{\wp''(\delta)}{\wp (x)-\wp (\delta )} +\frac{\wp'(\delta)^2}{(\wp (x)-\wp (\delta ))^2} ,
\end{align} 
the function $\Xi (x, E)$ is expressed as
\begin{align}
\Xi (x,E)= & c_0(E)+\sum_{i=0}^3 \sum_{j=0}^{l_i-1} b^{(i)}_j (E)\wp (x+\omega_i)^{l_i-j} \label{Fx1-1} \\
& + \sum _{i'=1}^M \left( d^{(i')}(E)(\wp (x+\delta _{i'}) +\wp (x- \delta _{i'}))+\frac{d^{(i')}_1(E)}{(\wp (x)-\wp (\delta _{i'}))}\right), \nonumber
\end{align}
where $c_0(E)$, $b^{(i)}_j (E)$, $d^{(i')}(E)$ and $d^{(i')}_1(E)$ are polynomials in $E$ which do not share a common divisor.
At $x= \delta _{i'}$ we have the expansion
\begin{equation}
\Xi (x,E) =  \frac{d^{(i')}(E)}{(x-\delta _{i'} )^{2}} + \frac{d^{(i')}_1(E)}{\wp '(\delta _{i'}) (x-\delta _{i'} )} + (\mbox{holomorphic at } x=\delta _{i'}).
\end{equation}
By substituting this expansion into Eq.(\ref{prodDE}), we obtain the equality $d^{(i')}_1(E) =0$ upon observing the coefficient of $1/(x-\delta _{i'})^4$.
Hence we obtain the expression (\ref{Fx1}).

We express the function $\Xi (x,E)$ in descending order of powers of $E$. From Lemma \ref{lemma:Xi}, the top term is constant, hence the degrees of the coefficients in Eq.(\ref{Fx1}), other than $c_0 (E)$, are strictly less than the degree of the function $\Xi (x,E)$ in $E$.
Therefore $c_0(E) \neq 0$, $\deg _E b^{(i)}_j (E) < \deg _E c_0 (E)$ and $\deg _E d^{(i')} (E) < \deg _E c_0 (E)$ for all $i$, $i'$, $j$.
By multiplying by a constant, $c_0(E)$ is normalized to be monic.
Thus we obtain (ii).
\end{proof}

If there exists an odd-order differential operator 
$A= \left( d/dx \right)^{2g+1} +  \! $ $ \sum_{j=0}^{2g-1}\! $ $ b_j(x) \left( d/dx \right)^{2g-1-j} $ such that $[A, -d^2/dx^2+q(x)]=0$, then $q(x)$ is called the algebro-geometric finite-gap potential. Note that the equation  $[A, -d^2/dx^2+q(x)]=0$ is equivalent to the function $q(x)$ being a solution to a stationary higher-order KdV equation.

Now we construct the commuting operator, $A$, for the operator $-d^2/dx^2+v(x)$ ($v(x)$ defined in Eq.(\ref{fingapDE})) by using an expansion of the function $\Xi (x,E)$ in $E$.
Write
\begin{equation}
\Xi(x,E) = \sum_{i=0}^{g} a_{g-i}(x) E^i. \label{Xiag}
\end{equation}
It follows from Proposition \ref{prop:prodFG} that $a_0(x)=1$.
Since the function $\Xi (x,E)$ in Eq.(\ref{Xiag}) satisfies the differential equation (\ref{prodDE}), we obtain the following relations by equating the coefficients of $E^{g-j}$:
\begin{equation}
a'''_j(x)-4v(x)a'_j(x)-2v'(x)a_j(x)+4a'_{j+1}(x)=0.
\label{a'''v}
\end{equation}

\begin{thm} \label{thm:PhiA}
Assume that $l_0 , l_1, l_2 , l_3 \in \Zint _{\geq 0}$, values $\delta _1 ,\dots, \delta _M$ satisfy Eq.(\ref{eq:ds}), $\delta _j \not \equiv \omega _i$ mod $2\omega_1 \Zint \oplus 2\omega_3 \Zint$ $(0\leq i\leq 3, \; 1\leq j\leq M)$ and  $\delta _j \pm \delta _{j'} \not \equiv 0$  mod $2\omega_1 \Zint \oplus 2\omega_3 \Zint$ $(1\leq j< j' \leq M)$. Let $v(x)$ be the function defined in Eq.(\ref{fingapDE}).
Define the $(2g+1)$st-order differential operator $A$ by
\begin{equation}
A= \sum_{j=0}^{g} \left\{ a_j(x)\frac{d}{dx}-\frac{1}{2} \left( \frac{d}{dx} a_j(x) \right) \right\} \left( - \frac{d^2}{dx^2} +v(x) \right) ^{g-j}, \label{Adef}
\end{equation}
where the $a_j(x)$ are defined in Eq.(\ref{Xiag}).
Then the operator $A$ commutes with the operator $H=-d^2/dx^2  +v(x) $. In other words, the function $v(x)$ is an algebro-geometric finite-gap potential.
\end{thm}
\begin{proof}
The commutativity of the operators $A$ and $H$ follow from Eq.(\ref{a'''v}). See also \cite[Theorem 3.1]{Tak3}.
\end{proof}
Upon setting
\begin{align}
& Q(E)= \Xi (x,E)^2\left( E- v(x)\right) +\frac{1}{2}\Xi (x,E)\frac{d^2\Xi (x,E)}{dx^2}-\frac{1}{4}\left(\frac{d\Xi (x,E)}{dx} \right)^2, \label{constFG}
\end{align}
it is shown similarly to Eq.(\ref{pfconst}) that $Q(E)$ is independent of $x$. By definitions of $\Xi (x,E)$ and $Q(E)$, $Q(E)$ is a monic polynomial in $E$ of degree $2g+1$. The following proposition is proved by reviewing \cite[Proposition 3.2]{Tak3}: 
\begin{prop} \label{prop:algrel}
Let $H$ be the operator $-d^2/dx^2  +v(x)$, $A$ be the operator defined by Eq.(\ref{Adef}) and $Q(E)$ be the polynomial defined in Eq.(\ref{constFG}). Then
\begin{equation} 
A^2+Q(H)=0. \label{algrel}
\end{equation}
\end{prop}

We now relate the present work to Picard's potential. Let $q(x)$ be an elliptic function. If the differential equation $(-d^2/dx^2 +q(x))f(x)=Ef(x)$ has a meromorphic fundamental system of solutions with respect to $x$ for all values of $E$, then $q(x)$ is called a Picard potential (see \cite{GW1}). It is known that, under the condition that $q(x)$ is an elliptic function, $q(x)$ is a Picard potential if and only if $q(x)$ is an algebro-geometric potential (see \cite{GW2} and the references therein). Hence the function $v(x)$ defined in Eq.(\ref{fingapDE}) with Eq.(\ref{eq:ds}) is a Picard potential. It is possible to prove directly that $v(x)$ is a Picard potential by combining Lemma \ref{prop:locmonod} and the apparency of singularities at $x=\pm \delta _{i'}$ $(i'=1,\dots ,M)$ ensured by Eq.(\ref{eq:ds}).

\subsection{Monodromy and hyperelliptic integral} \label{sec:hypell}

We obtain an integral representation of solutions to the differential equation (\ref{fingapDE}), and express the monodromy in terms of a hyperelliptic integral.
Throughout this subsection, we assume that $l_0, l_1, l_2 ,l_3 \in \Zint _{\geq 0}$, $\delta _1 ,\dots ,\delta _M$ satisfy Eq.(\ref{eq:ds}),  $\delta _j \not \equiv \omega _i$ mod $2\omega_1 \Zint \oplus 2\omega_3 \Zint$ $(0\leq i\leq 3, \; 1\leq j\leq M)$ and  $\delta _j \pm \delta _{j'} \not \equiv 0$  mod $2\omega_1 \Zint \oplus 2\omega_3 \Zint$ $(1\leq j< j' \leq M)$.

An integral representation of solutions is obtained in Proposition \ref{prop:Linteg}. Namely, the function
\begin{equation}
\Lambda ( x,E)=\sqrt{\Xi (x,E)}\exp \int \frac{ \sqrt{-Q(E)}dx}{\Xi (x,E)},
\label{integ1FG}
\end{equation}
is a solution to the differential equation (\ref{fingapDE}).

Assume that the value $E_0$ satisfies $Q(E_0)=0$. Then it follows from Proposition \ref{prop:Q0F} that the function $\Lambda (x,E_0)$ is doubly-periodic up to signs, i.e., $\Lambda (x +2\omega _k,E_0) /\Lambda (x,E_0) \in \{\pm 1\}$ $(k=1,3)$.
In \cite[Theorem 3.7]{Tak3} the monodromy of solutions to Heun's equation for the case $l_0,l_1,l_2,l_3 \in \Zint$ is calculated in terms of a hyperelliptic integral. Similarly, we can calculate the monodromy of solutions to Eq.(\ref{fingapDE}) in terms of a hyperelliptic integral.

\begin{prop} (c.f. \cite[Theorem 3.7]{Tak3}) \label{thm:conj3} 
Assume that $E_0$ satisfies $Q(E_0)=0$. Then there exist $q_1, q_3 \in \{0,1\}$ such that $\Lambda (x+2\omega _k,E_0)=(-1)^{q_k} \Lambda (x,E_0)$ and
\begin{equation} 
\Lambda (x+2\omega _k,E)=(-1)^{q_k} \Lambda (x,E) \exp \left( -\frac{1}{2} \int_{E_0}^{E}\frac{ \int_{0+\varepsilon }^{2\omega _k+\varepsilon }\Xi (x,\tilde{E})dx}{\sqrt{-Q(\tilde{E})}} d\tilde{E}\right)
\label{analcontP}
\end{equation}
for $k=1,3$ with $\varepsilon $ denoting a constant chosen in order to avoid passing through the poles in the integration. 
\end{prop}
\begin{proof}
This proposition is proved by analogous argument to the proof of \cite[Theorem 3.7]{Tak3}.\end{proof}

We express Eq.(\ref{analcontP}) more explicitly. Since the function $\wp (x)^n $ is written as a linear combination of the functions $\left( \frac{d}{dx} \right) ^{2j} \wp (x)$ $(j=0, \dots ,n)$, the function $\Xi (x,E)$ can be expressed as
\begin{align}
\Xi (x,E)= & \: c(E)+\sum_{i=0}^3 \sum_{j=0}^{l_i-1 } a^{(i)}_j (E)\left( \frac{d}{dx} \right) ^{2j} \wp (x+\omega_i)  \label{FFx} \\
& \quad \quad \quad +\sum _{i'=1}^M d^{(i')}(E)(\wp (x+\delta _{i'}) +\wp (x- \delta _{i'})) . \nonumber
\end{align}
Set
\begin{equation}
a(E)=\sum _{i=0}^3 a^{(i)} _0 (E)+2 \sum _{i'=1}^M d^{(i')}(E).
\label{polaE}
\end{equation}
From Proposition \ref{thm:conj3} we have
\begin{equation}
\Lambda (x+2\omega _k,E)=(-1)^{q_k}\Lambda (x,E) \exp \left( -\frac{1}{2} \int_{E_0}^{E}\frac{ -2\eta _k a(\tilde{E}) +2\omega _k c(\tilde{E}) }{\sqrt{-Q(\tilde{E})}} d\tilde{E}\right) ,
\label{hypellint}
\end{equation}
for $k=1,3$, where $\eta _k =\zeta (\omega _k )$ $(k=1,3)$.
If $Q(E')\neq 0$, then the functions $\Lambda (x,E')$ and $\Lambda (-x,E')$ are a basis of the space of solutions to Eq.(\ref{fingapDE}) (see Proposition \ref{prop:indep}). Thus, if $Q(E')\neq 0$, then the monodromy matrix of solutions to Eq.(\ref{fingapDE}) on the basis $(\Lambda (x,E') , \Lambda (-x,E'))$ with respect to the cycle $x \rightarrow x+2\omega _k$ $(k=1,3)$, is diagonal and described by hyperelliptic integrals as Eq.(\ref{hypellint}).

\subsection{Bethe Ansatz and Hermite-Krichever Ansatz} \label{sec:HKA}

In this subsection we express a solution to Eq.(\ref{fingapDE}) in the form of the Bethe Ansatz and also in the form of the Hermite-Krichever Ansatz. The monodromy is described by the data of the Hermite-Krichever Ansatz (or the Bethe Ansatz).
Throughout this subsection, we will also assume that $l_0, l_1, l_2 ,l_3 \in \Zint _{\geq 0}$, $\delta _1 ,\dots ,\delta _M$ satisfy Eq.(\ref{eq:ds}), $\delta _j \not \equiv \omega _i$ mod $2\omega_1 \Zint \oplus 2\omega_3 \Zint$ $(0\leq i\leq 3, \; 1\leq j\leq M)$ and  $\delta _j \pm \delta _{j'} \not \equiv 0$  mod $2\omega_1 \Zint \oplus 2\omega_3 \Zint$ $(1\leq j< j' \leq M)$.

Set $l=2M+\sum _{i=0}^3 l_i $, $\tilde{l} _0 = 2M+l_0 $, $\tilde{l}_i =l_i$ $(i=1,2,3)$ and
\begin{equation}
\Psi _g (x)=\prod _{i'=1}^M (\wp (x) -\wp (\delta _{i'})).
\end{equation}
Assume that $Q (E')\neq 0$ and the dimension of the space of even doubly-periodic solutions to Eq.(\ref{prodDE}) is one. By Proposition \ref{prop:BA} (i), the function $\Lambda (x,E') $ in Eq.(\ref{integ1FG}) is expressed in the form of the Bethe Ansatz. Namely,
\begin{align}
& \Lambda (x,E') = \frac{C_0 \prod_{j=1}^l \sigma(x-t_j)}{\Psi _g (x) \sigma(x)^{\tilde{l}_0}\sigma_1(x)^{\tilde{l}_1}\sigma_2(x)^{\tilde{l}_2}\sigma_3(x)^{\tilde{l}_3}}\exp \left(cx \right), 
\label{eq:tilLFG}
\end{align}
for some $t_1, \dots , t_l$, $c$ and $C_0 (\neq 0)$.
It follows from Proposition \ref{prop:BA} (vi) that 
\begin{equation}
\left. \frac{d\Xi (x,E')}{dz}\right| _{z=z_j}=\frac{2\sqrt{-Q(E')}}{\wp'(t_j)},
\label{signpptjFG}
\end{equation}
where $z= \wp (x)$ and $z_j=\wp (t_j)$.

The function $\Lambda (x,E')$ is also expressed in the form of the Hermite-Krichever Ansatz. Recall that the function $\Phi _i(x,\alpha )$ $(i=0,1,2,3)$ defined in Eq.(\ref{Phii}) has periodicity described as Eq.(\ref{ddxPhiperiod}).
\begin{prop}
(i) The function $\Lambda (x,E)$ in Eq.(\ref{integ1FG}) is expressed as
\begin{align}
& \Lambda (x,E) = \frac{\exp \left( \kappa x \right) }{\Psi _g (x) } \left( \sum _{i=0}^3 \sum_{j=0}^{\tilde{l}_i-1} \tilde{b} ^{(i)}_j \left( \frac{d}{dx} \right) ^{j} \Phi _i(x, \alpha ) \right)
\label{LalphaFG}
\end{align}
for $\alpha $, $\kappa$ and $\tilde{b} ^{(i)}_j$ $(i=0,\dots ,3, \: j= 0,\dots ,\tilde{l}_i-1)$, or 
\begin{align}
& \Lambda  (x,E) = \frac{\exp \left( \bar{\kappa } x \right)}{\Psi _g (x) }  \left( \bar{c} +\sum _{i=0}^3 \sum_{j=0}^{\tilde{l}_i-2} \bar{b} ^{(i)}_j \left( \frac{d}{dx} \right) ^{j} \wp (x+\omega _i) +\sum_{i=1}^3 \bar{c}_i \frac{\wp '(x)}{\wp (x)-e_i} \right)
\label{Lalpha0FG}
\end{align}
for $\bar{\kappa }$, $\bar{c}$, $\bar{c}_i$ $(i=1,2,3)$ and $\bar{b} ^{(i)}_j$ $(i=0,\dots ,3, \: j= 0,\dots ,\tilde{l}_i-2)$.
If $\Lambda (x,E)$ is expressed as Eq.(\ref{LalphaFG}), then
\begin{align}
& \Lambda (x+2\omega _k,E) = \exp (-2\eta _k \alpha +2\omega _k \zeta (\alpha ) +2 \kappa \omega _k ) \Lambda (x,E) , \quad  (k=1,3). \label{ellintFG} 
\end{align}
(ii) There exist polynomials $P_1(E), \dots ,P_6 (E)$ such that, if $P_2(E') \neq 0$, then the function $\Lambda (x,E')$ in Eq.(\ref{integ1FG}) is written in the form of Eq.(\ref{LalphaFG}), and the values $\alpha $ and $\kappa $ are expressed as
\begin{equation}
 \wp (\alpha ) =\frac{P_1 (E')}{P_2 (E')}, \; \; \; \wp ' (\alpha ) =\frac{P_3 (E')}{P_4 (E')} \sqrt{-Q(E')} , \; \; \kappa  =\frac{P_5 (E')}{P_6 (E')} \sqrt{-Q(E')}.
\label{P1P6}
\end{equation}
If $P_2(E') = 0$, then the function $\Lambda (x,E')$ in Eq.(\ref{integ1FG}) is expressed in the form of Eq.(\ref{Lalpha0FG}).
\end{prop}
\begin{proof}
(i) follows from Theorem \ref{thm:alpha}. Note that $\Lambda (x,E) \Psi _g (x)$ is a solution to Eq.(\ref{eq:Hg}).

(ii) is proved by quite a similar argument to that of the proof of Proposition \ref{prop:P6HK}. We provide a sketch of the proof of (ii).

We assume that $Q (E')\neq 0$ and the dimension of the space of solutions to Eq.(\ref{prodDE}), which are even doubly-periodic for fixed $E'$, is one. For the case $Q(E')=0$, or the case when the dimension of the space of solutions to Eq.(\ref{prodDE}), which, for fixed $E'$, are even and doubly-periodic, is more than one, (ii) is shown by considering a continuation on parameter $E$.

It follows from Eq.(\ref{eq:tilLFG}) that
\begin{align} 
& \Lambda (x+2\omega _k ,E') = \exp \left( 2\eta _k \left( - \sum_{j=1}^l t_{j} + \sum_{i=1}^3 l_{i}\omega_{i} \right)  +2 \omega _k \left( c - \sum _{i=1}^3 l_i \eta_i \right) \right) \Lambda (x,E')   
\end{align}
for $k=1,3$. By comparing with Eq.(\ref{ellintFG}), we have
\begin{align}
& \alpha \equiv \sum_{j=1}^l t_{j} - \sum_{i=1}^3 l_{i}\omega_{i} \quad \quad  (\mbox{mod }2\omega_1 \Zint \oplus 2\omega_3 \Zint), \\
& \kappa = -\zeta \left(\sum_{j=1}^{l} t_{j} - \sum_{i=1}^3 l_i \omega_i \right) + \sum_{j=1}^{l} \zeta (t_{j})- \sum_{i=1}^3 l_i \eta_i + \delta _{l_0,0} \frac{\sqrt{-Q(E')}}{\Xi (0,E')}.
\end{align}
It follows from expressing $\wp (\sum_{j=1}^l t_{j} - \sum_{i=1}^3 l_{i}\omega_{i})$ as a combination of $\wp (t_j)$ and $\wp '(t_j)$ $(j=1,\dots ,l)$, and applying Eq.(\ref{signpptjFG}) together with the expression
\begin{equation}
\Xi (x,E) \Psi _g(x) ^2 =\frac{D \prod_{j=1}^{l}(\wp(x)-\wp(t_{j}))}{(\wp(x)-e_1)^{l_1}(\wp(x)-e_2)^{l_2}(\wp(x)-e_3)^{l_3}}
\end{equation}
for $D\neq 0$ that $\wp (\alpha )$ is expressed as a rational function in $E$. 
We can similarly obtain expressions for $\wp '(\alpha )$ and $\kappa $ in the form of Eq.(\ref{P1P6}).

The condition $P_2 (E') =0$ is equivalent to the condition $\alpha \equiv 0$ mod $2\omega_1 \Zint \oplus 2\omega_3 \Zint$. If $\alpha \equiv 0$ (resp. $\alpha \not \equiv 0$), then the function $\Lambda (x,E')$ is expressed as Eq.(\ref{LalphaFG}) (resp. Eq.(\ref{Lalpha0FG})).
Thus we obtain (ii).
\end{proof}

\subsection{Hyperellptic-ellptic reduction formulae} \label{sec:red}
We obtain hyperelliptic-elliptic reduction formulae by comparing two expressions of monodromies.
The following argument is analogous to the one in \cite[\S 3]{Tak4}.

By comparing Eq.(\ref{hypellint}) and Eq.(\ref{ellintFG}), we have 
\begin{align}
 -\eta _k \left( 2\alpha +\int_{E_0}^{E}\frac{a(\tilde{E})}{\sqrt{-Q(\tilde{E})}} d\tilde{E} \right) +\omega _k \left( 2(\zeta (\alpha )+\kappa )+ \int_{E_0}^{E}\frac{c(\tilde{E}) }{\sqrt{-Q(\tilde{E})}} d\tilde{E}\right) &\\
  =\pi \sqrt{-1} (q_k+2n_k), & \nonumber
\end{align}
for $k=1,3$ and integers $n_1$ and $n_3$.
By Legendre's relation $\eta _1 \omega _3 - \eta _3 \omega _1 =\pi\sqrt{-1}/2$, it follows that
\begin{align}
& \alpha +\frac{1}{2}\int_{E_0}^{E}\frac{a(\tilde{E})}{\sqrt{-Q(\tilde{E})}} d\tilde{E} =-(q_1+2n_1)\omega_3 +(q_3+2n_3)\omega _1, \label{alpE} \\
& \zeta (\alpha )+\kappa + \frac{1}{2}\int_{E_0}^{E}\frac{c(\tilde{E}) }{\sqrt{-Q(\tilde{E})}} d\tilde{E} =-(q_1+2n_1)\eta_3 +(q_3+2n_3)\eta _1. \label{alpEz}
\end{align}
We set $\xi =\wp (\alpha )$. By a similar argument to that of \cite[Proposition 2.4]{Tak4}, it may be proved that $\alpha \rightarrow 0$ $($mod $2\omega_1 \Zint \oplus 2\omega_3 \Zint )$ as $E \rightarrow \infty$. Combining with the relation $\int (1/\wp'(\alpha )) d\xi =\int d\alpha$, we have
\begin{equation}
\int _{\infty} ^{\xi } \frac{d \tilde{\xi }}{\sqrt{4 \tilde{\xi } ^3-g_2 \tilde{\xi } -g_3}} = \alpha = -\frac{1}{2} \int _{\infty}^{E} \frac{a(\tilde{E})}{\sqrt{-Q(\tilde{E})}}d\tilde{E}.
\label{alpint}
\end{equation}
Note that $Q(E)$ is a polynomial of degree $2g+1$, while $a(E)$ is a polynomial of degree $g$. Hence Eq.(\ref{alpint}) represents a formula which reduces a hyperelliptic integral of the first kind to an elliptic integral of the first kind. The transformation of variables is given by $\xi =P_1(E) /P_2 (E)$ for polynomials $P_1 (E)$ and $P_2 (E)$ (see Eq.(\ref{P1P6})).
Let $\alpha _0$ denote the value of $\alpha $ at $E=E_0$, where $E_0$ is the value satisfying $Q(E_0)=0$. It follows from Eq.(\ref{alpE}) that $\alpha _0 =-(q_1+2n_1)\omega_3 +(q_3+2n_3)\omega _1$ and 
\begin{align}
& \alpha - \alpha _0 +\frac{1}{2}\int_{E_0}^{E}\frac{a(\tilde{E})}{\sqrt{-Q(\tilde{E})}} d\tilde{E} =0.
\label{aaaQ}
\end{align}
If $\alpha _0 \equiv 0$ (mod $2\omega _1 \Zint \oplus 2\omega _3 \Zint$), then $\zeta (\alpha -\alpha _0 )= \zeta (\alpha ) + (q_1+2n_1)\eta_3  - (q_3+2n_3)\eta _1$.
Combining with Eqs.(\ref{alpEz}, \ref{aaaQ}), we have
\begin{align}
& \kappa = -\frac{1}{2} \int  _{E_0}^{E} \frac{c(\tilde{E})}{\sqrt{-Q(\tilde{E})}}d\tilde{E} + \zeta \left(  \frac{1}{2}\int_{E_0}^{E}\frac{a(\tilde{E})}{\sqrt{-Q(\tilde{E})}} d\tilde{E} \right) . \label{kap0}
\end{align}
If $\alpha _0 \not\equiv 0$ (mod $2\omega _1 \Zint \oplus 2\omega _3 \Zint$), then $\zeta (\alpha _0 )= -(q_1+2n_1)\eta_3 +(q_3+2n_3)\eta _1$ and 
\begin{align}
 \kappa  & = -\frac{1}{2} \int  _{E_0}^{E} \frac{c(\tilde{E})}{\sqrt{-Q(\tilde{E})}}d\tilde{E} + \int _{\wp (\alpha _0)} ^{\xi } \frac{ \tilde{\xi } d \tilde{\xi }}{\sqrt{4\tilde{\xi }^3-g_2 \tilde{\xi }-g_3}} . \label{kap123} 
\end{align}
Note that $Q(E)$ is a polynomial of degree $2g+1$, $c(E)$ is a polynomial of degree $g+1$ and $\kappa$ is expressed as $\kappa =\sqrt{-Q(E)}P_5 (E)/P_6 (E)$ for polynomials $P_5(E)$ and $P_6 (E)$ (see Eq.(\ref{P1P6})). Hence Eq.(\ref{kap123}) represents a formula which reduces a hyperelliptic integral of the second kind to an elliptic integral of the second kind, and the transformation of variables is also given by $\xi =P_1(E) /P_2 (E)$.

The following proposition describes the asymptotic behavior of $\wp (\alpha )$ and $\kappa $ as $E \rightarrow \infty$, which is proved in a similar manner to \cite[Proposition 3.2]{Tak4}.
\begin{prop} \label{prop:asym} (c.f. \cite[Proposition 3.2]{Tak4})
As $E \rightarrow \infty$, we have $\alpha \sim \frac{1}{2\sqrt{- E}} ( 4M+\sum_{i=0}^3 l_i(l_i+1)) $, $\wp (\alpha ) \sim -4E/(4M+\sum_{i=0}^3 l_i(l_i+1))^2$ and $\kappa \sim \sqrt{- E} (1-2/(4M+ \sum_{i=0}^3 l_i(l_i+1)))$.
\end{prop}

In \cite{Tak4}, following Maier \cite{Mai}, twisted Heun polynomials and theta-twisted Heun polynomials are introduced. We can extend the notions of twisted Heun polynomials and theta-twisted Heun polynomials to our potential to express the transformation of variables $\xi =P_1(E) /P_2 (E)$ and the value $\kappa =\sqrt{-Q(E)}P_5 (E)/P_6 (E)$.

\section{Examples on finite-gap potential} \label{sec:exa}

We here consider in detail several examples on finite-gap potential discussed in section \ref{sec:FGP}. The results below partially overlap with those of Smirnov \cite{Smi2}.

\subsection{The case $M=1$, $l_0=l_1=l_2=l_3=0$}

The differential equation is written as
\begin{align}
&  \left( -\frac{d^2}{dx^2} + 2(\wp (x-\delta _{1}) + \wp (x+\delta _{1}))\right) f(x) =Ef(x). \label{fingapDE10000} 
\end{align}
Set $b_1 = \wp (\delta _1)$. Then the condition that the regular singular points $x =\pm \delta _1$ of Eq.(\ref{fingapDE10000}) are apparent is given by
\begin{equation}
\prod _{i=1}^3 (b_1^2 -2 e_i b_1 -2 e_i^2 +g_2/4) =0
\end{equation}
(see Eq.(\ref{eq:f0b10})). This equation is equivalent to $\wp (2\delta _1)= e_i$ for some $i \in \{1,2,3 \}$, which is solved by $\delta _1 \equiv \omega _i/2$ mod $\omega_1 \Zint \oplus \omega_3 \Zint$.
By the shift $x \rightarrow x+\delta _1$, Eq.(\ref{fingapDE10000}) is written as 
\begin{align}
&  \left( -\frac{d^2}{dx^2} + 2(\wp (x) + \wp (x+\omega _{i}))\right) f(x) =Ef(x), \end{align}
whose potential is the Treibich-Verdier potential for the case $(l_0 ,l_1, l_2, l_3) =(1,1,0,0)$, $(1,0,1,0)$ or $(1,0,0,1)$.

We derive the functions that have appeared in section \ref{sec:FGP} for the case $(b_1^2 -2 e_i b_1 -2 e_i^2 +g_2/4) =0$ $(i \in \{1,2,3 \})$.
The functions $\Xi (x,E)$ and $Q(E)$ are given by
\begin{align}
& \Xi (x,E)=E-3e_i+\wp (x-\delta _{1}) + \wp (x+\delta _{1}), \\
& Q(E)= (E-4e_i)(E^2-2e_iE+g_2-11e_i^2) .
\end{align}
Hence the genus of the associated curve $\nu ^2=-Q(E)$ is one, and a third-order commuting operator is constructed from $\Xi (x,E)$ (see Theorem \ref{thm:PhiA}). The function $\Lambda (x, E)$ defined by Eq.(\ref{integ1FG}) is a solution to Eq.(\ref{fingapDE10000}), and the monodromy formula corresponding to Eq.(\ref{hypellint}) is given by
\begin{equation}
\Lambda (x+2\omega _k,E)=\Lambda (x,E) \exp \left( -\frac{1}{2} \int_{4e_i}^{E}\frac{ -4\eta _k +2\omega _k (\tilde{E}-3e_i) }{\sqrt{-Q(\tilde{E})}} d\tilde{E}\right) ,
\label{hypellint10000}
\end{equation}
for $k=1,3$.
The function $\Lambda (x, E)$ admits an expression in the form of the Hermite-Krichever Ansatz as
\begin{equation}
\Lambda (x,E) = \frac{\exp \left( \kappa x \right) }{\wp (x)- \wp (\delta _1) } \left( \tilde{b}_0 \Phi _0 (x, \alpha ) + \tilde{b}_1  \frac{d}{dx} \Phi _0 (x, \alpha ) \right)
\label{Lalpha10000}
\end{equation}
for generic $E$, and the values $\alpha $ and $\kappa $ satisfy
\begin{equation}
\wp (\alpha )= e_i -\frac{E^2-2e_iE+g_2-11e_i^2}{4(E-4e_i)}, \quad \kappa = \frac{1}{2}\sqrt{\frac{-(E^2-2e_iE+g_2-11e_i^2)}{(E-4e_i)}}.
\end{equation}
The monodromy is written by using the values $\alpha $ and $\kappa $ (see Eq.(\ref{ellintFG})). By comparing the two expressions of monodromy, we obtain
\begin{align}
& \int _{\infty} ^{\xi } \frac{d \tilde{\xi }}{\sqrt{4 \tilde{\xi } ^3-g_2 \tilde{\xi } -g_3}} = - \int _{\infty}^{E} \frac{d\tilde{E}}{\sqrt{-Q(\tilde{E})}}, \label{alpint10000} \\
&  \kappa = -\frac{1}{2} \int  _{E_0}^{E} \frac{\tilde{E}-3e_i}{\sqrt{-Q(\tilde{E})}}d\tilde{E} + \int _{e_i} ^{\xi } \frac{ \tilde{\xi } d \tilde{\xi }}{\sqrt{4\tilde{\xi }^3-g_2 \tilde{\xi }-g_3}} , \label{kap12310000} 
\end{align}
for the transformation 
\begin{equation}
\xi = e_i -\frac{E^2-2e_iE+g_2-11e_i^2}{4(E-4e_i)},
\end{equation}
where $E_0$ satisfies $E_0^2-2e_iE_0+g_2-11e_i^2=0$, and these formulae are related to the Landen transformation. Note that our results are compatible with the one of the Treibich-Verdier potential for the case $(l_0, l_1, l_2 ,l_3)=(1,1,0,0)$ (see \cite{Tak4}).

\subsection{The case $M=1$, $l_0=1$, $l_1=l_2=l_3=0$} \label{sec:exaM11000}

The differential equation is written as
\begin{align}
&  \left( -\frac{d^2}{dx^2} + 2(\wp (x) + \wp (x-\delta _{1}) + \wp (x+\delta _{1}))\right) f(x) =Ef(x). \label{fingapDE11000} 
\end{align}
Set $b_1 = \wp (\delta _1)$. Then the condition that the regular singular points $x =\pm \delta _1$ of Eq.(\ref{fingapDE11000}) are apparent is given by
\begin{equation}
(b_1^4-g_2 b_1^2/2-g_3 b_1 -g_2^2/48) (b_1^2-g_2/12) =0.
\end{equation}

We will first obtain the functions of the present study for the case $b_1^4-g_2 b_1^2/2-g_3 b_1 -g_2^2/48=0$.
 
The functions $\Xi (x,E)$ and $Q(E)$ are given by
\begin{align}
& \Xi (x,E)=E-6b_1+\wp (x) +\wp (x-\delta _{1}) + \wp (x+\delta _{1}), \\
& Q(E)= E^3-12b_1E^2+9(2b_1^2+g_2/4)E+126b_1^3-39g_2b_1/2-27g_3/4 ,
\end{align}
and $Q(E)$ is factorized as $Q(E)=(E-\alpha _1)(E-\alpha _2)(E-\alpha _3)$, where
$\alpha _i=3(8b_1^3+8e_ib_1^2-2(8e_i^2+g_2)b_1+3e_ig_2-12e_i^3)/(g_2-12e_i^2)$ $(i=1,2,3)$.
The genus of the associated curve $\nu ^2=-Q(E)$ is one, and a third-order commuting operator is constructed from $\Xi (x,E)$ (see Theorem \ref{thm:PhiA}). The function $\Lambda (x, E)$ defined by Eq.(\ref{integ1FG}) is a solution to Eq.(\ref{fingapDE11000}), and the monodromy formula corresponding to Eq.(\ref{hypellint}) may be expressed in the form
\begin{equation}
\Lambda (x+2\omega _k,E)=- \Lambda (x,E) \exp \left( -\frac{1}{2} \int_{\alpha _2}^{E}\frac{ -6\eta _k +2\omega _k (\tilde{E}-6b_1) }{\sqrt{-Q(\tilde{E})}} d\tilde{E}\right) ,
\label{hypellint11000}
\end{equation}
for $k=1,3$.
The function $\Lambda (x, E)$ admits an expression in the form of the Hermite-Krichever Ansatz as
\begin{equation}
\Lambda (x,E) = \frac{\exp \left( \kappa x \right) }{\wp (x)- \wp (\delta _1) } \left( \tilde{b}_0 \Phi _0 (x, \alpha ) + \tilde{b}_1  \frac{d}{dx} \Phi _0 (x, \alpha ) + \tilde{b}_2  \left(\frac{d}{dx} \right)^2 \Phi _0 (x, \alpha ) \right)
\label{Lalpha11000}
\end{equation}
for generic $E$, and the values $\alpha $ and $\kappa $ satisfy
\begin{align}
& \wp (\alpha )= e_i - \frac{(E-\alpha _i)(E-9b_1+9e_i/2+\alpha _i/2))^2}{9(E-7b_1)^2} \quad \quad (i=1,2,3)\\
& \quad \quad = -\frac{E^3-18b_1E^2+99b_1^2E-126b_1^3+(-93b_1/2+9E/2)g_2-27g_3}{9(E-7b_1)^2}, \nonumber \\
& \kappa = \frac{2}{3(E-7b_1)}\sqrt{-Q(E)}.
\end{align}
The monodromy is written by using the values $\alpha $ and $\kappa $ (see Eq.(\ref{ellintFG})). By comparing the two expressions of monodromy, we obtain
\begin{align}
& \int _{\infty} ^{\xi } \frac{d \tilde{\xi }}{\sqrt{4 \tilde{\xi } ^3-g_2 \tilde{\xi } -g_3}} = -\frac{3}{2} \int _{\infty}^{E} \frac{d\tilde{E}}{\sqrt{-Q(\tilde{E})}}, \label{alpint11000} \\
&  \kappa = -\frac{1}{2} \int  _{e_i}^{E} \frac{\tilde{E}-6b_1}{\sqrt{-Q(\tilde{E})}}d\tilde{E} + \int _{\alpha _i} ^{\xi } \frac{ \tilde{\xi } d \tilde{\xi }}{\sqrt{4\tilde{\xi }^3-g_2 \tilde{\xi }-g_3}} , \quad \quad (i=1,2,3), \label{kap12311000} 
\end{align}
for the transformation 
\begin{equation}
\xi =  -\frac{E^3-18b_1E^2+99b_1^2E-126b_1^3+(-93b_1/2+9E/2)g_2-27g_3}{9(E-7b_1)^2}.
\end{equation}

We now consider the case $b_1^2-g_2/12=0$.
The functions $\Xi (x,E)$ and $Q(E)$ are given by
\begin{align}
& \Xi (x,E)=(E^2-3g_2/2)+\frac{(2b_1g_2+3g_3)E+g_2^2+18b_1g_3}{2b_1g_2+3g_3}\wp (x) \\
& \quad \quad \quad +\frac{2(2b_1g_2+3g_3)E-(g_2^2+18b_1g_3)}{2(2b_1g_2+3g_3)}(\wp (x-\delta _{1}) + \wp (x+\delta _{1})), \nonumber \\
& Q(E)= (E-6b_1)(E+6b_1)(E-3e_1)(E-3e_2)(E-3e_3).
\end{align}
Hence the genus of the associated curve $\nu ^2=-Q(E)$ is two, and a fifth-order commuting operator is constructed from $\Xi (x,E)$ (see Theorem \ref{thm:PhiA}). The function $\Lambda (x, E)$ defined by Eq.(\ref{integ1FG}) is a solution to Eq.(\ref{fingapDE11000}), and the monodromy formula corresponding to Eq.(\ref{hypellint}) is given by
\begin{equation}
\Lambda (x+2\omega _k,E)=\Lambda (x,E) \exp \left( -\frac{1}{2} \int_{6b_1}^{E}\frac{ -6\eta _k \tilde{E} +2\omega _k (\tilde{E}^2-3g_2/2) }{\sqrt{-Q(\tilde{E})}} d\tilde{E}\right) ,
\label{hypellint11000-2}
\end{equation}
for $k=1,3$.
The function $\Lambda (x, E)$ admits an expression in the form of the Hermite-Krichever Ansatz as Eq.(\ref{Lalpha11000}) for generic $E$, and the values $\alpha $ and $\kappa $ satisfy
\begin{equation}
 \wp (\alpha )=-\frac{E^3-27g_3}{9(E^2-3g_2)} , \quad \kappa = \frac{2}{3}\sqrt{\frac{-(E^3-9g_2E/4-27g_3/4)}{(E^2-3g_2)}}.
\end{equation}
The monodromy may be written upon using the values $\alpha $ and $\kappa $ (see Eq.(\ref{ellintFG})). By comparing the two expressions of monodromy, we obtain 
\begin{align}
& \int _{\infty} ^{\xi } \frac{d \tilde{\xi }}{\sqrt{4 \tilde{\xi } ^3-g_2 \tilde{\xi } -g_3}} = - \frac{3}{2} \int _{\infty}^{E} \frac{\tilde{E}d\tilde{E}}{\sqrt{-Q(\tilde{E})}}, \label{alpint11000^2} \\
&  \kappa = -\frac{1}{2} \int  _{e_i}^{E} \frac{\tilde{E}^2-3g_2/2}{\sqrt{-Q(\tilde{E})}}d\tilde{E} + \int _{3e_i} ^{\xi } \frac{ \tilde{\xi } d \tilde{\xi }}{\sqrt{4\tilde{\xi }^3-g_2 \tilde{\xi }-g_3}} , \label{kap12311000-2} 
\end{align}
for the transformation 
\begin{equation}
\xi = -\frac{E^3-27g_3}{9(E^2-3g_2)}.
\end{equation}
These formulae reduce hyperelliptic integrals of genus two to elliptic integrals. Note that our results are similar to that of the Treibich-Verdier potential for the case $(l_0, l_1, l_2 ,l_3)=(2,0,0,0)$ (see \cite{Tak4}), and these two potentials may be related by an isospectral deformation (see \cite{Smi2}).

\subsection{The case $M=1$, $l_0=2$, $l_1=l_2=l_3=0$}

The differential equation is written as
\begin{align}
&  \left( -\frac{d^2}{dx^2} + 6\wp (x) + 2(\wp (x-\delta _{1}) + \wp (x+\delta _{1}))\right) f(x) =Ef(x). \label{fingapDE12000} 
\end{align}
Set $b_1 = \wp (\delta _1)$. Then the condition that the regular singular points $x =\pm \delta _1$ of Eq.(\ref{fingapDE12000}) are apparent is given by
\begin{equation}
b_1^6-53g_2 b_1^4/100 -17g_3 b_1^3/25+19g_2^2 b_1^2/400+11g_2 g_3 b_1/100+ g_2^3/1600+g_3^2/25=0.
\end{equation}

Upon setting 
\begin{align}
& H^{(0)}(E)=E^2-\frac{160b_1^3-24g_2b_1-16g_3}{12b_1^2-g_2}E-\frac{800b_1^4-124g_2b_1^2-32g_3b_1+3g_2^2}{12b_1^2-g_2}, \\
& H^{(i)}(E)=E-\{800b_1^5+800e_ib_1^4+(-384g_2+320e_i^2)b_1^3+(60g_2-656e_i^2)e_ib_1^2 \\
& +(30g_2^2+144e_i^2g_2-512e_i^4)b_1+(-23g_2^2+12g_2e_i^2+256e_i^4)e_i\} /\{(12e_i^2-g_2)(12b_1^2-g_2)\}, \nonumber
\end{align}
for $i=1,2,3$,  the functions $\Xi (x,E)$ and $Q(E)$ are given by
\begin{align}
& \Xi (x,E)=E^2-10b_1E-\frac{3}{4}\frac{800b_1^4-124b_1^2g_2-32b_1g_3+3g_2^2}{12b_1^2-g_2} +9\wp (x)^2 \\
& \quad +3(E-4b_1)\wp (x) +\left( E-\frac{-72b_1^3+14b_1g_2+12g_3}{12b_1^2-g_2} \right) (\wp (x-\delta _{1})  + \wp (x+\delta _{1}) ) , \nonumber \\
& Q(E)= H^{(0)}(E)H^{(1)}(E)H^{(2)}(E)H^{(3)}(E).
\end{align}
The genus of the associated curve $\nu ^2=-Q(E)$ is two, and a fifth-order commuting operator is constructed from $\Xi (x,E)$ (see Theorem \ref{thm:PhiA}). The function $\Lambda (x, E)$ defined by Eq.(\ref{integ1FG}) is a solution to Eq.(\ref{fingapDE12000}), and the monodromy formula corresponding to Eq.(\ref{hypellint}) is written as
\begin{align}
& \Lambda (x+2\omega _k,E)=- \Lambda (x,E) \cdot \label{hypellint12000} \\
& \quad \cdot \exp \left( \int_{\alpha _2}^{E}\frac{ \eta _k\left( 5\tilde{E}-\frac{8(2b_1g_2+3g_3)}{12b_1^2-g_2}\right) -\omega _k \left(\tilde{E}^2-10b_1\tilde{E}-\frac{3(200b_1^4-34b_1^2g_2-8b_1g_3+g_2^2)}{12b_1^2-g_2}\right) }{\sqrt{-Q(\tilde{E})}} d\tilde{E}\right) , \nonumber 
\end{align}
for $k=1,3$, where $\alpha _i$ satisfies $H^{(i)}(\alpha _i)=0$ $(i=1,2,3)$.
The function $\Lambda (x, E)$ admits an expression in the form of the Hermite-Krichever Ansatz as
\begin{equation}
\Lambda (x,E) = \frac{\exp \left( \kappa x \right) }{\wp (x)- \wp (b_1) } \left( \sum_{j=0}^{3} \tilde{b}_j \left( \frac{d}{dx} \right) ^{j} \Phi _0 (x, \alpha ) \right)
\label{Lalpha12000}
\end{equation}
for generic $E$, and the values $\alpha $ and $\kappa $ satisfy
\begin{align}
& \wp (\alpha )= e_i - \frac{H^{(i)}(E)Ht^{(i)}(E)^2}{25H^{(0)}(E)Ht^{(0)}(E)^2}, \quad \quad (i=1,2,3), \\
& \kappa = \frac{4H\theta (E)}{5Ht^{(0)}(E)}\sqrt{-\frac{H^{(1)}(E)H^{(2)}(E)H^{(3)}(E)}{H^{(0)}(E)}},
\end{align}
where
\begin{align}
& H\theta (E)= E+\frac{-420b_1^3+59g_2b_1+36g_3}{2(12b_1^2-g_2)}, \quad Ht^{(0)}(E)= E+\frac{-900b_1^3+83g_2b_1+12g_3}{5(12b_1^2-g_2)}, \\
& Ht^{(i)}(E)=E^2+E \{400b_1^5+400e_ib_1^4+(108g_2-3440e_i^2)b_1^3+(-120g_2+1472e_i^2)e_ib_1^2 \\
&  +(-18g_2^2+468e_i^2g_2-256e_i^4)b_1+(13g_2^2-336g_2e_i^2+704e_i^4)e_i\}/\{(12e_i^2-g_2)(12b_1^2-g_2)\} \nonumber \\
& +\{1600e_ib_1^5+(-2800g_2+35200e_i^2)b_1^4+(6352g_2-84800e_i^2)e_ib_1^3 \nonumber \\
& +(-24g_2^2+2688e_i^2g_2-28672e_i^4)b_1^2 +(-1152g_2^2+16240g_2e_i^2-17920e_i^4)e_ib_1 \nonumber \\
& +41g_2^3-328e_i^2g_2^2-4096e_i^4g_2+19712e_i^6\}/\{2(12e_i^2-g_2)(12b_1^2-g_2)\}, \quad (i=1,2,3). \nonumber
\end{align}
The monodromy may be written using the values $\alpha $ and $\kappa $ (see Eq.(\ref{ellintFG})). By comparing the two expressions of monodromy, we obtain
\begin{align}
& \int _{\infty} ^{\xi } \frac{d \tilde{\xi }}{\sqrt{4 \tilde{\xi } ^3-g_2 \tilde{\xi } -g_3}} = -\frac{1}{2} \int _{\infty}^{E} \frac{5\tilde{E}-\frac{8(2b_1g_2+3g_3)}{12b_1^2-g_2}}{\sqrt{-Q(\tilde{E})}}d\tilde{E}, \label{alpint12000} \\
&  \kappa = -\frac{1}{2} \int  _{e_i}^{E} \frac{\tilde{E}^2-10b_1\tilde{E}-\frac{3(200b_1^4-34b_1^2g_2-8b_1g_3+g_2^2)}{12b_1^2-g_2}}{\sqrt{-Q(\tilde{E})}}d\tilde{E} + \int _{\alpha _i} ^{\xi } \frac{ \tilde{\xi } d \tilde{\xi }}{\sqrt{4\tilde{\xi }^3-g_2 \tilde{\xi }-g_3}} , \label{kap12312000} 
\end{align}
for the transformation 
\begin{equation}
\xi =  e_i - \frac{H^{(i)}(E)Ht^{(i)}(E)^2}{25H^{(0)}(E)Ht^{(0)}(E)^2} , \quad \quad (i=1,2,3).
\end{equation}
These formulae reduce hyperelliptic integrals of genus two to elliptic integrals, which may not be reduced to the case of the Treibich-Verdier potential.

\subsection{The case $M=2$, $l_0=l_1=l_2=l_3=0$} \label{sec:exaM20000}

The differential equation is written as
\begin{align}
&  \left( -\frac{d^2}{dx^2} + 2\sum_{i'=1}^2 (\wp (x-\delta _{i'}) + \wp (x+\delta _{i'}))\right) f(x) =Ef(x). \label{fingapDE20000} 
\end{align}

The equation for the apparency of singularity is written as
\begin{align}
& \wp '(2\delta _1 )+\wp '(\delta _1 +\delta _2) +\wp '(\delta _1 -\delta _2) =0, \quad \wp '(2\delta _2 )+\wp '(\delta _2 +\delta _1) +\wp '(\delta _2 -\delta _1) =0. \label{20000d}
\end{align}
We set $\alpha = \delta _1 +\delta _2$ and $\beta =\delta _1 -\delta _2$. Then Eq.(\ref{20000d}) are equivalent to 
\begin{align}
& \wp '(\alpha +\beta  )+\wp '(\alpha - \beta ) +2\wp '(\alpha ) =0, \quad \wp '(\alpha +\beta  )+\wp '(\beta - \alpha ) +2\wp '(\beta ) =0. \label{20000ab}
\end{align}
From the relation
\begin{equation}
\wp '(x +y )+\wp '(x -y ) =-\frac{\wp '' (y)\wp ' (x)}{(\wp (x)-\wp (y))^2}-2\frac{\wp ' (y)^2\wp ' (x)}{(\wp (x)-\wp (y))^3},
\end{equation}
we have that Eq.(\ref{20000ab}) is equivalent to
\begin{align}
& (A(\alpha , \beta )+ 2)\wp '(\alpha ) =0, \quad (A(\beta , \alpha )+ 2)\wp '(\beta ) =0, \label{eq:20000A} \\
& A(x,y)= -\frac{\wp '' (y)}{(\wp (x)-\wp (y))^2}-2\frac{\wp ' (y)^2}{(\wp (x)-\wp (y))^3}.
\end{align}

Thus solutions to Eq.(\ref{eq:20000A}) are divided into four cases; the case $\wp '(\alpha ) =\wp '(\beta ) =0$, the case $\wp '(\alpha ) =A(\beta , \alpha )+ 2=0$, the case $\wp '(\beta ) =A(\alpha , \beta )+ 2=0$, and the case $A(\alpha , \beta )+ 2=A(\beta , \alpha )+ 2=0$.

We first consider the case $\wp '(\alpha ) =\wp '(\beta ) =0$. Since $\wp '(x)=0$ is equivalent to $x \equiv \omega _1, \omega _2, \omega _3$ mod $2\omega_1 \Zint \oplus 2\omega_3 \Zint$, we have $\alpha , \beta \equiv \omega _1, \omega _2, \omega _3$ mod $2\omega_1 \Zint \oplus 2\omega_3 \Zint$. By considering the condition $\delta _1 \not \equiv 0 \not \equiv \delta _2$, we have
$(\delta _1, \delta _2)\equiv \pm ((\omega _i +\omega _j) /2, (\omega _i - \omega _j)/2 )$ mod $2\omega_1 \Zint \oplus 2\omega_3 \Zint$ for $1\leq i\neq j \leq 3$.
For this case, we have $\wp (2\delta _1)=\wp (2\delta _2)$ and
\begin{align}
\Xi(x,E)= & E-3(\wp (2\delta _1)+ \wp (\delta _1 +\delta _2)+ \wp (\delta _1 -\delta _2)) \\
& +\wp(x + \delta _1)+\wp(x - \delta _1)+\wp(x + \delta _2)+\wp(x - \delta _2). \nonumber
\end{align}
The degree of the polynomial $Q(E)$ is three, and the genus of the associated curve $\nu ^2=-Q(E)$ is one.

Secondly, we consider the case $\wp '(\alpha ) =A(\beta , \alpha )+ 2=0$. It follows from $\wp '(\alpha ) =0$ that $\alpha \equiv \omega _i$ mod $2\omega_1 \Zint \oplus 2\omega_3 \Zint$ for some $i\in \{1,2,3\}$. Then we have $A(\beta , \omega _i )+ 2=0$, which may be written as
$ 6e_i^2 -g_2 /2=2(e_i -\wp (\beta ))^2$. The solutions of this equation are given by $2\beta \equiv \omega _i $  mod $2\omega_1 \Zint \oplus 2\omega_3 \Zint$. Hence $(\delta _1 , \delta _2) \equiv ( \omega _j/2 \pm \omega _i/4, \omega _j /2 \mp \omega _i/4)$, $-( \omega _j/2 \pm \omega _i/4, \omega _j /2 \mp \omega _i/4)$  mod $2\omega_1 \Zint \oplus 2\omega_3 \Zint$ for $1\leq i\leq 3$, $0\leq j \leq 3$.
For this case, we have $\wp (2\delta _1)=\wp (2\delta _2)$ and
\begin{align}
\Xi(x,E)= & E-3(\wp (2\delta _1)+ \wp (\delta _1 +\delta _2)+ \wp (\delta _1 -\delta _2)) \\
& +\wp(x + \delta _1)+\wp(x - \delta _1)+\wp(x + \delta _2)+\wp(x - \delta _2). \nonumber
\end{align}
The degree of the polynomial $Q(E)$ is three, and the genus of the associated curve $\nu ^2=-Q(E)$ is one.
The case  $\wp '(\beta ) =A(\alpha , \beta )+ 2=0$ can be treated similarly.

We now consider the case $A(\alpha , \beta )+ 2=A(\beta , \alpha )+ 2=0$. It follows from a direct derivation that $\wp (\alpha ) +\wp (\beta )= 2e_i$ and $\wp (\alpha )\wp (\beta )=2e_i^2-g_2/12$ for some $i \in \{1,2,3\}$. The function $\Xi (x,E)$ is given by
\begin{equation}
\Xi (x,E) = c(E) + \sum _{j=1}^2 d^{(j)} (E) (\wp(x + \delta _j)+\wp(x - \delta _j)),
\end{equation}
where 
\begin{align}
& d^{(1)} (E) = E-2e_i +3\wp (2\delta _2), \quad d^{(2)} (E) = E-2e_i +3\wp (2\delta _1), \\
& \textstyle c(E)= (d^{(1)} (E)+d^{(2)} (E))(\frac{E}{2}-3e_i)-\frac{3}{2}(d^{(1)} (E) \wp (2\delta _1)+d^{(2)} (E) \wp (2\delta _2)).
\end{align} 
The degree of the polynomial $Q(E)$ is five, and the genus of the associated curve $\nu ^2=-Q(E)$ is two.

\section{Concluding remarks} \label{sec:rmk}

We have shown in sections \ref{sec:HK} and \ref{sec:P6} that solutions of the linear differential equation that produces the sixth Painlev\'e equation have integral representations and that they are expressed in the form of the Hermite-Krichever Ansatz. Furthermore we got a procedure for obtaining solutions of the sixth Painlev\'e equation (see Eq.(\ref{eq:P6ellip})) for the cases $\kappa _0 , \kappa _1 , \kappa _t, \kappa _{\infty} \in \Zint +1/2$ by fixing the monodromy, and we presented explicit solutions for the cases $(\kappa _0 , \kappa _1 , \kappa _t, \kappa _{\infty} )=(1/2, 1/2, 1/2, 1/2)$ and $(1/2, 1/2, 1/2, 3/2)$.

By B\"acklund transformation of the sixth Painlev\'e equation (see \cite{TOS} etc.), Hitchin's solution (i.e., solutions for the case $(\kappa _0 , \kappa _1 , \kappa _t, \kappa _{\infty} )=(1/2, 1/2, 1/2, 1/2)$) is transformed to the solutions for the case $(\kappa _0 , \kappa _1 , \kappa _t, \kappa _{\infty} ) \in O_1 \cup O_2$, where
\begin{align}
& O_1= \left\{ (\kappa _0 , \kappa _1 , \kappa _t, \kappa _{\infty} ) | 
\kappa _0 , \kappa _1 , \kappa _t, \kappa _{\infty} \in \Zint +1/2 \right\}, \\
& O_2 = \left \{(\kappa _0 , \kappa _1 , \kappa _t, \kappa _{\infty} ) \left| 
\begin{array}{ll}
\kappa _0 , \kappa _1 , \kappa _t, \kappa _{\infty} \in \Zint \\
\kappa _0 + \kappa _1 + \kappa _t + \kappa _{\infty}  \in 2 \Zint 
\end{array}
\right. \right\}. 
\end{align}
Note that solutions for the case $(\kappa _0 , \kappa _1 , \kappa _t, \kappa _{\infty} )=(0, 0, 0, 0)( \in O_2)$ are already known and are called Picard's solution.

For the case $(\kappa _0 , \kappa _1 , \kappa _t, \kappa _{\infty} ) \in O_1$, solutions of the linear differential equation are investigated by our method, and solutions of the sixth Painlev\'e equation follow from them.
On the other hand, for the case $(\kappa _0 , \kappa _1 , \kappa _t, \kappa _{\infty} ) \in O_2$, we cannot obtain results on integral representation and the Hermite-Krichever Ansatz by our method, although solutions of the sixth Painlev\'e equation are obtained in principle by B\"acklund transformation from the case $(\kappa _0 , \kappa _1 , \kappa _t, \kappa _{\infty} ) =(1/2, 1/2, 1/2, 1/2) $.
Note that the condition $(\kappa _0 , \kappa _1 , \kappa _t, \kappa _{\infty} ) \in O_2$ corresponds to the condition $l_0, \dots ,l_3 \in \Zint +1/2$, $l_0+ l_1 +l_2 + l_3 \in 2 \Zint $ (see Eq.(\ref{eq:kili})).

Now we propose a problem to investigate solutions and their monodromy of the linear differential equation (Eq.(\ref{Hgkr1}) with the condition (\ref{pgkrap})) for the cases $l_0, \dots ,l_3 \in \Zint +1/2$, $l_0+ l_1 +l_2 + l_3 \in 2 \Zint $. In partiuclar, how can we investigate solutions and their monodromy of the linear differential equation for the case $\kappa _0 = \kappa _1 = \kappa _t= \kappa _{\infty} =0$ (i.e. $l_0=l_1=l_2=l_3 =-1/2$)?\footnote{Note added in 2008: Results on this subject were obtained by Takemura K., Integral representation of solutions to Fuchsian system and Heun's equation, {\it J. Math. Anal. Appl.} {\bf 342} (2008), 52--69.}

In section \ref{sec:FGP}, we produced new examples of finite-gap potential and invetsigated properties of them.
More properties should be clarified in near furute.
For instance, it is not immediate to calculate the genus of the associated curve.
As is seen from examples in sections \ref{sec:exaM11000} and \ref{sec:exaM20000}, the genus depends on the solution of the equations which determine the position of apparent singularities (i.e. Eq.(\ref{eq:ds})).
Related results were obtained by Treibich \cite{Tre} for the case $M=1$, and they are to be simplified and generalized. 

To find finite-gap potential, we considered only the case $r_{i'}=2$ and $s_{i'}=0$ (see Eq.(\ref{Inopotent})) in section \ref{sec:FGP}.
We propose a problem for a study of finite-gap potential for the cases $r_{i'} \neq 2$ for some $ i'$.

\appendix
\section {Elliptic functions} \label{sect:append}
This appendix presents the definitions of and the formulas for the elliptic functions.

The Weierstrass $\wp$-function, the Weierstrass sigma-function and the Weierstrass zeta-function with periods $(2\omega_1, 2\omega_3)$ are defined as follows:
\begin{align}
& \wp (z)= \frac{1}{z^2}+  \sum_{(m,n)\in \Zint \times \Zint \setminus \{ (0,0)\} } \left( \frac{1}{(z-2m\omega_1 -2n\omega_3)^2}-\frac{1}{(2m\omega_1 +2n\omega_3)^2}\right),  \\
& \sigma (z)=z\prod_{(m,n)\in \Zint \times \Zint \setminus \{(0,0)\} } \left(1-\frac{z}{2m\omega_1 +2n\omega_3}\right) \nonumber \\
& \; \; \; \; \; \; \; \; \; \; \; \; \; \; \cdot \exp\left(\frac{z}{2m\omega_1 +2n\omega_3}+\frac{z^2}{2(2m\omega_1 +2n\omega_3)^2}\right), \nonumber \\
& \zeta(z)=\frac{\sigma'(z)}{\sigma (z)}. \nonumber
\end{align}
Setting $\omega_2=-\omega_1-\omega_3$ and 
\begin{align}
& e_i=\wp(\omega_i), \; \; \; \eta_i=\zeta(\omega_i), \; \; \; \; (i=1,2,3)
\end{align}
yields the relations
\begin{align}
& e_1+e_2+e_3=\eta_1+\eta_2+\eta_3=0, \; \; \; \label{eq:Leg} \\
& \eta _1 \omega _3- \eta _3 \omega _1 = \eta _3 \omega _2- \eta _2 \omega _3 = \eta _2 \omega _1- \eta _1 \omega _2 = \pi\sqrt{-1} /2, \nonumber \\
& \wp(z)=-\zeta'(z), \; \; \; (\wp'(z))^2=4(\wp(z)-e_1)(\wp(z)-e_2)(\wp(z)-e_3). \nonumber 
\end{align}
The periodicity of functions $\wp(z)$, $\zeta (z)$ and $\sigma (z)$ are as follows:
\begin{align}
& \wp(z+2\omega_i)=\wp(z), \; \; \; \zeta(z+2\omega_i)=\zeta(z)+2\eta_i ,\; \; \; \; (i=1,2,3), \label{periods} \\
& \sigma (z+2\omega _i) = - \sigma (z) \exp (2\eta _i (z + \omega _i)), \; \; \; \frac{\sigma (z+t+2\omega _i )}{\sigma (z+2\omega _i)}= \exp(2\eta _i t) \frac{\sigma (z+t)}{\sigma (z)} .\nonumber
\end{align}
The constants $g_2$ and $g_3$ are defined by
\begin{equation}
g_2=-4(e_1e_2+e_2e_3+e_3e_1), \; \; \; g_3=4e_1e_2e_3.
\end{equation}
The co-sigma functions $\sigma_i(z)$ $(i=1,2,3)$ and co-$\wp$ functions $\wp_i(z)$ $(i=1,2,3)$ are defined by
\begin{align}
& \sigma_i(z)=\exp (-\eta_i z)\frac{\sigma(z+\omega_i)}{\sigma(\omega _i)}, \; \; \; \wp_i(z) = \frac{\sigma_i(z)}{\sigma(z)}, \label{eq:sigmai}
\end{align}
and satisfy
\begin{align}
& \wp_i(z) ^2 =\wp(z)-e_i, \quad \quad \quad (i,i' =1,2,3) \label{rel:sigmai} \\
& \wp _i (z+2\omega _{i'}) = \exp (2(\eta _{i'} \omega _i -\eta _{i} \omega _{i'}) ) \wp _i (z) = (-1)^{\delta _{i,i'}} \wp _i (z). \nonumber
\end{align}

{\bf Acknowledgments.}
The author would like to thank Prof. Hidetaka Sakai for fruitful discussions, and Prof. F. Gesztesy, Prof. A. V. Kitaev and the referee for valuable comments.
He is partially supported by a Grant-in-Aid for Scientific Research (No. 15740108) from the Japan Society for the Promotion of Science.

\end{document}